\documentclass[12pt]{amsart}
\usepackage{amscd,amsfonts,amsmath,amsxtra,amssymb}

\usepackage{psfig}
\usepackage[all]{xy}

\def\supp{\mathop{\hbox{$\sup$}}\limits}
\def\sepsec{\vskip2cm}
\def\sepsub{\vskip7mm}
\def\sepsubsub{\vskip4mm}
\def\sepprop{\vskip4mm}

\newtheorem{sub}{}[section]
\newtheorem{subsub}{}[sub]

\newcommand{\M}{{\mathbb M}}

\newcommand{\C}{{\mathbb C}}
\newcommand{\Q}{{\mathbb Q}}
\renewcommand{\P}{{\mathbb P}}

\newcommand{\E}{{\mathbb E}}

\renewcommand{\M}{{\mathbb M}}

\newcommand{\ka}{{\mathcal A}}
\newcommand{\kb}{{\mathcal B}}

\newcommand{\ke}{{\mathcal E}}
\newcommand{\kf}{{\mathcal F}}

\newcommand{\ki}{{\mathcal I}}

\newcommand{\kk}{{\mathcal K}}
\newcommand{\kl}{{\mathcal L}}

\newcommand{\ko}{{\mathcal O}}
\newcommand{\kp}{{\mathcal P}}

\newcommand{\kx}{{\mathcal X}}

\newcommand{\kz}{{\mathcal Z}}

\def\lra{\longrightarrow}
\def\sigg{\mathop{\hbox{$\displaystyle\sum$}}\limits}
\def\som{\mathop{\hbox{$\displaystyle\bigoplus$}}\limits}
\def\ot{\otimes}

\def\bw{{\bf W}}
\def\bg{{\bf G}}
\def\ws{W^{s}(G,\Lambda)} 
\def\wss{W^{ss}(G,\Lambda)} 
\def\bws{{\bf W}^{s}({\bf G},\tilde{\Lambda})} 
\def\bwss{{\bf W}^{ss}({\bf G},\tilde{\Lambda})} 
\def\bvs{{\bf V}^{s}({\bf G}_L\times G_R, \bar{\Lambda})} 
\def\bvss{{\bf V}^{ss}({\bf G}_L\times G_R, \bar{\Lambda})}
\def\and{\quad\text{and}\quad}
\def\txtem#1{{\quad\text{\em #1}\quad}}
\def\txte#1{{\quad\text{#1}\quad}}
\def\dsp{\displaystyle}
\def\ov#1{\overline{#1}}
\def\codim{\mathop{\rm codim}\nolimits}
\def\coker{\mathop{\rm coker}\nolimits}
\def\Hom{\mathop{\rm Hom}\nolimits}
\def\GL{\mathop{\rm GL}\nolimits}

\def\Ext{\mathop{\rm Ext}\nolimits}
\def\EExt{\mathop{\mathcal Ext}\nolimits}

\def\Aut{\mathop{\rm Aut}\nolimits}
\def\End{\mathop{\rm End}\nolimits}
\def\imm{\mathop{\rm Im}\nolimits}
\def\lra{\longrightarrow}
\def\som{\mathop{\hbox{$\displaystyle\bigoplus$}}\limits}
\def\sigg{\mathop{\hbox{$\displaystyle\sum$}}\limits}

\def\supp{\mathop{\hbox{$\sup$}}\limits}

\def\hfl#1#2{\smash{\mathop{\ \hbox to 12mm{\rightarrowfill}}
\limits^{\scripstyle#1}_{\scripstyle#2} \ }}
\def\hflb#1#2{\smash{\mathop{\hbox to 12mm{\leftarrowfill}}
\limits^{\scripstyle#1}_{\scriptstyle#2}}}

\def\pline#1{<\hskip-3.5pt#1\hskip-3.5pt>}

\def\m#1{{\hbox{$#1$}}}

\def\ot{\otimes}
\def\og{\leavevmode\raise.3ex\hbox{$\scriptscriptstyle\langle\!\langle$}}
\def\fg{\leavevmode\raise.3ex\hbox{$\scriptscriptstyle\,\rangle\!\rangle$}}

\def\nsp{\lbrace 0\rbrace}

\def\txtem#1{{\quad\text{\em #1}\quad}}
\def\txte#1{{\quad\text{#1}\quad}}
\def\dsp{\displaystyle}
\def\flinc{\ar@{^{(}->}}
\def\fleq{\ar@{=}}
\def\flon{\ar@{->>}}
\def\fmaps{\ar@{|-{>}}}

\begin{document}

\title[{\tiny Quotients by non-reductive groups}]{Moduli spaces of decomposable morphisms of 
sheaves and quotients by non-reductive groups}

\author[{J.M.Dr\'{e}zet}]{Jean--Marc Dr\'{e}zet}
\address{
Institut de Math\'{e}matique\\
UMR 7586 du CNRS, Aile 45-55, 5$^e$ \'{e}tage\\ 
2, place Jussieu\newline   
F-75251 Paris Cedex 05, France}
\email{drezet@math.jussieu.fr}
\author[{G.Trautmann}]{G\"unther Trautmann}
\address{
Universit\"at Kaiserslautern\\
Fachbereich Mathematik\\
Erwin-Schr\"odinger-Stra{\ss}e\newline
D-67663 Kaiserslautern}
\email{trm@mathematik.uni-kl.de}

\begin{abstract}
We extend the methods of geometric invariant theory to actions of non--reductive
groups in the case of homomorphisms between decomposable sheaves 
whose automorphism groups are non--reductive.
Given a linearization of the natural action of the group $\Aut(E)\times\Aut(F)$
on Hom(E,F), a homomorphism is called stable if its orbit with respect to the 
unipotent radical is contained in the stable locus with respect to the natural 
reductive subgroup of the automorphism group. We encounter effective numerical 
conditions for a linearization such that the corresponding open set of 
semi-stable homomorphisms admits a good and projective quotient in the sense
of geometric invariant theory, and that this quotient is in addition a geometric
quotient on the set of stable homomorphisms. 
\end{abstract}

\maketitle
%\today
\tableofcontents
\newpage
\section{Introduction}

Let $X$ be a projective algebraic variety over the field of complex
numbers. Given two coherent sheaves $\ke, \kf$, on $X$ the algebraic
group $G=\Aut(\ke)\times \Aut(\kf)$ acts naturally on the affine space
$W=\Hom(\ke, \kf)$ by $(g,h). w=h\circ w\circ g^{-1}$.
If two morphisms are in the same $G$--orbit then they have isomorphic
cokernels and kernels. Therefore it is natural to ask for good quotients of
such actions in the sense of geometric invariant theory. 
\vskip2mm

\begin{sub}{\sc\small\small Morphisms of type $(r,s)$.}\rm\\
In general $\ke$ and $\kf$ will be decomposable such that $G$ is not
reductive. More specifically let $\ke$ and $\kf$ be direct sums 
\[
\ke=\underset{1\leq i\leq r}{\oplus}M_i\otimes \ke_i\quad\text{ and }\quad 
\kf = \underset{1\leq l\leq s}{\oplus} N_l\otimes\kf_l ,
\]
where $M_i$ and $N_l$ are finite dimensional vector spaces and $\ke_i$,
$\kf_l$ are simple sheaves, i.e. their only endomorphisms are the homotheties,
and such that $\Hom(\ke_i, \ke_j)=0$ for $i>j$ and $\Hom(\kf_l, \kf_m)=0$ for
$l>m$. In this case we call homomorphisms $\ke \to \kf$ of type $(r,s)$. Then
the groups $Aut(\ke)$ and $\Aut(\kf)$ can be viewed as groups of
matrices of the following type. The group $\Aut(\ke)$, say, is the group of
matrices
\[
\left (
\begin{array}{ccccc}
g_1 & 0 & \cdots & 0\\
u_{21} & g_2 & & \vdots\\
\vdots & \ddots & \ddots & 0\\
u_{r1} & \cdots & u_{r, r-1} & g_r
\end{array}
\right)
\]
where $g_i\in \GL(M_i)$ and $u_{ji} \in \Hom(M_i, M_j\otimes \Hom(\ke_i,
\ke_j))$.

In the literature on moduli of vector bundles and coherent sheaves many
quotients of spaces $\P \Hom(\ke, \kf)$ of type $(1,1)$ by the reductive group
$\Aut(\ke)\times \Aut(\kf)$ have been investigated, see for example \cite{dr1},
\cite{ell_P_str}, 
\cite{elling}, \cite{karpov}, \cite{miro-trm}. 
The moduli spaces described in this way are the simplest ones, and this allows
to test in these cases some conjectures that are expected to be true on
more general moduli spaces of sheaves (cf. \cite{dr1b}, \cite{tjotta}).
We think that the moduli spaces of morphisms of type $(r,s)$ will be as  
useful to treat other less simple moduli problems of sheaves. In fact, if one
wants to use the spaces $\Hom(\ke,\kf)$ as parameter spaces for moduli spaces
of sheaves, which are as close as possible to the moduli spaces, the higher
types $(r,s)$ are unavoidable.

The homomorphisms in a Beilinson complex of a bundle on projective $n$--space,
for example, have in general arbitrary type $(r,s)$ depending on the
dimensions of the cohomology spaces of the bundle. In several papers, see \cite{miro}, \cite{oko} for example,
semi--stable sheaves or ideal sheaves of subschemes of projective spaces, are
represented as quotients of injective morphisms of type $(r,s)$, and one
should expect that the moduli spaces of such sheaves are isomorphic to a good
quotient of an open subset of the corresponding space of homomorphisms. In
some cases of type $(2,1)$ this has been verified for semi--stable sheaves on
$\P_2$ in \cite{dr2}.

In case of type $(r,s)$ there are
good and projective quotients if one restricts the action to the reductive
subgroup
\[
G_{red}=\prod\ \GL(M_i)\times\prod\GL(N_l).
\]
This has been shown recently by A. King in \cite{king}.
The quotient problem for $\Hom(\ke,\kf)$ of type $(r,s)$ with respect to the
full group $\Aut(\ke)\times \Aut(\kf)$ is however the generic one and indispensable.

Unfortunately the by now standard geometric invariant theory (GIT) doesn't
provide a direct answer for these quotient problems in case $\Aut(\ke)\times
\Aut(\kf)$ is not reductive. There are several papers dealing with the action
of an arbitrary algebraic group like \cite{faunt}, \cite{faunt2}, \cite{dix},
\cite{dixray} and older ones, but their results are insufficient for the above
problem. The conditions of \cite{faunt} are close to what we need, but they
don't allow a concrete description of the set of semi--stable points in our
case and they don't guarantee good or projective quotients, see remark
\ref{Fauntleroy}.
\end{sub}
\vskip2mm

\begin{sub}{\sc\small\small The main idea}\rm\\
Our procedure is very close to standard GIT and we finally reduce the
problem of the quotient to the one of a reductive group action. 
We introduce polarizations $\Lambda\in\Q^{r+s}$ of tuples of rational numbers
for the action of $G$ on the affine space in analogy
to the ones of A. King in \cite{king}, which are refinements of the
polarizations by ample line bundles on the projective space $\P W$, and then 
introduce open sets
$\ws\subset \wss$ of stable and semi--stable points depending on $\Lambda$ and
study the quotient problem for these open subsets. There are chambers in
$\Q^{r+s}$ such that the polarizations in one chamber define the same open set,
in accordance with the chamber structure in Neron--Severi spaces of
polarizations in the reductive case, see f.e. \cite{dolg} I. Dolgachev - Y. Hu
and \cite{thadd} M. Thaddeus. However, in contrast to the reductive case, good
quotients $\wss//G$ don't exist for all polarizations, see \ref{prem_ex}. As a main
achievement we are providing numerical conditions on the polarizations,
depending on the dimensions of the spaces $M_i$ and $N_l$, under which such
quotients exist. The main step for that is to embed the group actions $G\times
W\to W$ into an action $\bg\times\bw\to\bw$ of a reductive group $\bg$ and to 
compare the open sets $\wss$ and $\bwss$, where $\widetilde{\Lambda}$ is a
polarization for the $\bg$--action associated to $\Lambda$.
\end{sub}
\vskip2mm

\begin{sub}{\sc\small\small Construction of quotients by non reductive groups.}
\rm\\
To be more precise, a polarization $\Lambda$ is a tuple $(\lambda_1,
 \ldots,\lambda_r, \mu_1, \ldots, \mu_s)$ of positive rational numbers, called
 weights of the factors $M_i\otimes \ke_i$ and $N_l\otimes \kf_l$
 respectively, which satisfy
$\sum\lambda_i m_i=\sum \mu_l n_l=1$, where $m_i, n_l$ denote the dimensions of 
the spaces of the same name. We use then the numerical
criterion of A. King, \cite{king}, as definition for semi--stability with
respect to the reductive group $G_{red}$. An element $w\in W$ is $(G_{red},
\Lambda)$--stable if for any proper choice of subspaces $M'_i\subset M_i,\quad
N'_l\subset N_l$ such that 
$w$ maps $\oplus(M'_i\ot \ke_i)$ into $\oplus(N'_l\ot \kf_l)$,
we have $\sum\lambda_i m'_i<\sum \mu_l n'_l$, or
semi--stable if equality is allowed. Let $W^s(G_{red},\Lambda)\subset
W^{ss}(G_{red}, \Lambda)$ denote the set of stable and semi--stable points so
defined. If $H\subset G$ is the unipotent radical of $G$, which is generated by
the homomorphisms $\ke_i\to\ke_j$ and $\kf_l\to\kf_m$ for $i<j$ and $l<m$, we
say that $w$ is $(G, \Lambda)$--\hbox{(semi--)}stable if $h.w$ is $(G_{red},
\Lambda)$--(semi--)stable for any $h\in H$, see \ref{defstab}. We thus have open
subsets $\wss\subset W^{ss}(G_{red},\Lambda)$ and $\ws\subset W^{ss}(G_{red},
\Lambda)$. 
\bigskip

{\em The main result of our paper is that there are sufficient
numerical and effective bounds for the polarizations $\Lambda$ such that $\wss$
admits a good and even projective quotient $\wss//G$ and that in addition $\ws$
admits a geometric quotient, which is smooth and quasi--projective, see
proposition \ref{goodqu} and the results \ref{sone}, \ref{equality}, and
section 8.}

The definitions of good and geometric quotients are recalled in
\ref{constquot}.
By using correspondences between spaces of morphisms, called {\em mutations},
it is possible to deduce from our results other polarizations such that there
exists a good projective quotient (see \cite{dr4}, \cite{other}). 
%\newpage

All this is achieved by embedding the action $G\times W\to W$ into
an action ${\mathbf G}\times {\mathbf W}\to {\mathbf W}$ of a reductive group
and then imposing conditions for the equality $\wss=W\cap \bwss$, where
$\widetilde{\Lambda}$ is the associated polarization. The quotient is then the
quotient of the saturated subvariety ${\mathbf G}\wss\subset \bwss$. The
quotient will be projective if $\overline {{\mathbf G}.W}\smallsetminus
{\mathbf G}.W$ doesn't meet $\bwss$. Also for this, numerical conditions can be 
found in section 8.

The idea of embedding the non--reductive action $G\times W\to W$ into the
action ${\mathbf G} \times {\mathbf W}\to {\mathbf W}$ is simply to replace the
$\ke_i$ by $\ke_1$ using the evaluation maps $\Hom(\ke_1, \ke_i)\otimes
\ke_1\to \ke_i$. It is explained in \ref{motiv} and \ref{Remmotiv} that this
is the outcome when we start to replace the sheaves $\ke_i$ step by step and
similarly for the sheaves $\kf_l$. Since we have to deal everywhere with the
dimensions of the vector spaces $\Hom(\ke_i, \ke_j)$ and $\Hom(\kf_l, \kf_m)$
which form the components of the unipotent group $H$, we have translated the
whole setup into an abstract multilinear setting and related actions by
technical reasons. This gives more general results although we have only
applications in the theory of sheaves. The reader should always keep in mind
the motivation in \ref{motiv}.

The results obtained in the simplest case (morphisms of type (2,1) or (1,2)) 
are stated in \ref{result-2-1}. They are characteristic for the general case in
which only the conditions are more complicated.
\end{sub}
\vskip2mm

\begin{sub}{\sc\small\small Remark on finite generatedness}\rm\\ 
One would expect that the quotients of $W$ could be obtained by first forming
the quotient $W/H$ with respect to the unipotent radical $H$ and then in a
second step a quotient of $W/H$ by $G/H\cong G_{red}$. However, the actions of
unipotent groups behave generally very badly, \cite{gr-pf}, and we are not able
to prove that the algebra $\C[W/H]$ is finitely generated. This would
be an essential step in a direct construction of the quotient. Of course, the
main difficulty also in this paper arises from the presence of the group
$H$. The counterexample of M.Nagata, \cite{naga}, also shows that the finite
generatedness depends on the dimensions of the problem. So from a
philosophical point of view we are determining bounds for the dimensions
involved under which we can expect local affine $G$-invariant coordinate rings 
which are finitely generated, and thus to obtain good quasi--projective 
quotients, even so the bounds might not be the best. 
The simple examples \ref{prem_ex}, \ref{remgeoqu} show that a good quotient
$\wss//G$ might not exist if the conditions are not fulfilled.
\end{sub}
\vskip2mm

\begin{sub}{\sc\small\small Morphisms of type (2,1)}\label{result-2-1}\rm\\ 
In this case the homomorphisms of sheaves are of the type
\[m_1\ke_1\oplus m_2\ke_2\lra n_1\kf_1 ,\]
where we use the notation $m\ke$ for $\C^m\otimes\ke$.
For this type a polarization is given by a pair $(\lambda_1,\lambda_2)$ of positive
rational numbers such that \ $\lambda_1m_1+\lambda_2m_2=1$. It is determined
by the rational number \ $t=m_2\lambda_2$ \ which lies in $[0,1]$. Writing 
$W^{ss}(t)$ for $W^{ss}$ and $W^s(t)$ for $W^s$ for the moment, our results 
depend on constants $c(k)$ defined as follows : Let
$$\tau : \Hom(\ke_1,\kf_1)^*\otimes \Hom(\ke_1,\ke_2)
\longrightarrow \Hom(\ke_2,\kf_1)^*$$
be the linear map induced by the composition map \
$\Hom(\ke_2,\kf_1)\otimes\Hom(\ke_1,\ke_2)\to\Hom(\ke_1,\kf_1)$ , and
$$\tau_k = \tau\otimes I_{\C^k} : 
\Hom(\ke_1,\kf_1)^*\otimes (\Hom(\ke_1,\ke_2)\otimes\C^k) 
\longrightarrow \Hom(\ke_2,\kf_1)^*\otimes\C^k .$$
Let $\kk$ be the set of proper linear subspaces
\ $K\subset \Hom(\ke_1,\ke_2)\otimes\C^{k}$ \
such that for every proper linear subspace $F\subset\C^{k}$, $K$ is not
contained in  $\Hom(\ke_1,\ke_2)\otimes F$. Let
$$c(k) = \supp_{K\in{\kk}}
(\frac{\codim(\tau_k(\Hom(\ke_1,\kf_1)^*\otimes K)}
{\codim(K)}).$$

\sepprop

\begin{subsub}\label{theo_main}{\bf Theorem: }
There exists a good projective quotient  \ $W^{ss}(t)//G$ \ and a geometric
quotient \ $W^s(t)/G$ \ if
\[t \ > \ \frac{m_2\dim(\Hom(\ke_1,\ke_2))}{\dim(\Hom(\ke_1,\ke_2))+m_1} 
\ \ \ \ {\rm and} \ \ \ \
t \ > \ \dim(\Hom(\ke_1,\ke_2)).c(m_2)\frac{m_2}{n_1}. \]
\end{subsub}

\sepprop

In the case of morphisms $m_1\ko(-2)\oplus m_2\ko(-1)\lra n_1\ko $
on projective spaces the constants have been computed in \cite{other} and we 
obtain the more explicit result :

\sepprop

\begin{subsub}\label{theo_main2}{\bf Theorem: } Let $n\geq 2$ be an integer.
There exists a good projective quotient  \ $W^{ss}(t)//G$ \ and a geometric
quotient \ $W^s(t)/G$ \ in the case of morphisms \hfil\break $m_1\ko(-2)\oplus
m_2\ko(-1)\lra n_1\ko $ \ on the projective space $\P_n$ if
\[ \begin{array}{lcl}
t & > & \dsp\frac{(n+1)m_2}{(n+1)m_2+m_1} , \\[3ex] 
t & > & \dsp\frac{(n+1)m_2^2(m_2-1)}{2n_1(m_2(n+1)-1)}\txte{if}2\leq 
m_2\leq n+1,\\[2.5ex]
t & > & \dsp\frac{(n+1)^2m_2}{2(n+2)n_1}\txte{if}m_2>n+1.
\end{array}
\]  
\end{subsub}
\end{sub}

\vskip5mm

\begin{sub}{\sc\small\small Construction of fine moduli spaces of torsion free
sheaves}\label{int_fine}\rm\\ 
In section \ref{Fine_Mod} we construct smooth projective fine moduli spaces of
torsion free coherent sheaves on $\P_n$ using morphisms
\[(*)  \ \ \ \ \
\ko(-2)\ot\C^2\lra\ko(-1)\oplus(\ko\ot\C^k) ,\]
(for \ $(n+1)(n+2)/2<k<(n+1)^2$ ). More precisely we prove that for all
polarizations, semi-stable morphisms are injective outside a closed subvariety
of codimension $\geq 2$, hence their cokernels are torsion free sheaves. A
generic morphism is injective and its cokernel is locally free. In this case we can construct 
\[q \ = \ \frac{(n+1)(n+2)}{2} - \left[\frac{n+k+1}{2}\right]\]
distinct smooth projective moduli spaces $\M_1,\cdots,\M_q$ of such morphisms, of dimension
\hfil\break
\m{2(n-1)+k((n+1)^2-k)}. Moreover, all the $\M_i$ are birational to each other.
For $1\leq i\leq q$, we construct a coherent sheaf
$\ke_i$ on $\M_i\times\P_n$, flat over $\M_i$, such that for every closed point
$z\in\M_i$, $\ke_{iz}$ is isomorphic to the cokernel of the morphism $(*)$
corresponding to $z$. We prove that $\M_i$ is a fine moduli space of
torsion free sheaves with universal sheaf $\ke_i$. In particular, this means
that for every closed point $z\in \M_i$, the Koda\"\i ra-Spencer map
\[T_z\M_i\lra\Ext^1(\ke_{iz},\ke_{iz})\]
is bijective, and for any two distinct closed points $z_1,z_2\in\M_i$, the
sheaves $\ke_{iz_1}$, $\ke_{iz_2}$ are not isomorphic. 
\end{sub}
\vskip2mm

\begin{sub}{\sc\small\small Open problems}\label{quest}\rm\\ 
Even in the simplest case of morphisms of type (2,1) we do not know what all
the polarizations are for which a good quotient $W^{ss}//G$
exists. More generally it would be interesting to find all the saturated open
subsets $U$ of $W$ such that a good quotient (quasiprojective or not) 
\m{U//G} exists, or
all the open subsets $U$ such that a geometric quotient $U/G$ exists.
The corresponding problem for reductive groups has been studied in \cite{mumf},
1.12, 1.13, and in \cite{bbs1}, \cite{bbs2}.
\end{sub}
\vskip2mm

\begin{sub}{\sc\small\small Organization of the paper}\label{Org}\rm\\ 
In section 2 we describe our problem in
terms of multilinear algebra. 

In section 3 we recall results of A. King,
\cite{king}. The reductive group actions considered in this paper, the action
of $G_{red}$ on $W$ and that of ${\mathbf G}$ on ${\mathbf W}$, are 
particular cases of \cite{king}. We also discuss the relation of
$\Lambda$--(semi--)stability in $W$ with that in the projective space
$\P W$. But we cannot work solely on the projective niveau, because the
embedding $W\subset {\mathbf W}$ is not linear. 

After defining
$G$--(semi--)stability for the non--reductive group in section 4 we describe
the embedding in section 5 and introduce the associated polarizations. 

Section 6 contains the step of constructing the quotient $\wss//G$ using the
GIT--quotient $\bwss//{\mathbf G}$ of A. King. 

Sections 7 and 8 are the hard parts of the
paper. Here the conditions of the weights which define good polarizations are
derived. It seems that the constants appearing in these estimates had not been
considered before. 

In section 9 we are investigating a few examples in order to test
the strength of the bounds. Here we restrict ourselves to small type $(2,1),
(2,2), (3,1)$ in order to avoid long computations of the constants which give
the bounds for the polarizations. What we discover in varying the
polarizations are flips between the moduli spaces, as one has to expect from
the general results on the variation of linearizations of group actions, cf. 
\cite{reid}, \cite{dolg}, \cite{thadd}. In example \ref{ex1} we have a very
simple effect of a flip, but in example \ref{ex2gen} the chambers of the
polarizations look already very complicated.

In section 10 we define new fine moduli spaces of torsion free sheaves using
our moduli spaces of morphisms.

\vskip2mm

{\sc\small Acknowledgement}. The work on this paper was supported by DFG. 
The first author wishes to thank the University of
Kaiserslautern, where the work was started, for its hospitality.
\end{sub}

\sepsec
%\newpage
\section{The moduli problem for decomposable homomorphisms}

Let $\ke=\oplus\; \ke_i\otimes M_i$ and $\kf=\oplus\; \kf_l\otimes N_l$ 
be semi--simple sheaves as in the
introduction. In order to describe the action of $G=\Aut(\ke)\times \Aut(\kf)$
on $W=\Hom(\ke, \kf)$ in greater detail we use the abbreviations
\[
\begin{array}{lcl}
H_{li} & = & \Hom(\ke_i, \kf_l) \\
A_{ji} & = & \Hom(\ke_i, \ke_j) \\
B_{ml} & = & \Hom(\kf_l, \kf_m)\ ,
\end{array}
\]
such that we are given the natural pairings
\[
\begin{array}{lclclp{1cm}ll}
H_{lj} & \otimes & A_{ji} & \to & H_{li} & & \text{ for } & i\leq j\\
A_{kj}  & \otimes & A_{ji} & \to & A_{ki} & & \text{ for } & i\leq j\leq k\\
B_{ml}& \otimes & H_{li} & \to & H_{mi} & & \text{ for } & l\leq m\\
B_{nm} & \otimes & B_{ml} & \to & B_{nl} & & \text{ for } & l\leq
m\leq n. 
\end{array}
\]
The group $G$ consists now of pairs $(g, h)$ of matrices
\[
g=\left(  \begin{array}{cccc}
g_1 & 0 & \cdots & 0\\
u_{21} & g_2 & & \vdots\\
\vdots & \ddots & \ddots & 0\\
u_{r1} & \cdots & u_{r, r-1} & g_r
  \end{array}
\right)\quad \text{ and }\quad
h=\left(
  \begin{array}{cccc}
h_1 & 0 & \cdots & 0\\
v_{21} & h_2 & & \vdots\\
\vdots & \ddots & \ddots & 0\\
v_{s1} & \cdots & v_{s,s-1} & h_s
  \end{array}
\right)
\]
with diagonal elements $g_i\in \GL(M_i),\  h_l\in \GL(N_l)$ and $u_{ji}\in 
\Hom(M_i, M_j\otimes A_{ji}),\\
v_{ml}\in \Hom(N_l, N_m \otimes B_{ml})$.

Similarly a homomorphism $w\in \Hom(\ke,\kf)$ is represented by a matrix
$w=(\varphi_{li})$ of homomorphisms $\varphi_{li}\in \Hom(M_i, N_l\otimes
H_{li})=\Hom(H_{li}^\ast \otimes M_i, N_l)$. Using the natural pairings, the 
left action $(g, h).w=h wg^{-1}$ of $G$ on $W$ is described by the matrix product
\[
\left(
  \begin{array}{cccc}
h_1 & 0 & \cdots & 0\\
v_{21} & h_2 & & \vdots\\
\vdots & \ddots & \ddots & 0\\
v_{s1} & \cdots & v_{s,s-1} & h_s
  \end{array}
\right)\ \circ\ 
\left(
\begin{array}{ccc}
\varphi_{11} & \cdots &\varphi_{1r}\\[1.5ex]
\vdots & & \vdots\\[1.5ex]
\varphi_{s1} & \cdots & \varphi_{sr}\end{array}\right)\ \circ\ 
\left(
  \begin{array}{cccc}
g_1 & 0 & \cdots & 0\\
u_{21} & g_2 & & \vdots\\
\vdots & \ddots & \ddots & 0\\
u_{r1} & \cdots & u_{r, r-1} & g_r
  \end{array}
\right)^{-1}\raisebox{-5ex}{,}
\]
where the compositions $v_{ml} \circ \varphi_{li}$ and $\varphi_{lj}\circ
u_{ji}$ are compositions as sheaf homomorphisms but can also be
interpreted as compositions induced by the pairings of the vector spaces
above. Thus the group $G$, the space $W$ and the action are already determined
by the vector spaces $A_{ji}, B_{ml}, H_{li}$ and the pairings between
them. Therefore, in the following we define $G, W$ and the actions 
$G\times W\to W$ by abstractly given vector spaces and pairings. The resulting 
statements can then be applied to systems of sheaves by specifying the spaces
as spaces of homomorphisms as above.
\sepsub

\begin{sub}\label{abstr}{\sc\small The abstract setting}\rm

Let $r, s$ be positive integers and let for $1\leq i\leq j\leq r,\ \ 1\leq l\leq
m \leq s$ finite dimensional vector spaces $A_{ji}, B_{ml}, H_{li}$ be given,
where we assume that $A_{ii}=\C$ and $B_{ll}=\C$. Moreover we suppose that we
are given linear maps, called {\it compositions},
\[
\begin{array}{lclll}
H_{lj}\otimes A_{ji}&\to & H_{li} & \quad\text{ for } & 1\leq i\leq j\leq r,\quad
1\leq l\leq s\\
A_{kj}\otimes A_{ji}&\to & A_{ki} & \quad\text{ for } & 1\leq i\leq j\leq k\leq
r\\
B_{ml} \otimes H_{li}&\to & H_{mi} & \quad\text{ for } & 1\leq i\leq r,\quad 1\leq l\leq m\leq s\\
B_{nm}\otimes B_{ml}&\to & B_{nl} & \quad\text{ for } & 1\leq l\leq m\leq n\leq s.
\end{array}
\]

We assume that all these maps and the induced maps
\[
H_{li}^\ast \otimes A_{ji} \to H_{lj}^\ast \quad\text{ and }\quad H_{mi}^\ast \otimes
B_{ml}\to H_{li}^\ast
\]
are {\bf surjective}. This is the case when all the spaces are spaces of 
sheaf homomorphisms as above for which the sheaves $\ke_i$ and $\kf_l$
are line bundles on a projective space or each of them is a bundle
$\Omega^p(p)$.

We may and do assume that these pairings are the identities if $i=j$, $l=m$
etc.\ . Finally, we suppose that these maps verify the natural associative properties of 
compositions. This means that the induced diagrams 

%\[
%\begin{CD}
%A_{kj}\otimes A_{ji}\otimes A_{ih}@>>> A_{ki}\otimes A_{ih}\\
%@VVV @VVV\\
%A_{kj}\otimes A_{jh}@>>> A_{kh}
%\end{CD}\hskip1cm 
%\begin{CD}
%B_{on}\otimes B_{nm}\otimes B_{ml} @>>> B_{om}\otimes B_{nl}\\
%@VVV @VVV\\
%B_{on}\otimes B_{nl}  @>>> B_{ol} 
%\end{CD}
%\]
\[\xymatrix{
A_{kj}\ot A_{ji}\ot A_{ih}\ar[r]\ar[d] & A_{ki}\otimes A_{ih}\ar[d]\\
A_{kj}\ot A_{jh}\ar[r] & A_{kh}
} \ \ \ \ \
\xymatrix{
B_{on}\ot B_{nm}\ot B_{ml}\ar[r]\ar[d] & B_{om}\otimes B_{nl}\ar[d]\\
B_{on}\ot B_{nl}\ar[r] & B_{ol} 
}\]

%for $1\leq h\leq i\leq j\leq k\leq r$ and $1\leq l\leq m\leq n\leq o\leq s$ as well as 
%\[
%\begin{CD}
%H_{lk}\otimes A_{kj}\otimes A_{ji}@>>> H_{lj}\otimes A_{ji}\\
%@VVV @VVV\\
%H_{lk}\otimes A_{ki}@>>> H_{li}
%\end{CD}\hskip1cm
%\begin{CD}
%B_{nm}\otimes B_{ml}\otimes H_{li}@>>> B_{nl}\otimes H_{li}\\
%@VVV @VVV\\
%B_{nm}\otimes H_{mi}@>>> H_{ni}
%\end{CD}
%\]
\[\xymatrix{
H_{lk}\ot A_{kj}\ot A_{ji}\ar[r]\ar[d] & H_{lj}\ot A_{ji}\ar[d]\\
H_{lk}\ot A_{ki}\ar[r] &  H_{li}
}
\ \ \ \ \ 
\xymatrix{
B_{nm}\ot B_{ml}\ot H_{li}\ar[r]\ar[d] & B_{nl}\ot H_{li}\ar[d]\\
B_{nm}\ot H_{mi}\ar[r] & H_{ni}
}
\]

%for $1\leq i\leq j\leq k\leq r$ and $1\leq l\leq m\leq n\leq s$ and

%\[
%\begin{CD}
%B_{ml}\otimes H_{lj}\otimes A_{ji}@>>> H_{mj}\otimes A_{ji}\\
%@VVV @VVV\\
%B_{ml}\otimes H_{li}@>>> H_{mi}
%\end{CD}
%\]
\[\xymatrix{
B_{ml}\ot H_{lj}\ot A_{ji}\ar[r]\ar[d] & H_{mj}\ot A_{ji}\ar[d]\\
B_{ml}\ot H_{li}\ar[r] & H_{mi}
}
\]

are commutative for all possible combinations of indices.

In our setup we also let finite dimensional vector spaces $M_i$ for 
$1\leq i\leq r$ and $N_l$ for $1\leq l\leq s$ be given
and we consider finally the vector space
\[
W=\underset{i,l}\oplus\ \Hom(M_i, N_l\otimes 
H_{li})=\underset{i,l}\oplus\ \Hom(H_{li}^\ast \otimes M_i, N_l)
\]
where summation is over $1\leq i\leq r$ and $1\leq l\leq s$. This is the space 
of homomorphisms in the abstract setting. The group $G$ and its action on $W$ 
are now also given in the abstract setting as follows.
\end{sub}
\newpage
%\sepsub

\begin{sub}\label{group}{\sc\small The group $G$}\rm

We define $G$ as a product $G_L\times G_R$ of two groups where the left group
$G_L$ replaces $\Aut(\ke)$ and the right group $G_R$ replaces $\Aut(\kf)$ in our
motivation. Let $G_L$ be the set of matrices
\[
\left (
\begin{array}{ccccc}
g_1 & 0 & \cdots & 0\\
u_{21} & g_2 & & \vdots\\
\vdots & \ddots & \ddots & 0\\
u_{r1} & \cdots & u_{r, r-1} & g_r
\end{array}
\right)
\]

with $g_i \in \GL(M_i)$ and $u_{ji}\in \Hom(M_i, M_j\otimes
A_{ji})=\Hom(A^\ast_{ji}\otimes M_i, M_j)$. The
group law in $G_L$ is now defined as matrix multiplication where we define the
compositions $u_{kj} \ast u_{ji}$ naturally according to the given pairings as
the composition
\[
M_i\xrightarrow{u_{ji}} M_j\otimes A_{ji}\xrightarrow{u_{kj}\otimes id}
M_k\otimes A_{kj}\otimes A_{ji}\xrightarrow{id\otimes comp} M_k\otimes A_{ki}.
\]

Explicitly, if $g$ has the entries $g_i, u_{ji}$ and $g'$ has the entries $g'_i, u'_{ji}$
then the product
\[
g''= g'g
\]
in $G_L$ is defined as the matrix with the entries $g''_i=g'_i\circ g_i$ in
the diagonal and

\[
u''_{ki}= u'_{ki}\circ g_i+\sum\limits_{i<j<k} u'_{kj}\ast u_{ji}+(g'_k\otimes
id)\circ u_{ki}
\]
for $1\leq i<k\leq r$. The verification that this defines a group structure on
$G_L$ is now straightforward. 

As a set $G_L$ is the product of all the $\GL(M_i)$ and all $\Hom(M_i, M_j
\otimes A_{ji})$ for $i<j$ and thus has the structure of an affine
variety. Since multiplication is composed by a system of bilinear maps it is a
morphism of affine varieties. Hence $G_L$ is naturally endowed with the
structure of an algebraic group. The group $G_R$ is now defined in the same
way by replacing the spaces $M_i$ and $A_{ji}$ by $N_l$ and $B_{ml}$. Finally
$G=G_L\times G_R$ is defined as an algebraic group.
\end{sub}\sepsub

\begin{sub}\label{action}{\sc\small The action of $G$ and $W$}\rm

We will define a left action of $G_R$ and a right action of $G_L$ on $W$ such
that the action of $G$ on $W$ can be defined by $(g,h).w=h.w.g^{-1}$. Both
actions are defined as matrix products as described above in the case of sheaf
homomorphisms using the abstract compositions as in the definition of the
group law.

If $w$ has the entries $\varphi_{li}\in \Hom(H^\ast_{li}\otimes M_i, N_l)$ and
  $g\in G_L$ has the entries $g_i$ and $u_{ij}$ then $w.g$ is defined as the
  matrix product

\[
\left(
  \begin{array}{ccc}
\varphi'_{11} & \cdots & \varphi'_{1r}\\[1.6ex]
\vdots        &        & \vdots\\[1.6ex]
\varphi'_{s1} & \cdots & \varphi_{sr}
  \end{array}
\right)= \left(
\begin{array}{ccc}
\varphi_{11} & \cdots & \varphi_{1r}\\[1.6ex]
\vdots        &        & \vdots\\[1.6ex]
\varphi_{s1} & \cdots & \varphi_{sr}
\end{array}
\right) \ 
\left(
  \begin{array}{cccc}
g_1    & 0      & \cdots & 0\\
u_{21} & g_2    & &\vdots\\
\vdots & \ddots & \ddots & 0\\
u_{r1} & \cdots & u_{r,r-1} & g_r
  \end{array}\right)
\]

with 
\[
\varphi'_{li}= \varphi_{li}\circ g_i+\sum\limits_{i<j}\varphi_{lj}\ast
u_{ji}\quad (\text{if } i=r \text{ the last sum is } 0),
\]

where $\varphi_{lj}\ast u_{ji}$ is the composition
\[
M_i\to M_j\otimes A_{ji}\to N_l\otimes H_{lj}\otimes A_{ji}\to N_l\otimes
H_{li}
\]
or dually the composition
\[
H_{li}^\ast\otimes M_i\to H_{lj}^\ast \otimes A_{ji}^\ast \otimes M_i\to
H_{lj}^\ast \otimes M_j\to N_l.
\]

The left action of $G_R$ is defined in the same way. In the next two sections
we give an analysis of stability and semi-stability for the action of $G$ and
its natural reductive subgroup $G_{red}$. In the reductive case this is due
to A. King.
\end{sub}\sepsub

\begin{sub}\label{subgrp}{\sc\small Canonical subgroups of $G$}\rm

We let $H_L\subset G_L$ and $H_R\subset G_R$ be the maximal normal
unipotent subgroups of $G_L$ and $G_R$ defined by the condition that all
$g_i=id_{M_i}$ and all $h_l=id_{N_l}$. Then $H=H_L\times H_R$ is a maximal
normal unipotent subgroup of $G$. Similarly we consider the reductive
subgroups $G_{L, red}$ and $G_{R, red}$ of $G_L$ and $G_R$ defined by the
conditions $u_{ji}=0$ and $v_{ml}=0$ for all indices. Then 
$G_{red}=G_{L,red}\times G_{R, red}$ is a reductive subgroup of $G$ and it is 
easy to see
that $G/H\cong G_{red}$. The restricted action of $G_{red}$ is much simpler
and reduces to the natural actions of $\GL(M_i)$ on $M_i$ and $\GL(N_l)$ on 
$N_l$
\end{sub}
%\sepsec
\newpage
\section{Actions of reductive groups}
\begin{sub}\label{ResKing}{\sc\small Results of A. King}\rm

Let $Q$ be a finite set,\ $\Gamma\subset Q\times Q$\ a subset
such that the union of the images of the two projections of $\Gamma$ is $Q$. 
For each $\alpha\in Q$, let $m_\alpha$ be a
positive integer, $M_\alpha$ a vector space of dimension $m_\alpha$
and for each $(\alpha,\beta)\in\Gamma$, let 
$V_{\alpha\beta}$ be a finite dimensional nonzero vector space. Let
$$W_0 = \underset{(\alpha,\beta)\in\Gamma}{\oplus}
\Hom(M_\alpha\ot V_{\alpha\beta},M_\beta).$$
On $W_0$ we have the following action of the reductive group
$$G_0 = \prod_{\alpha\in Q}\GL(M_\alpha)$$
arising naturally in this situation. If \ $(f_{\beta\alpha})\in W_0$ \ 
and \ $(g_{\alpha})\in G_0$, then
$$(g_{\alpha}).(f_{\beta\alpha}) = (g_\beta\circ
f_{\beta\alpha}\circ(g_\alpha\ot id)^{-1}).$$
Let $(e_\alpha)_{\alpha\in Q}$ be a sequence of integers such that
$$\sum_{\alpha\in Q}e_\alpha m_\alpha = 0 .$$
To this sequence is associated the character $\chi$ of $G_0$ defined by
$$\chi(g) = \prod_{\alpha\in Q}
\det(g_\alpha)^{-e_\alpha}.$$
This character is trivial on the canonical subgroup of $G_0$ isomorphic to $\C^*$ (for every \ $\lambda\in \C^*$, the element $(g_\alpha)$
of $G_0$ corresponding to $\lambda$ is such that \ 
$g_\alpha = \lambda.id$ \ for each $\alpha$). 
This subgroup acts trivially on $W_0$.
A point $x\in W_0$ is called {\em $\chi$-semi-stable} if there exists an 
integer $n\geq 1$ and a polynomial $f\in\C\lbrack W_0\rbrack$ which is
$\chi^n$-invariant and such that $f(x)\not = 0$ ($f$ is called 
$\chi^n$-invariant if for every $w\in W_0$ and $g\in G_0$ we have
\ $f(gw) = \chi^n(g)f(w)$). The point $x$  is called {\em $\chi$-stable} if moreover

$\dim(G_0x) = \dim(G_0/\C^*)$ \ and if the action of $G_0$ on \
$\lbrace w\in W_0, f(w)\not = 0\rbrace$ \ is closed.

\bigskip

A. King proves in \cite{king} the following results :
\begin{enumerate}
\item [(1)] A point $x=(f_{\beta\alpha}) \in W_0$ is $\chi$-semi-stable
  (resp. $\chi$-stable) if and only if for each family 
  $(M'_\alpha),\ \alpha\in Q,$ of subspaces $M'_\alpha\subset M_\alpha$ 
  which is neither the trivial family (0) nor the given family $(M_\alpha)$
  and which satisfies
$$f_{\beta\alpha}(M'_\alpha\otimes V_{\alpha\beta})\subset
M'_\beta$$
for each $(\alpha,\beta)\in\Gamma$, we have

$$\sum_{\alpha\in Q}e_\alpha \dim(M'_\alpha) \ \leq \ 0 \ \ \rm{( resp. }
\ < 0 \ \rm{)}.$$

\item [(2)] Let $W_0^{ss}$ (resp. $W_0^s$) be the open subset of $W_0$
consisting of semi-stable (resp stable) points. Then there exists a good
quotient  
$$\pi : W^{ss}_0\lra M$$
by $G_0/\C^*$ which is a projective variety.

\item [(3)] The restriction of this quotient 
$$W^s_0\lra M^s = \pi(W^s_0)$$
is a geometric quotient and $M^s$ is smooth.
\end{enumerate}
\end{sub}\sepsub

\begin{sub}\label{polariz}{\sc\small Polarizations}\rm

The (semi-)stable points of $W_0$ remain the same if we replace 
$(e_\alpha)$ by $(ce_\alpha)$, $c$
being a positive integer. So the notion of (semi-)stability is fully
described by the reduced parameters $(\frac{e_\alpha}{t})$,
where 
$$t \ = \ \sigg_{\alpha\in Q, e_\alpha>0}e_\alpha m_\alpha.$$
So we can define the {\em polarization} of the action of $G_0$ on $W_0$ 
by any sequence 
$(c_\alpha)_{\alpha\in Q}$ of nonzero rational numbers such that
$$\sigg_{\alpha\in Q}c_\alpha m_\alpha = 0 \ , \ 
\sigg_{\alpha\in Q, c_\alpha>0}c_\alpha m_\alpha = 1 .$$
By multiplying this sequence by the smallest common denominator of the
$c_\alpha$ we obtain a sequence $(e_\alpha)$ of integers and the
corresponding character of $G_0$. Therefore the loci of stable and 
semi--stable points of $W_0$ with respect to $G_0$ and a polarization
$\Lambda_{0}=(c_{\alpha})$ are well defined and denoted by
\[
W_0^s(G_0,\Lambda_{0})\quad \text{and}\quad W_0^{ss}(G_0,\Lambda_{0}).
\]
\end{sub}
\sepsub

\begin{sub}\label{CondStab}{\sc\small Conditions imposed by the non-emptiness 
of the quotient}\rm

If $W_0^s$ is not empty, the $e_\alpha$ must satisfy some conditions. We will
derive this only in the three situations which occur in this paper. 
Polarizations satisfying these necessary conditions will be called 
{\em proper}. The first is that of the action of $G_{red}$ 
in \ref{subgrp} and the second is that of ${\mathbf G}$ and ${\mathbf W}$ in 
section 5, and the third is the case in between occurring in \ref{equal1}.

\sepsubsub

\begin{subsub}\label{first}First case\rm

Let $r,s$ be positive integers. We take
$$Q=\lbrace\alpha_1\ldots,\alpha_r,\beta_1,\ldots,\beta_s\rbrace,\qquad
\Gamma = \lbrace\alpha_1\ldots,\alpha_r\rbrace\times
\lbrace\beta_1,\ldots,\beta_s\rbrace.$$
This is the case of morphisms of type $(r,s)$. For $1\leq i\leq r$,
let \ $M'_{\alpha_i}=M_{\alpha_i}$ \ 
if \ $e_{\alpha_i}>0$, and $\lbrace
0\rbrace$ otherwise, and for $1\leq l\leq s$, let
$M'_{\beta_l}=M_{\beta_l}$. Then if one $e_{\alpha_i}$ is not positive,
we have
$$\sum_{\alpha\in Q}e_\alpha \dim(M'_\alpha) \ \geq \ 0 $$
and $(M'_\alpha) \not = (M_\alpha)$, 
so in this case no point of $W_0$ is stable.
So we obtain , if $W_0^s$ is non-empty, the conditions
$$e_{\alpha_i} > 0 \ , \ \text{for any\ } i,\quad \text{and}\quad
e_{\beta_l} < 0 \ , \ \text{for any\ } l.$$

A proper polarization is in this case a sequence $(\lambda_1,\ldots,\lambda_r,
-\mu_1,\ldots,-\mu_s)$ of rational numbers such that the $\lambda_i$ and 
the $\mu_l$ are positive and satisfy
$$\sigg_{1\leq i\leq r}\lambda_im_{\alpha_i} = 
\sigg_{1\leq l\leq s}\mu_lm_{\beta_l} = 1.$$
\end{subsub}

\sepsubsub

\begin{subsub}\label{second}Second case\rm

This case appears when we use a bigger reductive group to define the
quotient (this is the case of $\bf W$ later on).
Let $r,s$ be positive integers. Here we take
$$Q=\lbrace\alpha_1\ldots,\alpha_r,\beta_1,\ldots,\beta_s\rbrace,\qquad
\Gamma = \lbrace(\alpha_i,\alpha_{i-1}), 2\leq i\leq r,
(\alpha_1,\beta_s), (\beta_l,\beta_{l-1}), 2\leq l\leq s\rbrace.$$
Then the necessary conditions for $W_0^s$ to be non-empty are:
$$
\sum_{i\leq j\leq r}e_{\alpha_j}m_{\alpha_j} > 0\quad 
\text{for any\ } i,\ \ \text{ and\ }\sum_{1\leq l\leq m}e_{\beta_l}m_{\beta_l}
 < 0\quad \text{for any\ } m.
$$

To derive the first set of conditions we consider for any $i$ the family
$(M'_\gamma)$ for which $M'_{\alpha_j}=0$ if $i\leq j\leq r$ and
$M'_\gamma=M_\gamma$ for all other $\gamma\in Q$. Then $f_{\alpha \beta}
(M'_\alpha\otimes V_{\alpha\beta})\subset M'_\beta$ for any $f\in W_0$ and any
$(\alpha,\beta)\in \Gamma$. If $f$ is stable we obtain
\[
-\sum\limits_{i\leq j\leq r} e_{\alpha_j} m_{\alpha_j} =\sum\limits_{\gamma\in
  Q} e_\gamma \dim(M'_\gamma)<0
\]
Moreover, if the family $(M'_\gamma)$ is defined by $M'_{\alpha_j}=0$ for
$1\leq j\leq r, M'_{\beta_l}=0$ if $m<l\leq s$ and $M'_\gamma=M_\gamma$ else,
we obtain directly
\[
\sum\limits_{1\leq l\leq m} e_{\beta_l} m_{\beta_l} =\sum\limits_{\gamma\in Q}
e_\gamma\dim(M'_\gamma) < 0.
\]

A proper polarization in this case is then a sequence
$(\rho_1,\ldots,\rho_r,-\sigma_1,\ldots,-\sigma_s)$ of rational numbers 
satisfying
$$\sum_{1\leq i\leq r}\rho_im_{\alpha_i} =
\sum_{1\leq l\leq s}\sigma_lm_{\beta_l}=1.$$
and
$$\sum_{i\leq j\leq r}\rho_jm_{\alpha_j} > 0\quad \text{for any\ }i\quad
\text{and}\quad 
\sum_{1\leq l\leq m}\sigma_lm_{\beta_l} > 0\quad \text{for any\ }m.$$
We could also drop the normalization condition.
\end{subsub}\sepsubsub

\begin{subsub}\label{third}Third case\rm

This case is a combination of the first and second case. It appears in the
proof of the equivalence of semi--stability in \ref{conv2}. Here $Q$ is the 
same as in the previous cases and 
\[ \Gamma =\{(\alpha_i, \alpha_{i-1})\ 
,\ 2\leq i\leq r\ ,\ (\alpha_1, \beta_l)\ ,\ 1\leq l\leq s\}.  
\]
Now the necessary conditions for $W_0^s$ to be non--empty are: 
\[ 
\sum\limits_{i\leq j\leq r} e_{\alpha_j} m_{\alpha_j}>0\quad 
\text{for any}\ i,\quad \text{ and }\quad e_{\beta_l}<0\quad \text{ for any 
} l.  
\] 
The first condition follows as in the second case when we 
consider the family $(M'_\gamma)$ with $M'_{\alpha_j}=0$ for $i\leq 
j\leq r$ and $M'_\gamma=M_\gamma$ for all other $\gamma\in Q$.  The 
second condition follows when all $M'_\gamma$ are zero except 
$M'_{\beta_l}=M_{\beta_l}$ for one $l$.  Again a proper polarization 
in this case is a sequence $(\rho_1, \ldots, \rho_r, -\mu_1, \ldots, 
-\mu_l)$ with 
\[ \sum\limits_{1\leq i\leq r} \rho_i 
m_{\alpha_i}=\sum\limits_{1\leq l\leq s} \mu_l m_{\beta_l}=1 
\] 
and 
\[ 
\sum\limits_{i\leq j\leq r} \rho_j m_{\alpha_j}>0\quad\text{ for any } i 
\quad\text{ and } \mu_l>0\ \text{ for any } l.  
\]
 
\end{subsub}
\end{sub}
\sepsub

\begin{sub}\label{ProjAct}{\sc\small The action of $G_0$ on $\P(W_0)$}\rm

We suppose that we are in one of the first two preceding cases and that
there exist stable points in $W_0$. Let $P$ be a nonzero
homogeneous polynomial, $\chi^n$-invariant for some positive integer
$n$. The $\chi^n$--invariance implies that $P$ has degree $n.t$
where in case 1 (action of $G_{red}$ on $W$)
\[
t = \sum_{1\leq i\leq r}e_{\alpha_i}m_{\alpha_i},
\]
and in case 2 (action of $\bf G$ on $\bf W$)
\[
t = \sum_{1\leq i\leq r}ie_{\alpha_i}m_{\alpha_i}-
\sum_{1\leq l\leq s}(s-l)e_{\beta_l}m_{\beta_l}.
\]

To see this let $\lambda\in\C^\ast $ and let $g$ be given by
$g_{\alpha_i}=\lambda^{-1} id$ and $g_{\beta_l}=id$   in the first case and
by $g_{\alpha_i}=\lambda^{-i}id$ and $g_{\beta_l}=\lambda^{l-s}id$ in the
second case. Then $g x=\lambda x$ and $\chi^n(g)=\lambda^{nt}$ in both cases,
such that $P(\lambda x)=\lambda^{nt} P(x)$.

Now we will see that there exists a $G_0$-line bundle $\mathcal L$ on
$\P(W_0)$ such that the set $W_0^{ss}$ of semi-stable points is
exactly the set of points over $\P(W_0)^{ss}(G_0,{\mathcal L})$, which is the set of semi-stable points in the sense of Geometric Invariant
Theory corresponding to
\[
{\mathcal L} = {\mathcal O}_{\P(W_0)}(t),
\]
cf. \cite{mumf}, \cite{news}, \cite{p_v}.
Here the action of $G_0$ on $\kl$ is the natural action multiplied by
$\chi$. More precisely, the action of $G_0$ on $W_0$ induces an
action of this group on $S^tW_0$ and on $S^t W_0^*$ by:
\[
(g.F)(w) = F(g^{-1}w)
\]
for all $g\in G_0$, $w\in W_0$ and $F\in S^tW_0^*$, viewed as an
homogeneous polynomial of degree $t$ on $W_0$. The line bundle space 
$L$ of $\mathcal L$ is acted on by $G_0$ in the same way : if
$\xi\in L_{<w>}$ then \ 
$g.\xi\in L_{<gw>}$ \ is the form
on \ $<gw>^{*\otimes t} = L_{<gw>}$ \ given by \
$(g.\xi)(y) = \xi(g^{-1}y)$. We modify now the action of $G_0$ on
$L$ (resp. $S^tW_0^*$) by multiplying with $\chi(g)$ :
\[
g*\xi = \chi(g)g.\xi \ \ \ {\rm for \ } \xi\in L_{<w>},\ \text{ or }\ g*F = \chi(g)g.F \ \ \ {\rm for \ } F\in H^0(\P(W_0),{\mathcal L})
= S^tW_0^*.
\]

Now  $P\in H^0(\P(W_0),{\mathcal L}^{\otimes n})$ \ is an invariant section
if and only if $P$ is a homogeneous polynomial of degree $tn$ which
satisfies
$$P(gw) = \chi^n(g)P(w).$$
From the definition of semi-stable points in $W_0$ and $\P(W_0)$
with respect to the modified $G_0$-structure on 
${\mathcal L} = {\mathcal O}_{\P(W_0)}(t)$, we get immediately
\sepprop

\begin{subsub}\label{sslocus}{\bf Lemma:}
Assume that \ $W_0^s(G_0,\Lambda_{0})\not = \emptyset$ \ and let $t$ be defined 
as above in the two cases of $W_0$. Then the set $W_0^{ss}(G_0,\Lambda_{0})$ 
is the cone of the 
set $\P(W_0)^{ss}(G_0,{\mathcal O}_{\P(W_0)}(t))$ as defined in G.I.T.
\rm

\sepprop

There are two definitions of stable points in $\P(W_0)$, the classical one,
given in \cite{mumf}, \cite{news}, and a more recent one, given in \cite{p_v}.
If we take D.~Mumford's definition,
the cone of the set of stable points in $\P(W_0)$ does not coincide
with $W_0^s$  because every point of $\P(W_0)$ has a stabilizer of positive
dimension. In fact there is a subgroup of $G_0/\C^\ast$ of positive dimension
which acts trivially on $\P(W_0)$. In the first case for example such a group
is given by $g_{\alpha_i}=\lambda id$ and $g_{\beta_l}=\mu id$ with
$\lambda, \mu\in\C^\ast$. 
%Then the cone over $\P(W_0)^s
%(G_0/\C^\ast\times\C^\ast, \ko_{\P(W_0)}(t))$ coincides with $W_0^s$. 
If we want to keep the coincidence between the sets of stable points for one 
and the same group, we would have to consider the
action of a smaller reductive group in order to eliminate additional
stabilizers. We will do this in \ref{detgrp} only in the first case. If we
take the definition of V.L.~Popov and E.G.~Vinberg, then we obtain that the set 
$W_0^{s}(G_0,\Lambda_{0})$ is exactly the cone of the set 
$\P(W_0)^{s}(G_0,{\mathcal O}_{\P(W_0)}(t))$
\end{subsub}
\end{sub}

\sepsub
\begin{sub}\label{detgrp}{\sc\small The group $G'$}\rm

Let $G$ and $W$ be as in section 2 and let 
$\Lambda=(\lambda_1,\ldots,\lambda_r,-\mu_1,\ldots,-\mu_s)$ be a 
proper polarization as in \ref{first} for the action of $G_{red}$ on $W$. 
%It would suffice to consider the subgroup
%of $G_{red}$ consisting of elements $((g_i), (h_l))$ satisfying
%\[\prod_{1\leq i\leq r}det(g_i) = \prod_{1\leq l\leq s}det(h_l) = 1.\]
It is then convenient to use the subgroup $G'_{red}$ of
$G_{red}$ consisting of elements $((g_i),(h_l))$ satisfying
\[
\prod_{1\leq i\leq r}det(g_i)^{a_{i1}} = 
\prod_{1\leq l\leq s}det(h_l)^{b_{sl}} = 1, \text{ where } 
a_{ji}=\dim(A_{ji})\text{ and } b_{ml}=\dim(B_{ml}).
\]

We consider the action of $G'_{red}$ on $\mathcal L$ induced by the 
modified $\chi$-action of $G_{red}$.  Now the set $W^s(G_{red},\Lambda)$ of 
$\chi$-stable points of $W$ is exactly the cone over the locus 
$\P(W)^s(G'_{red},\kl)$ of stable points of $\P(W)$ in the sense of 
Geometric Invariant Theory.
\end{sub}

\sepsec
%\newpage
\section{Semi--stability in the non--reductive case}

Let $G$ and $W$ be as in section 2. A character $\chi$ on $G_{red}$ as in
King's setup can be extended to a character of $G$. Also the 
modified action of $G_{red}$ on $\mathcal L$ can be extended to an action of $G$.
 Let $G'$ be the subgroup of $G$ defined by the same equations as for $G'_{red}$. 
It contains $H$ and  $G'_{red}$, and we have $G'/H\simeq G'_{red}$.

In the case of the action of $G_{red}$ on $W$  a proper polarization 
is given by a sequence 
$\lambda_1,\ldots,\lambda_r,\mu_1,\ldots,\mu_s$ of positive rational 
numbers such that 
\[ 
\sigg_{1\leq i\leq r}\lambda_im_i = \sigg_{1\leq l\leq s}\mu_ln_l = 1.  
\] 
More precisely, the polarization is exactly 
the sequence $(\lambda_1,\ldots,\lambda_r,-\mu_1,\ldots,-\mu_s)$.  The 
parameter $\lambda_i$ (resp.  $\mu_l$) will be called the weight of 
the vector space $M_i$ (resp.  $N_l$).  We see that the dimension of 
the set of possible proper polarizations is $r+s-2$.  Let $t$ denote the 
smallest common denominator of the numbers $\lambda_i$ and $\mu_l$ and 
$\chi$ the character of $G_{red}$ defined by the sequence of integers 
$(-t\lambda_1,\ldots,-t\lambda_r,t\mu_1,\ldots,t\mu_s)$.  Let 
\[ 
{\mathcal L} = {\mathcal O}_{\P(W)}(t)\quad \text{with\ \ } 
t=\sigg_{1\leq i\leq r}m_it\lambda_i.  
\] 
As we have seen, if we consider the modified action of $G_{red}$ on 
$\mathcal L$, then the 
$\chi$-semi-stable points of $W$ are exactly those over the 
semi-stable points of $\P(W)$ in the sense of Geometric Invariant 
Theory with respect to the action of $G_{red}/\C^*$ on $\mathcal L$.  
The $\chi^{tn}$-invariant polynomials are the $G_{red}$-invariant 
sections of ${\mathcal L}^n$.

We are now going to define a notion of (semi--)stability for the points of $W$
with respect to the given action of the non--reductive group $G$. Let
$H\subset G$ be the above unipotent group, see also \ref{subgrp}.

\sepprop 

\begin{sub}\label{defstab}{\bf Definition:} A point $w\in W$ is called 
$G$--semi--stable (resp. $G$--stable) with respect to the (proper) polarization 
$\Lambda=(\lambda_1, \ldots, \lambda_r, -\mu_1, \ldots, -\mu_s)$ if every 
point of $Hw$ is $G_{red}$--semi--stable (resp.  $G_{red}$--stable) 
with respect to this polarization.

We denote these sets by $W^{ss}(G,\Lambda)$ resp. $W^s(G,\Lambda)$.

\rm For many of the quotient problems for the spaces of homomorphisms between
$\oplus m_i\ke_i$ and $\oplus n_j\kf_j$ and their cokernel sheaves this is a
fruitful notion. In \ref{prem_ex} we investigate an example with an explicit
description of the open sets $\ws\subset W^s(G_{red}, \Lambda)$. This example
also shows that the existence of a good quotient depends on the choice of the 
polarization.
\end{sub}

\sepprop

\begin{subsub}\label{defexam}{\it Situation for type (2,1):}\rm\\ 
In the case of morphisms of 
type $(2,1)$ we have $\mu_1=1/n_1$ and the polarization is completely described
by the single parameter \ $t=m_2\lambda_2$. We must have \ $0<t<1$.

A polarization such that there exists integers $m'_1$, $m'_2$, $n'_1$, 
with $0<n'_1<n_1$, $0\leq m'_i\leq m_i$, such that \ $m'_1n_1-m_1n'_1$,
$m'_2n_1-m_2n'_1$ \ are not both 0, and that
$$\lambda_1m'_1+\lambda_2m'_2 \ = \ \frac{n'_1}{n_1}$$
is called {\em singular}. There are only finitely many singular polarizations,
corresponding to the values \ $0<t_1<t_2<\cdots<t_p<1$ \ of $t$. Let $t_0=0$,
$t_{p+1}=1$. If $\Lambda$, $\Lambda'$ are polarizations corresponding to 
parameters $t$, $t'$ such that for some $i\in\lbrace 0,\cdots,p\rbrace$ we have 
$t_i<t,t'<t_{i+1}$, then 
$$W^{ss}(G,\Lambda)=W^{ss}(G,\Lambda') \ \ \ {\rm and} \ \ \ 
W^{s}(G,\Lambda)=W^{s}(G,\Lambda') .$$
Hence there are exactly $2p+1$ notions of $G$-(semi-)stability in this case.
Moreover, if $m_1$, $m_2$ and $n_1$ are relatively prime, and $\Lambda$ is a non
singular polarization, we have \ $W^{ss}(G,\Lambda)=W^{s}(G,\Lambda)$.

In the general case of morphisms of type $(r,s)$, it is not difficult to see
that there are only finitely many notions of $G$-(semi-)stability.
\end{subsub}

\sepsubsub

\begin{subsub}\label{Fauntleroy}{\bf Remark:\ }\rm
In \cite{faunt} semi-stability is defined as follows : A point $w\in W$ is
semi-stable if there exists a positive integer $k$ and a $G'$-invariant
section $s$ of ${\mathcal L}^k$ such that \ $s(w)\not = 0$ (there is
also a condition on the action of $H$). It is
clear that a semi-stable point in the sense of Fauntleroy is also
$G$-semi-stable with respect to 
$(\lambda_1,\ldots,\lambda_r,-\mu_1,\ldots,-\mu_s)$.
It is proved in \cite{faunt} that there exists a {\em categorical
quotient} of the open subset of semi-stable points in the sense of
\cite{faunt}, but it is not clear that all $G$-semi-stable points
are semi-stable. Moreover, in the general situation of \cite{faunt} 
there is no way to impose conditions which would imply that the categorical 
quotient is a good quotient or even projective. Using definition \ref{defstab}
we are able to derive a criterion for the existence of a good and 
projective quotient of $W$ under the action of $G$.
\end{subsub}

\sepsub

\begin{sub}\label{prem_ex}{\sc\small Existence and non--existence of good
    quotients, an example}\rm

We show here that we cannot expect that a good quotient $W^{ss}(\Lambda,G)//G$
will exist for any polarization $\Lambda$. 

We consider morphisms \ 
$2\ko(-2)\to\ko(-1)\oplus\ko$ \ on $\P_2$. There are 3 notions of
\hbox{$G$-(semi-)}stability in this case, two corresponding to non singular
polarizations. For one of the non singular polarizations the quotient 
$W^{s}(\Lambda,G)/G$ exists and for the other we prove the inexistence of a 
good quotient $W^{s}(\Lambda,G)//G$.

Let $V$ be a complex vector space of dimension 3, and $\P_2=\P V$.
Let \hfil\break \m{W=\Hom(2\ko(-2),\ko(-1)\oplus\ko)} \ on $\P_2$. 
A polarization for the action of $G$
on $W$ is a triple $(1/2,-\mu_1,-\mu_2)$ of positive rational numbers such that \
$\mu_1+\mu_2=1$. As in \ref{defexam} a polarization depends only on $\mu_1$.
There is only one singular polarization, corresponding to $\mu_1=1/2$. Hence
if we consider only non singular polarizations
there are only two notions of $G$-(semi-)stability, the first one corresponding
to polarizations such that $\mu_1>1/2$ and the second to polarizations such that
$\mu_1<1/2$. In both cases semi--stable points are already stable.
We are going to show that in the first case $\ws$ has a geometric quotient
which is projective and smooth and that in the second case $\ws$ doesn't even
admit a good quotient. 

The elements $x\in W$ and $g\in G$ are written as matrices
\[
x=\left(
  \begin{array}{ll}
z_1 & z_2\\
q_1 & q_2
  \end{array}\right)\qquad \text{ and }\qquad g= \left(\sigma ,\left(
  \begin{array}{ll}
\alpha & 0\\
z & \beta
  \end{array}\right)\right)
\]
where $z_1, z_2\in V^\ast, q_1, q_2 \in S^2 V^\ast, \sigma\in \GL(2)$, 
$\alpha,\beta\in\C^*$ and $z\in V^*$.

\sepsubsub

\begin{subsub}{The case \ $\mu_1> 1/2$}\label{case1}\rm\\
In this case $\ws$ has a geometric quotient which is the universal cubic
$Z\subset \P V\times \P S^3V^\ast$ of the Hilbert scheme of plane cubic curves
in $\P_2=\P V$. The quotient map is given by $x\mapsto (\langle z_1\wedge
z_2\rangle, \langle z_1q_2-z_2q_1\rangle)$

{\em Remark}: If $\mu_1>3/4$, then $\mu_1>3\mu_2$ and the conditions of
\ref{theo_main} (in the dual case $(1,2)$)for a good and projective quotient
to exist in this case are satisfied. 

The proof is done in several steps.
\vskip3mm

(1) {\bf claim 1}: Let $x\in W$ be as above. Then 
\begin{enumerate}
\item [(i)] $x\in W^s(G_{red}, \Lambda)$ if and only if $z_1\wedge z_2\neq 0$
  in $\Lambda^2 V^\ast$ and $q_1, q_2$ are not both zero.
\item [(ii)] $x\in \ws$ if and only if $z_1\wedge z_2\neq 0$ and
  $\det(x)=z_1q_2-z_2q_1\neq 0$ in $S^3 V^\ast$.
\end{enumerate}

{\em Proof of claim 1}. (i) follows easily from the criterion (1) in \ref{ResKing}. As
  for (ii) let $x\in W^s(G_{red}, \Lambda)$ with $\det(x)\neq 0$. Then
$\det(h.x)=\det(x)\neq 0$ 
for any $h=\left(\begin{smallmatrix}1 & 0\\z & 1\end{smallmatrix}\right)$
which implies that also $h.x\in W^s(G_{red}, \Lambda)$. Let conversely $x\in\ws$. 
Then $\det(x)\neq 0$ because otherwise there is a linear form $z\in
  V^\ast$ with $q_1=zz_1$ and $q_2=zz_2$ and with
  $h=\left(\begin{smallmatrix}1 & 0\\-z & 1\end{smallmatrix}\right)$  the
  element $h.x$ is the matrix $\left(\begin{smallmatrix}z_1 & z_2\\0 &
  0\end{smallmatrix}\right)$ which is not in $W^s(G_{red}, \Lambda)$.
\vskip3mm

(2) By the result of A. King in \ref{ResKing}, (3), there is a geometric
quotient $W^s(G_{red}, \Lambda)/G_{red}$ which is smooth and projective.

{\bf claim 2}: $W^s(G_{red}, \Lambda)/G_{red}\cong\P(Q^\ast\otimes
S^2V^\ast)$. Here $Q^\ast=\Omega^1$(1) is the dual of the tautological
quotient bundle over $\P V$. (The dimension of this quotient variety is 13
while $\dim W=18$ and $\dim G_{red}/\C^\ast=5)$.

To verify claim 2 we consider the map
\[
x=\left(
  \begin{array}{ll}
z_1 & z_2\\
q_1 & q_2
  \end{array}\right)\overset{\alpha}{\mapsto} (\langle z_1\wedge z_2\rangle,
\langle z_1\otimes q_2-z_2\otimes q_1\rangle)
\]
from $W^s(G_{red}, \Lambda)$ to $\P V\times \P(V^\ast\otimes S^2
V^\ast)\subset\P(V^\ast\otimes S^2 V^\ast\otimes \ko_{\P V})$ where we identify
$\P\Lambda^2 V^\ast$ with $\P V$ via $\langle z_1\wedge
z_2\rangle\leftrightarrow\langle a\rangle,\ z_1(a)=z_2(a)=0$. Then each
$\alpha(x)\in\P(Q^\ast_{\langle a\rangle}\otimes S^2V^\ast)$ 
because $Q^\ast_{\langle a\rangle}\subset V^\ast$ is the subspace of
forms vanishing in $\langle a\rangle$. It follows immediately that $\alpha$ is
a morphism
\[
W^s(G_{red}, \Lambda)\to \P(Q^\ast\otimes S^2 V^\ast)
\]
which is surjective and $G_{red}$--equivariant. It induces a morphism of the
geometric quotient to $\P(Q^\ast\otimes S^2 V^\ast)$ which is even
bijective. Since both, the quotient and the target are smooth, this is an
isomorphism.
\vskip3mm

(3) Since $Q^\ast\subset V^\ast\otimes \ko_{\P V}$ we have an induced
homomorphism $Q^\ast\otimes S^2 V^\ast\to S^3V^\ast\otimes \ko_{\P V}$. It is
the middle part of the canonical exact sequence
\[
0\to \Lambda^2 Q^\ast\otimes V^\ast\to Q^\ast\otimes S^2V^\ast\to
S^3V^\ast\otimes \ko_{\P V}\xrightarrow{ev}\ko_{\P V}(3)\to 0
\]
of vector bundles on $\P V$. Let $\kz$ be the kernel of $ev$. From the left
part of the sequence we obtain the affine bundle
\[
\P(Q^\ast\otimes S^2V^\ast)\smallsetminus \P(\Lambda^2Q^\ast\otimes
V^\ast)\xrightarrow{\beta}\P(\kz)\subset\P V\times \P S^3V^\ast.
\]
Here $\P(\kz)=Z$ is nothing but the universal cubic and the fibres of $\beta$
are isomorphic to $V^\ast$.

{\bf claim 3}: $\ws\subset W^s(G_{red},\Lambda)$ is the inverse image of
$\P(Q^\ast\otimes S^2V^\ast)\smallsetminus \P(\Lambda^2Q^\ast\otimes V^\ast)$
under $\alpha$ and $\alpha_{|\ws}$ is a geometric quotient with respect to
$G_{red}$.

{\em Proof of claim 3}. $z_1\otimes q_2-z_2\otimes q_1$ belongs to
  $\Lambda^2Q^\ast_{\langle a\rangle} \otimes V^\ast$ if and only if $
z_1q_2-z_2q_1 =0$, see (ii) of claim 1.

\vskip3mm

(4) Let now $\pi=\beta\circ\alpha$ be the morphism $\ws\to Z$, given by
$x\mapsto (\langle a\rangle, \langle z_1q_2-z_2q_1\rangle)$, where
$z_1(a)=z_2(a)=0$. It is obviously $G$--equivariant and its fibres coincide
with the $G$--orbits. Since $\alpha$ is a geometric quotient and $\beta$ is an
affine bundle, then $\pi$ is also a geometric quotient.
\vskip3mm

{\bf Remark}: The variety $Z$ is isomorphic to the moduli space
$M=M_{\P_2}(3m+1)$ of stable coherent sheaves on $\P_2$ with Hilbert
polynomial $\chi\kf(m)=3m+1$. This had been verified by J. Le Potier in
\cite{potier2}. The space $\ws$ is a natural parametrization of 
$M$ because any $\kf\in M$ can be presented in an extension sequence
$0\to\ko_C\to\kf\to\C_p\to 0$ where $C$ is the cubic curve supporting $\kf$
and $p\in C$, and then $\kf$ has a resolution
\[
0\to 2\ko(-2)\xrightarrow{x} \ko(-1)\oplus\ko\to \kf\to 0.
\]
This resolution is the Beilinson resolution as can easily be
verified. Moreover, $x$ is $(G, \Lambda)$--stable if and only if $\kf$ is
stable. (If $p$ is a smooth point of $C$, then $\kf$ is the line bundle 
$\ko_C(p)$ and if $p$ is a
singular point of $C$, then $\kf$ is the unique Cohen--Macaulay module on $C$
with the given polynomial). There is an obvious universal family 
${\mathbb \kf}$ on $\ws\times_H \P V$ which defines a $G$--equivariant 
morphism $\ws\to
  M$ and then a bijective morphism $Z\to M$, which by smoothness, is an
  isomorphism. One knows that $M$ carries a universal family $\ke$. This family
  can be obtained as the non--trivial extension
\[
0\to \ko_{Z\times_H Z}\to \ke\to \ko_\Delta\to 0,
\]
where $H=\P S^3V^\ast$ and $Z\times_H Z\subset Z\times \P V,$ or can be obtained
as the descent
of the family ${\mathbb \kf}$. More details can be found in \cite{frei}.
\end{subsub}

\sepsubsub

\begin{subsub}{The case \ $\mu_1< 1/2$}\label{case2}\rm\\
We suppose now that the polarization $\Lambda$ is such that \ $\mu_1<1/2$. In 
this case an element $x$ of $W$ is $G$-stable if and only if $z_1$, $z_2$ are not 
both zero, and if for every $z\in
V^*$, \m{q_1-zz_1} and \m{q_2-zz_2} are linearly independent. 
\end{subsub}
%\sepprop
\begin{subsub}{\bf Proposition:}\label{nonex} For this polarization there does 
not exist a good quotient \hfil\break $\ws//G.$
\end{subsub}

\begin{proof}
Let $z_1$ be a non-zero element of $V^*$, let \ 
\m{q\in S^2V^*\backslash z_1V^*}, and let $x\in W$ be the matrix
$$\left(\begin{array}{cc}z_1 & 0\\ q & z_1^2\end{array}\right).$$
Then $x$ is stable. 

\medskip

{\bf Claim : }{\em The orbit $Gx$ is closed and if $y\in W^s(G,\Lambda)$ is such
that $\ov{Gy}$ meets $Gx$, then $y\in Gx$.}

\medskip

Before proving the claim, we will show that it implies proposition \ref{nonex}.
The stabilizer of a generic point in $W^s(G,\Lambda)$ is isomorphic to $\C^*$ :
it consists of pairs of homotheties $(\lambda,\lambda)$. It follows that if \ 
\m{M=W^s(G,\Lambda)//G} exists, then all the fibers of the quotient morphism \ 
\m{\pi:W^s(G,\Lambda)\to M} \ are of dimension at least \m{\dim(G)-1}. Now
suppose that the claim is true. Then this implies that \
\m{\pi^{-1}(\pi(x))=Gx}. But the stabilizer $G_x$ of $x$ has dimension 2 : it
consists of pairs
$$(\left(\begin{array}{cc}\alpha & 0\\ \beta & \alpha\end{array}\right),
\left(\begin{array}{cc}\alpha & 0\\ \beta z_1 & \alpha\end{array}\right))$$
with $\alpha\in\C^\ast$, $\beta\in\C$, and hence has dimension 2. It follows
that \hfil\break \m{\dim(\pi^{-1}(\pi(x)))<\dim(G)-1}, a contradiction.

\medskip

{\em Proof of the claim.} Let $y\in W^s(G,\Lambda)$ such that $x\in \ov{Gy}$.
Let
$$y \ = \ \left(\begin{array}{cc}z & z_2\\ q_1 & q_2\end{array}\right) .$$
Then $z_1$ is contained in the vector space spanned by $z$ and $z_2$. Hence by
replacing $y$ with an element of $Gy$ we can assume that $z=z_1$ and that
$z_2=0$ if $z_2$ is a multiple of $z_1$.

According to lemma \ref{util_lemm} there exists a smooth irreducible curve $C$,
$x_0\in C$, and a morphism
$$\theta : C\backslash\lbrace x_0\rbrace\lra G $$
such that
$$\xymatrix@R=4pt{
\ov{\theta} : C\backslash\lbrace x_0\rbrace \ar[r] & W \\
{\ \ \ \ \ \ \ \ \ }t \fmaps[r] & \theta(t)y \\
}$$
can be extended to \ \m{\ov{\theta}:C\to W}, with \ \m{\ov{\theta}(x_0)=x}.
We can write, for \ \m{t\in C\backslash\lbrace x_0\rbrace},
$$\ov{\theta}(t) \ = \ 
\left(\begin{array}{cc}a(t)z_1+b(t)z_2 & c(t)z_1+d(t)z_2\\ q_1(t) & q_2(t)
\end{array}\right) $$
with
$$(1) \ \ \ \ \ \ \ \ \ \ 
q_1(t) \ = \ \lambda(t)\big(a(t)q_1+b(t)q_2+u(t)z_1\big) ,$$
$$(2) \ \ \ \ \ \ \ \ \ \ 
q_2(t) \ = \ \lambda(t)\big(c(t)q_1+d(t)q_2+u(t)z_2\big) ,$$
where $\lambda$, $a$, $b$, $c$, $d$ are morphisms \ 
\m{C\backslash\lbrace x_0\rbrace\to\C} \ and \ 
\m{u :C\backslash\lbrace x_0\rbrace\to V^*}. The morphisms $\lambda$, $a$, $b$, 
$c$, $d$ can be extended to morphisms \ \m{C\to\P_1=\C\cup\lbrace
0,\infty\rbrace}, denoted by $\ov{\lambda}$, $\ov{a}$, $\ov{b}$, $\ov{c}$, 
$\ov{d}$ respectively, and $u$ extends to \ \m{\ov{u}:C\to\P(V^*\oplus\C)}. 
Now we use the fact that $\ov{\theta}$ is defined at $x_0$. The first
consequence is that \ \m{\ov{a}(x_0)=1}, \m{\ov{c}(x_0)=0}, and if $z_2\not=0$
then $\ov{b}$  and $\ov{d}$ also vanish at $x_0$. 
The second is that the morphisms 
\m{q_1,q_2:C\backslash\lbrace x_0\rbrace\to S^2V^*} \ can be extended to \
\m{\ov{q_1},\ov{q_2}:C\to S^2V^*}, and we have \ \m{\ov{q_1}(x_0)=q},
\m{\ov{q_2}(x_0)=z_1^2}. 

We will now consider three cases : \m{\ov{\lambda}(x_0)=0},
\m{\ov{\lambda}(x_0)=\infty}, \m{\ov{\lambda}(x_0)\in\C^*}.

\medskip

Suppose that \m{\ov{\lambda}(x_0)=0}. If $z_2\not=0$, then $(1)$ implies that
$\ov{q_1}(x_0)=q$ is a multiple of $z_1$, but this is not true. If $z_2=0$ then
$(2)$ implies that $q_2$ is a multiple of $z_1^2$ and $(1)$ implies then that 
$q$ is also a multiple of $z_1$, which is not true. Hence we cannot have
\m{\ov{\lambda}(x_0)=0}.

\medskip

Suppose that \m{\ov{\lambda}(x_0)=\infty}. If $z_2\not=0$, then $(1)$ implies
that
$$\xymatrix@R=4pt{
\mu : C\backslash\lbrace x_0\rbrace \ar[r] & S^2V^* \\
{\ \ \ \ \ \ \ \ \ }t \fmaps[r] & a(t)q_1+b(t)q_2+u(t)z_1 \\
}$$
and
$$\xymatrix@R=4pt{
\eta : C\backslash\lbrace x_0\rbrace \ar[r] & S^2V^* \\
{\ \ \ \ \ \ \ \ \ }t \fmaps[r] & c(t)q_1+d(t)q_2+u(t)z_1 \\
}$$
extend to morphisms $C\to S^2V^*$ which vanish at $x_0$. It follows from the
fact that \m{\mu(x_0)=0} that \
\m{u=\ov{u}(x_0)\in V^*}, and that \m{q_1=-uz_1}. Since \m{q_1\not=0} (by
$G$-stability of $y$), we have \m{u\not=0}. But since \ 
\m{\ov{c}(x_0)=\ov{d}(x_0)=0}, this contradicts the fact that \m{\eta(x_0)=0}.

If $z_2=0$ then we deduce from the fact that \m{\mu(x_0)=0} that \
\m{q_1\in\ \pline{q_2,V^*z_1}}, which contradicts the $G$-stability of $y$.

\medskip

It follows that we have \ \m{\delta=\ov{\lambda}(x_0)\in\C^*}. If $z_2\not=0$,
using the fact that \ \m{\ov{a}(x_0)=1} \ and
\m{\ov{b}(x_0)=\ov{c}(x_0)=\ov{d}(x_0)=0} \ we see that
\ \m{u=\ov{u}(x_0)\in V^*} \ and that \ \m{z_1^2=\delta uz_2}, which contradicts
the fact that \ \m{z_1\wedge z_2\not=0}.

Hence we have $z_2=0$. It follows from $(2)$ that \m{\ov{d}(x_0)\in\C^*} \ and
that \ \m{z_1^2=\delta\ov{d}(x_0)q_2}. By $(1)$ we see that
$$\xymatrix@R=4pt{
\epsilon : C\backslash\lbrace x_0\rbrace \ar[r] & S^2V^* \\
{\ \ \ \ \ \ \ \ \ }t \fmaps[r] & b(t)q_2+u(t)z_1 \\
}$$
extends to $C$ and that \ 
$$\epsilon(x_0) \ = \ \frac{1}{\delta}q-q_1 .$$
We have, if $t\not=x_0$
$$\epsilon(t) \ = \ z_1(\frac{b(t)}{\delta\ov{d}(x_0)}z_1+u(t)) .$$
It follows that $\epsilon(x_0)$ is a multiple of $z_1$ : \m{\epsilon(x_0)=z_1v}.
We have then 
$$q_1 \ = \ \frac{1}{\delta}q-z_1v$$
and
$$y \ = \ \left(\begin{array}{cc}z_1 & 0\\ q_1 & q_2\end{array}\right) \ = \
\left(\begin{array}{cc}z_1 & 0\\ \frac{1}{\delta}q-z_1v & 
\frac{1}{\delta\ov{d}(x_0)}z_1^2\end{array}\right) \ \in \ Gx$$
as claimed.

It remains to show that $Gx$ is closed. This can be proved easily by computing
the stabilizers of all the points in $W^s(G,\Lambda)$. We see then that $G_x$
has the maximal possible dimension, hence $Gx$ is closed.
\end{proof}

\sepprop

We now give a proof of the lemma used in the preceeding proposition :

\sepprop

\begin{subsub}\label{util_lemm}{\bf Lemma : } Let $W$ be a finite dimensional
vector space, $G$ a linear algebraic group acting algebraically on $W$, $y\in
W$ and \ \m{x\in\ov{Gy}\backslash Gy}. Then there exists a smooth curve $C$,
\m{x_0\in C} and a morphism 
$$\theta : C\backslash\lbrace x_0\rbrace\lra G$$
such that the morphism
$$\xymatrix@R=4pt{
\ov{\theta} : C\backslash\lbrace x_0\rbrace \ar[r] & W \\
{\ \ \ \ \ \ \ \ \ }t \fmaps[r] & \theta(t)y \\
}$$
extends to \ \m{\ov{\theta}:C\to W} \ and that \ \m{\ov{\theta}(x_0)=x}.
\end{subsub}

\begin{proof}
Let \ $n=\dim(W)$, $d=\dim(Gy)$. The generic $(n-d+1)$-dimensional affine
subspace \m{F\subset W} through $x$ meets \m{\ov{Gy}} on a curve, and meets
\m{\ov{Gy}\backslash Gy} in a finite number of points. Hence we can find a
curve \m{X\subset\ov{Gy}} that meets \m{\ov{Gy}\backslash Gy} only at $x$.
Taking the normalization of $X$ and substracting a finite number of points or
unnecessary components if needed, we obtain a morphism \m{\alpha:Z\to\ov{Gy}} 
(where $Z$ is a smooth curve) and a point \m{z_0\in Z} such that
\m{\alpha(z_0)=x} and \ \m{\alpha(Z\backslash\lbrace z_0\rbrace)\subset Gy}.
Consider now the restriction of $\alpha$
$$Z\backslash\lbrace z_0\rbrace\lra Gy\simeq G/G_y .$$
There exists a smooth curve $Z'$ and an etale surjective morphism \
\m{\phi:Z'\to Z\backslash\lbrace z_0\rbrace} \ such that the principal 
$G_y$-bundle \m{\phi^*\alpha^*G} on $Z'$ is locally trivial. By considering
completions $\ov{Z'}$, $\ov{Z}$ of $Z'$, $Z$ and an extension of $\phi$ to a
morphism \ \m{\ov{Z'}\to\ov{Z}} we obtain a smooth curve $Y$, $y_0\in Y$ and
a morphism \ \m{\beta: Y\to Z} \ such that \m{\beta(y_0)=z_0} and that the
principal $G_y$-bundle \ \m{\Gamma=\beta^*\alpha^*G} \ is defined on 
\m{Y\backslash\lbrace y_0\rbrace} and locally trivial. Let 
$U\subset Y$ be a nonempty open subset such that we have a $G_y$-isomorphism 
$$\gamma:\Gamma_{\mid U}\simeq U\times G_y$$
Then we can take \ \m{C=U\cup\lbrace y_0\rbrace}, \m{x_0=y_0}, and for \
\m{t\in C\backslash\lbrace x_0\rbrace=U}, we have
$$\theta(t) \ = \ \psi(\gamma^{-1}(t,e)),$$
where $\psi$ is the canonical morphism $\Gamma\to G$.
\end{proof}

\end{sub}

\sepsub

\begin{sub}{\sc\small More general counterexamples of inexistence of geometric
quotients}\label{remgeoqu}\rm

Let $W$ be the space of homomorphisms
\[
\ko(-2)\oplus\ko(-1)\to {\C}^{2n}\otimes \ko(1)
\]
over $\P_n$ and let the homomorphism $\phi_0\in W$ be given by the matrix
\[
\left(
  \begin{array}{cc}
z_0^2z_1 & z_1^2\\
\vdots & \vdots\\
z_0^2z_n & z_n^2\\
z_0z_1^2 & 0\\
\vdots & \vdots\\
z_0 z_n^2 & 0
  \end{array}
\right)
\]
where the $z_\nu$ are homogeneous coordinates. The stabilizer of $\phi_0$
contains ${\C}^\ast$ and the pairs
\[
\left(
  \begin{array}{cc}
1 & 0\\
a z_0 & 1
  \end{array}
\right)\ ,\qquad
\left(\begin{array}{cc}
I_n & -aI_n\\
0 & I_n
\end{array}
\right)
\]
in $\Aut(\ko(-2) \oplus \ko(-1))\times \GL({\C}^{2n})$ and thus has
dimension $\geq 2$. If $\Lambda=(\lambda_1, \lambda_2, -\mu_1)$ is a
polarization with $0<\lambda_1,\ 0<\lambda_2<\frac{1}{2}$, then it is easy to 
see that
$\phi_0$ is $\Lambda$--stable in the sense of \ref{defstab}. For example
$(m'_1, m'_2, n')=(0,1,n)$ is the dimension vector of a $\phi_0$--invariant
choice of subspaces with
$\lambda_1m'_1+\lambda_2m'_2-\mu_1n'=\lambda_2-1/2<0$. There are however
stable homomorphisms $\phi \in W$ with stabilizer ${\C}^\ast$. Therefore 
$\ws /G$ can never admit the structure of a geometric
quotient.  
We will see in \ref{sone} that a sufficient condition
for that in the case of this $W$ is $\lambda_2>(n+1)\lambda_1$ or
$\lambda_2>\frac{n+1}{n+2}$ because $\lambda_1+\lambda_2=1$.
\end{sub}
\sepsec
%\newpage
\section{Embedding into a reductive group action}

We will construct an algebraic reductive group $\bf G$, a finite
dimensional vector space $\bf W$ on which $\bf G$ acts algebraically,
and an injective morphism
$$\zeta : W\lra \bw$$ compatible with a morphism of groups $$\theta : 
G\lra \bg.$$ The traces of $\bf G$-orbits on $\zeta(W)$ will be 
exactly the $G$-orbits.  The space $\bf W$ is of the same type as 
those studied in \ref{ResKing}.  We will associate naturally to any 
polarization of the action of $G$ on $W$ a character $\chi$ of ${\bf 
G}/\C^{*}$, i.e.  a polarization of the action of $\bf G$ on $\bf W$.  
We will prove that in certain cases a point $w$ of $W$ is 
$G-$(semi-)stable with respect to the given polarization if and only 
if $\zeta(w)$ is $\chi$-(semi-)stable with respect to the associated 
polarization.  The existence of a good and projective quotient of the 
open set of $G$-semi-stable points will follow from this.  
\sepsub

\begin{sub}\label{motiv}{\sc\small Motivation in terms of sheaves}\rm

The idea for the embedding of $W$ into a space $\mathbf W$ with a 
reductive group action is to replace the sheaves $\ke_i$ in 
$\ke=\oplus (\ke_i\otimes M_i)$ by 
$\ke_1\otimes \Hom(\ke_1,\ke_i)$ and dually the sheaves $\kf_l$ in
$\kf = \oplus (\kf_l\otimes N_l)$ by
$\kf_s\otimes \Hom(\kf_l,\kf_s)^\ast$
and then to consider the induced composed homomorphisms $\gamma(\Phi)$ 
for $\Phi\in \Hom(\ke,\kf)=W$

$$\ke_1\otimes \Hom(\ke_1,\ke)\to\ke\to\kf\to\kf_s\otimes 
\Hom(\kf,\kf_s)^\ast$$

in the bigger space $\tilde{W}$ of all homomorphisms between 
$\ke_1\otimes \Hom(\ke_1,\ke)$ and $\kf_s\otimes \Hom(\kf,\kf_s)^\ast$. 
This space is naturally acted on by the reductive group
\[\tilde{G}=\GL(\Hom(\ke_1,\ke))\times \GL(\Hom(\kf,\kf_s)^\ast).\]
However it is not suitable enough for our purpose by two reasons.  It 
does not allow enough polarizations as in section 3 for direct sums in 
order to have consistency of (semi-)stability and, secondly the group 
actions $G\times W\to W$ and $\tilde{G}\times\tilde{W}\to \tilde{W}$ 
don't have consistent orbits.  Both insufficiencies are however 
eliminated when we consider the following enlargement of $\tilde{W}$.  
We set
$$P_{i}=\Hom(\ke_i,\ke)\ \ \text{and}\ \ Q_{l}=\Hom(\kf,\kf_l)^\ast,$$

and introduce the auxiliary spaces
$$\bw_L = \underset{1<i\leq r}{\oplus}\Hom(P_i\ot \Hom({\mathcal E}_{i-1},
{\mathcal E}_i),P_{i-1}), \ \
\bw_R = \underset{1\leq l<s}{\oplus}\Hom(Q_{l+1}\ot \Hom({\mathcal F}_l, 
{\mathcal F}_{l+1}),Q_l),$$ and define $$\bw = \bw_L\oplus 
\Hom({\mathcal E}_1\ot P_1,{\mathcal F}_s\ot Q_s)\oplus\bw_R .$$

There are distinguished elements 
$$(\xi_{2},\cdots ,\xi_{r})\in \bw_L,\quad (\eta_{1},\cdots , 
\eta_{s-1})\in \bw_R$$ whose components are the natural composition 
maps.  The embedding of $W$ into $\bw$ will be defined as the affine 
map

$$W\overset{\zeta}{\lra} \bw,\quad\quad \Phi\mapsto ((\xi_{2},\cdots 
,\xi_{r}), \gamma(\Phi), (\eta_{1},\cdots ,\eta_{s-1})),$$

where  $\gamma(\Phi)$ is the above composition for a given $\Phi\in W.$ 
The components of 
$\bw_L$ and $\bw_R$ will guarantee a compatible action of a reductive 
group and at the same time the possibility of choosing enough 
polarizations for this action.
\sepprop

\begin{subsub}\label{Remmotiv}{\bf Remark}: \rm One might hope to be able to do
  induction on $r$ and/or $s$ by simply replacing $M_{r-1}\otimes
  \ke_{r-1}\oplus M_r\otimes \ke_r$ by $(M_{r-1} \oplus M_r\otimes
  \Hom(\ke_{r-1},\ke_r))\otimes \ke_{r-1}$ and keeping the other $\ke_i$ for
  $i<r-1$. But then we drop the information about the homomorphisms $\ke_i\to
  \ke_r$. Therefore we are lead to replace all $\ke_i, i\geq 2$, by $\ke_1$ at
  a time, i.e. by
\[
P_1\otimes \ke_1=(M_1\oplus M_2\otimes A_{21}\oplus\cdots\oplus M_r\otimes
A_{r1})\otimes \ke_1,
\]
where $A_{ji}=\Hom(\ke_i,\ke_j)$. Moreover, in order to keep the information
of the homomorphisms $\ke_i\to \ke_j$ for $2\leq i\leq j$ we consider also the
spaces
\[
P_i=M_i\oplus M_{i+1}\otimes A_{i+1, j}\oplus\ldots\oplus M_r\otimes A_{ri}
\]
together with the maps $P_i\otimes A_{i,i-1}\to P_{i-1}$ in the following. The
reader may convince himself that only because of this the actions of the
original group is compatible with the action of the bigger reductive group. It
is a beautiful outcome that then we are able to compare the semi--stability
with respect to related polarizations in section 7.

\end{subsub}
\end{sub}
\sepsub

\begin{sub}\label{defbigW}{\sc\small The abstract definition of $\bw$} \rm

The above motivating definition of the space $\bw$ can immediately be 
turned into the following final definition using the spaces $H_{li}, 
A_{ji}\ \text{and\ } B_{ml}$ and the pairings between them.  For any 
possible $i$ and $l$ we introduce the spaces $$P_i = \underset{i\leq j\leq 
r}{\oplus}M_j\ot A_{ji}\quad \text{and\ \ } Q_l = \underset{1\leq m\leq l}
{\oplus}N_m\ot B^\ast_{lm},$$

and we denote by $p_i$ and $q_l$ their dimensions.  For $1<i$ and $l<s$ we let 
$$ P_i\otimes 
A_{i,i-1}\overset{\xi_i}{\lra} P_{i-1}\quad \text{and\ \ } Q_{l+1}\ot 
B_{l+1,l}\overset{\eta_l}{\lra} Q_l$$

be the canonical morphisms, defined as follows. On the component
$M_j\ot A_{ji}$ of $P_{i}$, the map $\xi_i$ is the map
$$(M_j\ot A_{ji})\ot A_{i,i-1}\lra M_j\ot A_{j,i-1}$$

induced by the composition map of the spaces $A$. 
The map $\eta_l$ is defined in the same way. As in \ref{motiv} we set

$$\bw_L = \underset{1<i\leq r}{\oplus}\Hom(P_i\ot A_{i,i-1},P_{i-1}),\quad 
\quad 
\bw_R = \underset{1\leq l<s}{\oplus}\Hom(Q_{l+1}\ot B_{l+1,l}, Q_l),$$

and
$$\bw = \bw_L\oplus \Hom(P_1,Q_s\ot H_{s1})\oplus \bw_R.$$ 
 
In order to define the embedding $\zeta$ we define the operator 
$\gamma$ as follows.  Given $w=(\phi_{li})\in W\ \text{with\ \ } 
\phi_{li}\in \Hom(M_{i}, N_{l}\ot H_{li})$, we let $$\gamma(w)\in 
\Hom(P_1, Q_s\ot H_{s1})=\Hom(P_1\ot H_{s1}^*, Q_s)$$

be the linear map defined by the matrix 
$(\gamma_{li}(w))$, for which each $\gamma_{li}(w)$
is the composed linear map
$$M_i\ot A_{i1}\lra N_l\ot H_{li}\ot A_{i1}\lra
N_l\ot H_{l1}\lra N_{l}\ot B_{sl}^*\ot H_{s1},$$

where the first map is induced by $\phi_{li}$, the second by the 
composition $H_{li}\ot A_{i1}\to H_{l1}$ 
and the third by the dual composition $H_{l1}\to B_{sl}^*\ot H_{s1}.$

The map $\zeta$ can now be defined by

$$W\overset{\zeta}{\lra} \bw,\quad\quad w\mapsto ((\xi_{2},\cdots 
,\xi_{r}), \gamma(w), (\eta_{1},\cdots ,\eta_{s-1})).$$ 
\sepprop

\begin{subsub}\label{inject}{\bf Lemma:} The linear map $\gamma$ is 
injective and hence the morphism $\zeta$ is a closed embedding of affine 
schemes.
\end{subsub}

\begin{proof} From the surjectivity assumptions in \ref{abstr} we find
that dually the composition  
$$ H_{li}\lra H_{l1}\ot A_{i1}^*\lra B_{sl}^*\ot H_{s1}\ot A_{i1}^{*}$$
is injective. Now it follows from the definition of $\gamma_{li}(w)$ 
that $\phi_{li}$ can be recovered from $\gamma_{li}(w)$, by 
shifting $A_{i1}$ to its dual.
\end{proof}\sepsub
\end{sub}
\sepsub

\begin{sub}\label{bigG}{\sc\small The new group $\bg$}\rm

We consider now the natural action on $\bw$ as described in 
\ref{ResKing} in the general situation, where the group is $$\bg = 
\bg_L\times \bg_R,\quad \text{with\ \ } \bg_L = \prod_{1\leq i\leq 
r}\GL(P_i), \ \ \ \bg_R = \prod_{1\leq l\leq s}\GL(Q_l).$$

To be precise, this action is described in components by 
$$ g_{i-1}\circ x_{i-1,i}\circ (g_{i}\ot id)^{-1}, \quad h_{s}\circ 
\psi\circ (g_{1}\ot id)^{-1}\quad \text{and\ \ } 
h_{l}\circ y_{l,l+1}\circ (h_{l+1}\ot id)^{-1},$$

with 
$$x_{i-1,i}\in \Hom(P_{i}\ot A_{i,i-1}, P_{i-1}), \quad 
\psi\in \Hom(P_{1}\ot H_{s1}^{*}, Q_{s}),\quad  
y_{l,l+1}\in \Hom(Q_{l+1}\ot B_{l+1,l}, Q_{l})$$ 

and with 
$$g_{i}\in \GL(P_{i}),\quad h_{l}\in \GL(Q_{l}).$$  
The first and third expression describe the natural actions of 
$\bg_{L}$ on $\bw_{L}$ and of $\bg_{R}$ on $\bw_{R}.$
 
There are also natural embeddings of $G_{L},\ G_{R},\ G$ 
into $\bg_{L},\ \bg_{R},\ \bg$ respectively.  For that it is enough to 
describe the embedding of $G_{L}$ in $\bg_{L}$.  Given an element 
$g\in G_{L},$ \[ g=\left( \begin{array}{cccc} g_1 & 0 & \ldots & 0\\
u_{21} & g_2 & & \vdots\\
\vdots & \ddots & \ddots & 0\\
u_{r1} & \ldots & u_{r, r-1} & g_r
  \end{array}
\right)
\]

with $g_{i}\in \GL(M_{i})\ \text{and\ } u_{ji}\in \Hom(M_{i},M_{j}\ot 
A_{ji})$ we define $\theta_{L,i}(g)\in \GL(P_{i})$ as the matrix
\[
\theta_{L,i}(g)=\left(  
\begin{array}{cccc}
\tilde{g}_i & 0 & \ldots & 0\\
\tilde{u}_{i+1,i} & \tilde{g}_{i+1} & & \vdots\\
\vdots & \ddots & \ddots & 0\\
\tilde{u}_{r,i} & \ldots & \tilde{u}_{r,r-1} & \tilde{g}_r
  \end{array}
\right)
\]

with respect to the decomposition of $P_{i}$ with the following 
components: $\tilde{g}_{j}=g_{j}\ot id$ on $M_{j}\ot A_{ji}$ and 
for $i\leq j\leq k$ the map $\tilde{u}_{kj}$ is the composition
$$M_{j}\ot A_{ji}\lra M_{k}\ot A_{kj}\ot A_{ji}\lra M_{k}\ot A_{ki},$$
where the second arrow is induced by the given pairing. In case $j=i$ 
we have $\tilde{g}_{i}=g_{i}$ and $\tilde{u}_{ki}=u_{ki}$.
Now we define the map 
$$G_{L}\overset{\theta_{L}}{\to}{\bf G_{L}}\quad \text{by}\quad
 g\mapsto (\theta_{L,1}g, \cdots, \theta_{L,r}g).$$

It is then easy to verify that $\theta_{L}$ is an injective group 
homomorphism and defines a closed embedding of algebraic groups. With 
this embedding we consider $G_{L}$ as a closed subgroup of $\bg_{L}$. 
In the same way we obtain a closed embedding $\theta_{R}$ of 
$G_{R}\subset \bg_{R}$. Finally we obtain the closed embedding $\theta = 
(\theta_{L},\theta_{R})$ of $G\subset \bg$.  \end{sub}

\sepprop

\begin{subsub}\label{stabilizer}{\bf Lemma:}
The subgroup $G_L\subset\bg_L$ (respectively $G_R\subset\bg_R$) is the 
stabilizer of the distinguished element $(\xi_2,\ldots,\xi_r)\in {\bf 
W}_L$ (respectively $(\eta_{1},\ldots,\eta_{s-1})\in \bw_R$)
\end{subsub}

\begin{proof} It is enough to prove the statement only for $G_{L}$ 
because of duality. The fact that $G_L$ stabilizes $(\xi_2,\ldots,\xi_r)$  
is an easy consequence of the properties of the composition maps. 
The converse can be proved by induction on $r$. It is trivial for
$r=1$. Suppose that $r\geq 2$ and that the statement is true for $r-1$. Let
$(\gamma_1,\ldots,\gamma_r)$ be an element of the stabilizer of
$(\xi_2,\ldots,\xi_r)$. When we  replace the space $W$ by $W'$, 
corresponding to the spaces $M_2,\ldots,M_r$ and the same spaces $N_{l}$ and 
similarly $\bw_{L}$ by ${\bf W'}_{L}$, then 
$(\gamma_2,\ldots,\gamma_r)$ is an element of the stabilizer of 
$(\xi_3,\ldots, \xi_r)$, so by the induction hypothesis it belongs to 
$G'_L$ and there exists an element \[ g'=\left( \begin{array}{cccc} 
g_2 & 0 & \cdots & 0\\
u_{32} & g_3 & & \vdots\\
\vdots & \ddots & \ddots & 0\\
u_{r2} & \cdots & u_{r, r-1} & g_r
  \end{array}
\right)
\]
such that $(\gamma_2,\ldots,\gamma_r)=\theta'_L(g')$.
Let now $\gamma_1\in \GL(P_1)$ have the components
$$M_i\ot A_{i1}\overset{y_{ji}}{\lra}M_j\ot A_{j1}\quad \text{for all}\quad 
1\leq i,j\leq r.$$
The identity $\gamma_1\circ\xi_2=\xi_2\circ\gamma_2$
then shows that $y_{ji}=0$ for $j<i,\ y_{ii}=g_i$ for $2\leq i$ and 
$y_{ji}=u_{ji}$ for $2\leq j<i$. Now let $g_1=y_{11}$, $u_{j1}=y_{j1}$, 
for $2\leq j\leq r$, which are
linear mappings $M_1\lra M_j\ot A_{j1}$. Then
\[
g=\left(  \begin{array}{cccc}
g_1 & 0 & \cdots & 0\\
u_{21} & g_2 & & \vdots\\
\vdots & \ddots & \ddots & 0\\
u_{r1} & \cdots & u_{r, r-1} & g_r
  \end{array}
\right)
\]
is an element of $G_L$ and we have $(\gamma_1,\ldots,\gamma_r)=\theta_L(g).$
\end{proof}

\sepprop

{\bf Remark}: since the action of $\bg_L$ on $\bw_L$ is linear, it is 
clear that we have an isomorphism $$\bg_L/G_L\ \simeq \ {\bf 
G}_L(\xi_2,\ldots,\xi_r),\quad \text{and similarly}\quad {\bf 
G}_R/G_R\ \simeq \ \bg_R(\eta_1,\ldots,\eta_{s-1}).$$ We will use this 
fact in section 8.  

\sepprop

Using the associativity of the composition maps it is again easy to 
verify that the actions of $G$ on $W$ and $\bg$ on $\bw$ are compatible, 
i.e. that the diagram 
\[
\begin{CD}
G\times W @ >>> W\\
@V\theta\times\zeta VV @ VV\zeta V\\
\bg\times \bw @ >>> \bw
\end{CD}
\]
is commutative, in which the horizontal maps are the actions. In 
addition we have the 
\sepprop

\begin{subsub}\label{orbits}{\bf Corollary:}
Let $w,w'\in W$. Then $w$ and $w'$ are in the same $G$-orbit in $W$ if and
only if $\zeta(w)$  and $\zeta(w')$ are in the same $\bf G$-orbit in 
$\bf W$.
\end{subsub}

\begin{proof} It follows from the compatibility of the actions that if 
$g.w=w'$ in $W$ then also $\theta(g).\zeta(w)=\zeta(w')$ in $\bw$ by the last 
diagram. Conversely, if 
${\bf g}\in \bg$ and ${\bf g}.\zeta(w)=\zeta(w')$ then g stabilizes 
$(\xi_{2},\cdots,\xi_{r}, \eta_{1},\cdots ,\eta_{s-1})$ by the 
definition of $\zeta$ in \ref{defbigW}.  By Lemma \ref{stabilizer} 
${\bf g}\in G$.  \end{proof}\sepsub

\begin{sub}\label{AssPolar}{\sc\small The associated polarization}\rm

In \ref{first} and \ref{second} we had introduced polarizations for 
the different 
types of actions of $G_{red}$ on $W$ and of $\bg$ on $\bw$.  In the 
following we will describe polarizations on $W$ and $\bw$ which are 
compatible with the morphism $\zeta : W\lra \bw$.  Their weight 
vectors are related by the following matrix equations and determine 
each other.  The entries of the matrices are just the dimensions of 
the spaces $A_{ji}$ and $B_{ml}$.
 
In the sequel we will use the following {\bf notation}: the dimension 
of a vector space will be the small version of its name. So 
$m_{i}=\dim(M_{i}),\quad n_{l}=\dim(N_{l}),\quad p_{i}=\dim(P_{i}),\quad  
q_{m} = \dim(Q_{m})\quad a_{ji}=\dim(A_{ji}),\quad  b_{ml}=\dim(B_{ml})$ etc.

A proper polarization of the action of $G$ on $W$ is a tuple 
$\Lambda=(\lambda_1,\ldots,\lambda_r,-\mu_1,\ldots,-\mu_s)$, where $\lambda_i$ 
and $\mu_l$ are positive rational numbers such that 
$$\underset{1\leq i\leq r}{\sum}\lambda_im_i 
= \underset{1\leq l\leq s} {\sum}\mu_ln_l = 1.$$

We define the new sequence of rational numbers 
$\alpha_1,\ldots,\alpha_r,\beta_1,\ldots,\beta_s$ by the conditions
 
\[
\left(\begin{array}{c}
\lambda_1\\
\vdots\\
\vdots\\
\lambda_r
\end{array}\right)=
\left(\begin{array}{cccc}
1    & 0      & \cdots & 0\\
a_{21} & 1    & \ddots   &\vdots\\
\vdots & \ddots & \ddots & 0\\
a_{r1} & \cdots & a_{r,r-1} & 1
\end{array}\right)
\left(\begin{array}{c}
\alpha_1\\
\vdots\\
\vdots\\
\alpha_r
\end{array}\right),\quad
\left(\begin{array}{c}
\mu_1\\
\vdots\\
\vdots\\
\mu_s
\end{array}\right)=
\left(\begin{array}{cccc}
1    &  b_{2,1}     & \cdots & b_{s1}\\
0 & 1    & \ddots   & \vdots\\
\vdots & \ddots & \ddots & b_{s,s-1}\\
0 & \cdots & 0 & 1
\end{array}\right)
\left(\begin{array}{c}
\beta_1\\
\vdots\\
\vdots\\
\beta_s
\end{array}\right).
\]

Then we have 
$$1=\underset{1\leq i\leq r}{\sum}\lambda_im_i = \underset{1\leq i\leq r}
{\sum}\alpha_ip_i
\quad \text{and}\quad
1=\underset{1\leq l\leq s}{\sum}\mu_ln_l = \underset{1\leq l\leq s}
{\sum}\beta_lq_l.$$
In particular the tuple 
$\tilde{\Lambda}=(\alpha_{1},\cdots,\alpha_{r},-\beta_{1},\cdots,-\beta_{s})$ 
is a polarization on $\bw$ such that $\alpha_i$ is the weight of $P_i$ 
and $-\beta_l$ the weight of $Q_l$.  It is called the {\em associated 
polarization} on $\bw$.  It is compatible with $\zeta$ in the 
following sense: If $M'_i\subset M_i$, and $N'_l\subset N_l$ are 
linear subspaces, and if the subspaces of $P_{i}$ and $Q_{l}$ are 
defined by 
$$P'_i = \underset{i\leq j}{\oplus}M'_j\ot A_{ji},\quad 
\text{and}\quad Q'_l = \underset{l\leq m}{\oplus}N'_l\ot B_{ml}^*$$ 
respectively then we have $$\underset{1\leq i\leq 
r}{\sum}\lambda_im'_i = \underset{1\leq i\leq 
r}{\sum}\alpha_ip'_i,\quad \text{and}\quad \underset{1\leq l\leq 
s}{\sum}\mu_ln'_l = \underset{1\leq l\leq s}{\sum}\beta_lq'_l.$$

If the set of stable points in $\bf W$ with respect to the
associated polarization is non-empty then by \ref{second} the 
weights satisfy the conditions
$$\underset{i\leq j\leq r}{\sum}\alpha_jp_j>0\quad \text{for any\ }i
\and 
\underset{1\leq l\leq m}{\sum}\beta_lq_l>0\quad \text{for any\ }m.$$

Equivalently the conditions may also be written as $$\underset{i\leq j\leq r}
{\sum}\alpha_jp_j>0\quad \text{for\ } 2\leq i\leq r 
\and 1-\underset{m\leq l\leq s}{\sum}\beta_lq_l>0\quad 
\text{for\ } 2\leq m\leq s.$$

Substituting the weights of the original polarization on $W$, we can reformulate 
these conditions. In the cases treated in the examples they reduce to 
the following
\end{sub}\sepsub

\begin{subsub}\label{weightcond}Weight conditions.\rm

Let $W$ be of type $(r,s)$ and let 
$\Lambda=(\lambda_1,\ldots,\lambda_r,-\mu_1,\ldots,-\mu_s)$ be a proper 
polarization of $W$ with positive $\lambda_i$ and $\mu_l$.  If the set 
$\bws$ of stable points of $\bw$ with respect to the associated 
polarization $\tilde{\Lambda}$ is non-empty, then in case of

type $(2,1)$:\qquad    $\lambda_2-a_{21}\lambda_1>0$,

type $(3,1)$:\qquad 
$\lambda_3-a_{32}\lambda_2+(a_{32}a_{21}-a_{31})\lambda_1>0,\quad
\lambda_1(m_1+a_{21}m_2+a_{31}m_3)<1$,  

type $(2,2)$:\qquad $\lambda_2-a_{21}\lambda_1>0,\quad \mu_1-b_{21}\mu_2>0$.
\end{subsub}\sepsub

\begin{sub}\label{invpoly}{\sc\small Comparison of invariant polynomials}\rm

In the following we assume that 
$\tilde{\Lambda}=(\alpha_1,\ldots,\alpha_r,-\beta_1,\ldots,-\beta_s)$ is the 
polarization on $\bw$ associated to the polarization 
$\Lambda=(\lambda_1,\ldots,\lambda_r,-\mu_1,\ldots,-\mu_s)$.  The 
semi--stable locus $\bwss$ with respect to this polarization is more 
precisely defined by the character $\kx$ associated to it as in 
\ref{ResKing}.  If $q$ is lowest common denominator of 
$\alpha_1,\ldots,\alpha_r,\beta_1,\ldots,\beta_s$, we have $$\kx({\bf 
g}) = (\underset{1\leq i\leq r}{\prod}det({\bf g}_i)^{-q\alpha_i}) 
(\underset{1\leq l\leq s}{\prod}det({\bf h}_l)^{q\beta_l})$$ for an 
element ${\bf g}\in \bg$ with components ${\bf g}_{i}$ and ${\bf 
h}_{l}$.  By the matrix relations between the polarizations $q$ is 
also a common denominator of 
$\lambda_1,\ldots,\lambda_r,\mu_1,\ldots,\mu_s$, such that, if $p$ 
denotes the lowest, we have $q=pu$ for some $u$.  The character $\chi$ 
with respect to the given polarization can be defined by 
$$\chi(g,h)=\underset{1\leq i\leq r}{\prod}det(g_i)^{-p\lambda_i} 
\underset{1\leq l\leq s}{\prod}det(h_l)^{p\mu_l},$$ where the $g_i$ 
resp.  $h_l$ are the diagonal components of $g$ resp.  $h$, see 
\ref{group}.  Now the relations between the polarizations imply by a 
straightforward calculation that $$\kx(\theta(g,h))=\chi(g,h)^u.$$ If 
F is a $\kx^m$-invariant polynomial on $\bw$ it follows that 
$$F(\zeta((g,h).w))=F(\theta(g,h).\zeta(w))=\chi(g,h)^{um}F(\zeta(w)),$$ 
i.e.  that $F\circ\zeta$ is a $\chi^{um}$-invariant polynomial on $W$.  
As a consequence we obtain the 
\end{sub}\sepprop

\begin{subsub}\label{ssrelat}{\bf Lemma:}
$\zeta^{-1}(\bwss)\subset \wss$,

i.e.  if $w\in W$ and $\zeta(w)$ is $\bg$-semi-stable in $\bw$ with 
respect to the polarization 
$\tilde{\Lambda}=(\alpha_1,\ldots,\alpha_r,-\beta_1,\ldots,-\beta_s)$ 
then $w$ is $G$-semi-stable in $W$ with respect to the polarization 
$\Lambda=(\lambda_1,\ldots,\lambda_r,-\mu_1,\ldots,-\mu_s)$ (in the 
sense of \ref{defstab}).
\end{subsub}

\begin{proof}There exists a $\kx^m$-invariant polynomial F on $\bw$ 
such that $F(\zeta(w))\neq 0$.  Then 
$$F(\zeta((g,h).w))=F(\zeta(w))\neq 0$$ for any element $(g,h)$ in the 
unipotent subgroup $H\subset G$.  This means that w is $G$-semi-stable.  
\end{proof}

\sepprop

\begin{subsub}\label{remdet}{\bf Remark:}\rm
\, When we consider the subgroup ${\bf G'}\subset \bg$ defined by the 
condition $$\det({\bf g}_1) = \det({\bf h}_s) = 1,$$ we have 
$\theta(G')\subset{\bf G'}$ as follows from the definition of $G'$ in 
\ref{detgrp}.  With respect to these groups the semi-stable points are 
those over the semi-stable loci in $\P(W)$ resp.  $\P(\bw)$, with 
respect to the line bundles $$\kl = \ko_{\P(W)}(t)\qquad\text{and\ \ 
}{\bf L} = \ko_{\P(\bw)}({\bf t}),$$ where $t$ and ${\bf t}$ is 
defined as in \ref{ProjAct} in the different cases endowed with the 
modified action defined by the characters.  However, we cannot compare 
$\P(W)$ and $\P(\bw)$ directly because the morphism $\zeta$ does not 
descend.

We need the analogous statement of Lemma \ref{ssrelat} also in the case of
stable points. For that is is more convenient to use the subspace criterion
(1) of A. King in the case of $G_{red}$ and $\bg$. This gives also another
proof in the semi--stable case.
\end{subsub}

\sepprop

\begin{subsub}\label{srelat} {\bf Lemma:} With the same notation as in the
  previous Lemma
\[
\zeta^{-1}(\bws) \subset \ws
\]

\begin{proof}
Let $w=(\phi_{li})$ be a point of $W$ with maps $M_i\otimes
H_{li}^\ast\xrightarrow{\phi_{li}} N_l$ and suppose that $w$ is not
$G$--stable with respect to the polarization $\Lambda$. We can assume that it
is not $G_{red}$--stable, too. Then there are linear subspaces $M'_i \subset
M_i$ and $N'_l\subset N_l$ for all $i$ and $l$ such that the family $((M'_i)),
(N'_l))$ is proper and such that
\[
\phi_{li} (M'_i\otimes H^\ast_{li})\subset N'_l\quad \text{ and }\quad 
\sum\limits_i \lambda_i m'_i-\sum\limits_l \mu_l n'_l\geq 0.
\]
With these subspaces we can introduce the subspaces $P'_i\subset P_i$ and
$Q'_l\subset Q_l$ as
\[
P'_i=\underset{i\leq j}\oplus M'_j\otimes A_{ji}\quad \text{ and }\quad
Q'_l=\underset{m\leq l}\oplus N'_m\otimes B^\ast_{lm}.
\]
They form a proper family of subspaces and satisfy
\[
\xi_i(P'_i\otimes A_{i,i-1})\subset P'_{i-1}\quad ,\quad \gamma(w)(P'_1\otimes
H^\ast_{s1})\subset Q'_s\quad ,\quad \eta_l(Q'_{l+1}\otimes B_{l+1,l})
\subset Q'_l
\]
for the possible values of $i$ and $l$. But by the definition of the spaces
and because $\widetilde{\Lambda}$ is the associated polarization, the formulas
of \ref{AssPolar} imply the dimension formula
\[
\sum\limits_i \alpha_i p'_i-\sum\beta_l q'_l=\sum\limits_i \lambda_i
m'_i-\sum\limits_l \mu_l n'_l\geq 0.
\]
This states that also $\zeta(w)$ is not $\bg$--stable.
\end{proof}
\sepprop

\rm In section 7 we will derive sufficient conditions for the equality \[ 
\zeta^{-1} (\bws)= \ws\quad \text{ and }\quad \zeta^{-1} (\bwss)= \wss .
\]     In the following section we show how this equality implies 
the existence of a good and projective quotient $\wss//G$ using the 
result for $\bwss//\bg$ from Geometric Invariant Theory.
\end{subsub}
\sepsec
%\newpage
\section {Construction and properties of the quotient}
\rm

 We keep the notation of the previous sections and let 
$\widetilde{\Lambda}$ be the polarization on $\bw$ associated to the 
polarization $\Lambda$ on $W$.  We do not require that they are proper 
here, but we will do that later for the examples.  In addition we 
introduce the saturation
\[ 
Z=\bg\zeta(W)\subset \bw 
\] 
of the image of $W$ with respect to the action of $\bg$.  
\sepsub

\begin{sub}{\sc\small Construction of the quotient}\label{constquot}

\begin{subsub}\label{goodqu} {\bf Proposition}: Let $W$ and $\bw$ together
with their $G$-- and $\bg$--structure be as in section 2 and 5, let $\Lambda$
be a polarization for $(W,G)$ and $\widetilde{\Lambda}$ be the associated
polarization for $(\bw,\bg)$.  

(1) If\quad $\zeta^{-1}(\bws)= \ws$, then there exists a geometric 
quotient \quad $W^{s}(G,\Lambda)\xrightarrow{} M^s$
of $W^{s}$ by $G$, which is a quasi--projective nonsingular variety.

(2) If in addition \quad
$\zeta^{-1}(\bwss)= \wss$\quad  and \quad 
$(\bar{Z}\smallsetminus Z)\cap \bwss=\emptyset$, 
then there exists a good quotient \quad$W^{ss}(G,\Lambda)\xrightarrow{\pi} M$,
such that $M$ is a normal projective variety, $M^s$ is an open subset of $M$,
and $\ws\to M^s$ is the restriction of $\pi$.
 
\vskip5mm

\rm We recall here the definition of a good and a geometric 
quotient of C.S. Seshadri, see \cite{news}, \cite {mumf}.  Let an 
algebraic group $G$ act on an algebraic variety or algebraic scheme $X$.  
Then a pair $(\varphi, Y)$ of a variety and a morphism 
$X\xrightarrow{\varphi} Y$ is called a good quotient if
\begin{enumerate}
\item [(i)] $\varphi$ is $G$--equivariant (for the trivial action of $G$ on
  $Y$),

\item [(ii)] $\varphi$ is affine and surjective,
\item [(iii)] If $U$ is an open subset of $Y$ then $\varphi^\ast$ is an
  isomorphism $\ko_Y(U)\approx \ko_Y(\varphi^{-1} U)^G$, where the latter
  denotes the ring of $G$--invariant functions,
\item [(iv)] If $F_1, F_2$ are disjoint closed and $G$--invariant subvarieties
  of $X$ then $\varphi(F_1), \varphi(F_2)$ are closed and disjoint.
\end{enumerate}

If in addition the fibres of $\varphi$ are the orbits of the action and all
have the same dimension, the
quotient $(\varphi, Y)$ is called a geometric quotient. 

As usual we write $X//G$ for a good quotient space and $X/G$ for a geometric
quotient space.

\begin{proof} 
We will prove the second statement first, assuming
that the conditions of (1) and (2) are satisfied.
We use the abbreviations $W^{ss}= \wss, {\bf W^{ss}}=\bwss$ and similarly $W^s, 
{\bf W^s}$ for the subsets of the stable points.
By the result of A. King, \ref{ResKing}, there exists a good projective
quotient of $\bw^{ss}$ by the reductive group $\bg$.  So there exists 
also a good and projective quotient of the closed invariant subvariety 
$\bar{Z}\cap \bw^{ss}$ which we denote by 
\[ \bar{Z}\cap{\bf W}^{ss}\xrightarrow{\pi_0} {M}. \]

By assumption (2) $\bg\zeta(W^{ss})=Z\cap {\bf W}^{ss}=\bar{Z}\cap {\bf W}^{ss}$.
We let $\pi$ be the composition 
\[ W^{ss}\xrightarrow{\zeta} \bg\zeta 
(W^{ss})\xrightarrow{\pi_0} M.  
\] 
We know already that $M$ is projective. We will then verify that $(\pi, M)$ 
is the good quotient of the proposition.  We consider first the 
commutative diagram 
%\[
%\begin{CD}
%\bg\times W^{ss} @ >\mu>> \bg\zeta(W^{ss})\\
%@VpVV @ VV{\pi_0}V\\
%W^{ss} @ >\pi>> M
%\end{CD}\ \raisebox{-8mm}{,}
%\]
\[\xymatrix{
\bg\times W^{ss}\ar[r]^-{\mu}\ar[d]^p & \bg\zeta(W^{ss})\ar[d]^{\pi_0}\\
W^{ss}\ar[r]^-{\pi} & M
}\]

in which $p$ is the projection and $\mu$ is defined by $({\bf g}, w)\mapsto
{\bf g} \zeta (w)$. There is an action of $G$ on $\bg\times W^{ss}$ by
$g.({\bf g}, w)=({\bf g}\theta (g)^{-1}, g.w)$ and it follows that $\mu$ is
$G$--equivariant. 

\medskip

{\bf Claim}: {\em The morphism $\mu$ is a geometric quotient of $\bg\times
  W^{ss}$ by $G$.}

{\em Proof of the claim:} We show first that the fibres of $\mu$ are the
$G$--orbits. So let $({\bf g}, w)\ ,\ ({\bf g'} , w')$ be two elements in
$\bg \times W^{ss}$ such that $\mu({\bf g}, w)=\mu({\bf g'}, w')$. Then
$\zeta (w)={\bf g}^{-1} {\bf g'} \zeta (w')$.
By Lemma \ref{stabilizer} $g={\bf g}^{-1}{\bf g'}\in G$ and $g.({\bf g}, w)=
({\bf g'}, w')$. The claim will be proved if we show that $\mu$ has local
sections. For this it suffices to use the remark following Lemma
\ref{stabilizer} and a local section of the quotient map $\bg\to {\bf G}/G.$

Now we are going to verify the 4 properties of a good quotient for
$\pi$. Clearly (i) is satisfied by the definition of $\pi$. 

Proof of (ii). It is clear that $\pi$ is surjective. The morphism $\pi$ is 
affine because
$\pi=\pi_0\circ \zeta$ and $\pi_0$ and $\zeta$ are affine.

Proof of (iii). Let $U\subset M$ be an open subset. Then
\[
\ko(U)\subset \ko(\pi^{-1}(U))^G
\]
since $\pi$ is $G$--invariant. Conversely let $f\in \ko(\pi^{-1}(U))^G$. The
$f\circ p\in\ko(\bg\times \pi^{-1}(U))^G$, and since $\mu$ is a geometric
quotient, $f\circ p$ descends to an $\bar{f}\in\ko(\mu(\bg\times 
\pi^{-1}(U)))$, which is $\bg$--invariant.  Now again $\bar{f}$ 
descends because $\pi_0$ is a good quotient.  This proves equality 
$\ko(U)=\ko(\pi^{-1}(U))^G$.

Proof of (iv). Let $F_1, F_2$ be disjoint, closed, $G$--invariant subvarieties
of $W^{ss}$. Then $p^{-1}(F_1), p^{-1}(F_2)$ are disjoint, closed and ${\bf
  G}$--invariant subvarieties of $\bg\times W^{ss}$. Since $\mu$ is a
good quotient, $\mu(p^{-1}(F_1)), \mu(p^{-1}(F_2))$ are disjoint, closed and
$\bg$--invariant in $\bg\zeta(W^{ss})$.  Finally, since $\pi_0$ is a 
good quotient, $\pi_0\circ\mu(p^{-1}(F_1)), 
\pi_0\circ\mu(p^{-1}(F_2))$ are disjoint and closed subvarieties of 
$M$.  But $\pi_0\circ\mu(p^{-1}(F_i))=\pi(F_i)$, which proves (iv).

The normality of $M$ follows from the fact that $\bg\zeta(W^{ss})$ is
smooth and $\pi_0$ is a good quotient, \cite{mumf}, with respect to the
reductive group $\bg$. That $\pi$ becomes a geometric quotient on the open
set $W^s$ of stable points follows from the fact that the $\bg$--orbits in 
${\bf G}\zeta(W^s)=Z\cap {\bf W}^s$ intersect $W^s$ in $G$--orbits. In 
particular the stabilizers
of $w$ in $G$ and of $\zeta(w)$ in ${\mathbf G}$ are isomorphic, such that all
orbits have the same dimension.

The proof of (1) is a modification of the above. In any case $\pi_0$ induces the
geometric quotient  $\bar{Z}\cap {\bf W}^s\overset{\pi_0}{\longrightarrow}M_0$ 
with $M_0$ open in $M$. Now ${\bf G}\zeta(W^s)=Z\cap {\bf W}^s$ is a  
$\pi_0$-saturated open subset of $\bar{Z}\cap {\bf W}^s$, such that we obtain a
geometric quotient ${\bf G}\zeta(W^s)\overset{\pi_0}{\longrightarrow}M^s$
with $M^s\subset M_0$ open. By the same arguments as above applied to the
diagram related to $\bg\times W^s\to \bg\zeta(W^s)$ we conclude that 
$W^s\overset{\pi}{\longrightarrow}M^s$ is a geometric quotient.                                              
\end{proof}

{\em Remarks}: 1) The idea of this proof comes from \cite{se}, and
has already been used in \cite{dr_lp} and \cite{dr2}.

2) If the second condition of (2) is not satisfied, we cannot even prove that
$\wss$ admits a good quasi--projective quotient, because $Z\cap {\bf W}^{ss}$
might not be saturated. Of course the projectivity of the quotient depends on
this condition.
\end{subsub}

\end{sub}

\sepsub

\begin{sub}\label{s-equi}{\sc\small S--equivalence}\rm

We suppose that the hypotheses of proposition \ref{goodqu} are satisfied, with 
polarization $\Lambda$ for $(W,G)$ and associated polarization 
$\widetilde\Lambda$ for $(\bw,\bg)$.

It is easy to define the {\em Jordan-H\"older filtration} of 
$\bg$-semi-stable elements of $\bw$ with respect to $\widetilde\Lambda$
(cf.\cite{king} for a more general situation).  
Using the preceding 
results we can also define a {\em Jordan-H\"older filtration} of a 
$G$-semi-stable element of $W$ with respect to $\Lambda$.  Let 
$w=(\phi_{li})\in \wss.$   
Then there exist a positive integer $p$, an element $h\in H$ and filtrations 
$$M_i^0=\lbrace 0\rbrace\subset M_i^1\subset\cdots\subset M_i^p=M_i,\quad 
N_l^0=\lbrace 0\rbrace\subset N_l^1\subset\cdots\subset N_l^p=N_l,$$ 
with 
$$\sum_i\lambda_i\dim(M_i^j) = \sum_l\mu_l\dim(N_l^j)$$ 
for each $j$, such that $h.w=(\phi_{li})$ satisfies 
$$\phi_{li}(H^\ast_{li}\otimes M_i^j)\subset N_l^j,$$ 
and that if 
$$\phi_{li}^j : H_{li}^*\otimes(M_i^j/M_i^{j-1})\lra N_l^j/N_l^{j-1}$$ 
is the induced morphism, then $(\phi_{li}^j)_{li}$ is $G$-stable with respect to 
$\Lambda$ for any $j$.  This filtration and $h$ need not be unique, 
but $p$ is unique and the $(\phi_{li}^j)$, too, up to the order and 
isomorphisms.  Conversely, an element of $W$ having such a filtration 
is $G$-semi-stable with respect to $\Lambda$.  We say that two 
elements $(\phi_{li})$ and $(\phi'_{li})$ of $\wss$ are {\em 
S-equivalent} if they have Jordan-H\"older decompositions 
$(\phi_{li}^j),\ ({{\phi'}_{li}}^j)$ respectively of the same length, 
and if there exists a permutation $\sigma$ of $\lbrace 
1,\ldots,p\rbrace$ such that $({{\phi'}_{li}}^j)$ is isomorphic to 
$(\phi_{li}^{\sigma(j)})$ for any $j$.

The following result is also easily deduced from \ref{goodqu}.
\sepprop

\begin{subsub}\label{equiv}{\bf Proposition}: Let $w, w'\in \wss$.  Then 
$\pi(w) = \pi(w')$ if and only if $w$ and $w'$ are S-equivalent.
\end{subsub}

It follows that the set of closed points of $M$ is exactly the set of
S-equivalence classes of elements of $W^{ss}$. 
\end{sub}

\sepsec
%\newpage
\section{Comparison of semi--stability}

We are going to investigate conditions for the weights of the polarizations
under which a\\ (semi--)stable point $w\in W$ is mapped to a (semi--)stable point
$\zeta(w)\in\bw$. For the estimates we need the following constants which
depend on the dimensions $m_i$ and the composition maps $H_{li}\otimes
A_{i1}\to H_{l1}$.
\sepsub

\begin{sub}\label{constants} {\sc\small Constants} \rm 

Let $\kk$ be the family of proper linear subspaces 
\[ K\subset\underset{2\leq i}\oplus\ M_i\otimes A_{i1} \] 
such that $K$ is not contained in $\underset{2\leq i}\oplus 
M'_i\otimes A_{i1}$ for any family $(M'_i)\neq (M_i)$ of subspaces.  
For any $l$ we let the map 
\[ \underset{2\leq i}\oplus\ M_i\otimes 
A_{i1}\otimes H^\ast_{l1}\xrightarrow{\delta_l}\underset{2\leq 
i}\oplus\ M_i\otimes H^\ast_{li} \] 
be induced by the maps 
$A_{i1}\otimes H^\ast_{l1}\to H^\ast_{li}$ associated to the 
composition maps, which are supposed to be surjective, see \ref{abstr}.  
We introduce the constant 
\[ c_l(m_2, \ldots, m_r)=\underset{K\in\kk} 
\sup\ \rho_l (K)\quad \text{with}\quad \rho_l (K)=
\dfrac{\codim(\delta_l (K\otimes H^\ast_{l1}))}{\codim(K)}.  \]

Similarly we define the constants $d_i(n_1,\ldots, n_{s-1})$ in the dual
situation. Let
\[
\underset{l<s}\oplus\ N_l^\ast\otimes 
H_{li}^\ast\overset{\delta^\vee_i}\longleftarrow 
\underset{l<s}\oplus\ 
N_l^\ast\otimes B_{sl}\otimes H_{si}^\ast \] 
be induced by the maps 
$B_{sl}\otimes H_{si}^\ast\to H_{li}^\ast$  and let $\kl$ be the 
family of proper subspaces 
\[ L\subset \underset{l<s}\oplus N_l^\ast \otimes B_{sl} \] 
which are not contained in $\underset{l<s}\oplus N'_l\otimes B_{sl}$ for any 
family $(N'_l)\neq (N_l^\ast)$ of subspaces.  Then we define 
\[ d_i(n)=d_i(n_1, \ldots, n_{s-1})=\underset{L\in\kl}\sup\ 
\dfrac{\codim(\delta^\vee_i(L\otimes H^\ast_{si}))}{\codim(L)} .\]

\end{sub}
\sepprop

\begin{subsub}\label{monoton} {\bf Lemma:} If $m_i\leq \bar{m}_i$ for all
 $i\geq 2$, then $c_l(m_2,\ldots, m_r)\leq c_l(\bar{m}_2,\ldots,\bar{m}_r)$.

\begin{proof} It will be sufficient to assume that $m_i=\bar{m}_i$ for all $i$
except one, $m_2 <\bar{m}_2$ say. Then let $\bar{M}_i$ be vector spaces of 
dimensions $\bar{m}_i$ and suppose that 
\[ \bar{M}_2= L_2\oplus 
M_2\quad \text{ and }\quad \bar{M}_i=M_i\ \text{ for } i\geq 3. 
\] 
For any $K\in\kk$ we consider the subspace 
\[ 
\bar{K} =(L_2\otimes A_{21})\oplus K\quad\subset\quad (\bar{M}_2\otimes 
A_{21})\oplus (\underset{2<j}\oplus M_j\otimes A_{j1}).
\] 
Then $\codim(\bar{K})=\codim(K)$ and also $\codim(\delta_l(\bar{K}\otimes 
H^\ast_{l1}))=\codim(\delta_l(K\otimes H^\ast_{l1}))$ because $\delta_l$ 
is a direct sum of the surjective operator $A_{j1}\otimes 
H^\ast_{l1}\to H^\ast_{l1}$ such that $\delta_l(L_2\otimes 
A_{21}\otimes H^\ast_{l1})$ equals $L_{2}\ot H^\ast_{l2}$ and 
$\delta_l(\bar{K}\otimes H^\ast_{l1})=(L_2\otimes H^\ast_{l2})\oplus 
\delta_l (K\otimes H^\ast_{l1})$.  Therefore 
$\rho_l(K)=\rho_l(\bar{K})$.  Once we have shown that also $\bar{K}$ 
belongs to the analogous family $\bar{\kk}$, the Lemma is proved.  To 
see this let $\bar{M}'_2\subset \bar{M}_2$ and $\bar{M}'_i = 
M'_i\subset M_i$ for $i\geq 3$ be subspaces such that \[ 
\bar{K}\subset \underset{2\leq i}\oplus \bar{M}'_i\otimes A_{i,1}.  \] 
Then in particular $$L_2\otimes A_{21}\subset \bar{M}'_2\otimes 
A_{21}$$ and thus $L_2\subset \bar{M}'_2$.  But then 
$\bar{M}'_2=L_2\oplus M'_2$ with $M'_2=\bar{M}'_2\cap M_2$ and it 
follows that \[ K\subset \underset{2\leq i}\oplus\ M'_i\otimes A_{i1}.  
\] Since $K\in \kk$ we obtain $M'_i=M_i$ for all $i$ and then also 
$\bar{M}'_2=\bar{M}_2$.
\end{proof}
\end{subsub}
\sepsub

\begin{sub}\label{conv1}{\sc\small Study of the converse I}\rm

Let $\Lambda =(\lambda_1, \ldots \lambda_r, -\mu_1, \ldots, -\mu_s)$ be a
polarization on $W$ and let $\tilde{\Lambda}=(\alpha_1, \ldots, \alpha_r,\\
-\beta_1, \ldots, -\beta_s)$ be the associated polarization on $\bw$
(the associated polarization has been defined in \ref{AssPolar}). We
had shown in \ref{ssrelat} and \ref{srelat} that if $w\in W$ and $\zeta(w)$ is
(semi--)stable in $\bw$ with respect to $\bg$ and $\tilde{\Lambda}$, 
then so is $w$ with respect to $G$ and $\Lambda$.  We are going to 
derive sufficient conditions for the converse, i.e.  whether 
$\zeta(w)$ is (semi--)stable if $w$ is (semi--)stable.

In the sequel we are going to use the following {\bf notation}: Given a family
$M'=(M'_i)$ of subspaces $M'_i\subset M_i$ we set
\[
P_i(M')=\underset{i\leq j}\oplus\ M'_j\otimes A_{ji}
\]
and call a subspace $P'_i\subset P_i$ saturated if there is such a family
with $P'_i=P_i(M')$. Note that in this case $\sum\limits_i\alpha_i
p'_i=\sum\limits_i \lambda_i m'_i$. Similarly we introduce the spaces
$Q_l(N')$ for a subfamily $N'=(N'_l)$ of $(N_l)$ and call them saturated.

Let $w=(\phi_{li})$ be given and assume that $\zeta(w)$ is not semi--stable
with respect to $\widetilde{\Lambda}$. Then there exist linear subspaces $P'_i
\subset P_i$ and $Q'_l\subset Q_l$ such that 
\[
\xi_i(P'_i\otimes A_{i,i-1})\subset P'_{i-1},\quad \gamma(w)(P'_1\otimes
H^\ast_{s1})\subset Q'_s,\quad \eta_l(Q'_{l+1}\otimes B_{l+1,l})\subset Q'_l
\]
and such that
\[
\sum\limits_i \alpha_ip'_i-\sum\limits_l\beta_lq'_l>0,
\]
where as before the small characters denote the dimension of the spaces. If
there were subspaces $M'_i\subset M_i$ and $N'_l\subset N_l$ with
$P'_i=P_i(M')$ and $Q'_l=Q_l(N')$ as in \ref{srelat}, then 
$\gamma(w)(P'_1\otimes H^\ast_{s1})\subset Q'_s$ would imply that 
$\varphi_{li}(M'_i\otimes H^\ast_{li})\subset N'_l$ and we would have 
\[ \sum\limits_i\lambda_im'_i-\sum\limits_l\mu_l 
n'_l=\sum\limits_i\alpha_i p'_i-\sum\limits_l\beta_l q'_l>0 ,
\] 
and $w$ would not be semi--stable. In the following we are going to construct 
families $M'',\ N''$ of
subspaces $M''_i\subset M_i$ and $N''_l\subset N_l$ such that 
$P''_i=P_i(M'')$ and $Q''_l=Q_l(N'')$ are as 
close to $P'_i, Q'_l$ as possible and such that there is a useful 
estimate for 
\[ \sum\limits_i\lambda_i m''_i-\sum\limits_l \mu_l n''_l. \] 

{\em Step 1}:\hskip4mm \rm We can assume that $P'_i$ has a decomposition 
\[ P'_i=M'_i\oplus X_i \quad \text{in}\quad 
M_i\oplus (\underset{i<j}\oplus\ M_j\otimes A_{ji}) \] and such that 
$X_{r}=0$.  To derive this, we remark that for a subspace $S$ of a 
direct sum $E\oplus F$ of vector spaces there exists a linear map 
$E\xrightarrow{u} F$ such that the isomorphism $\binom{1\ 0}{u\ 1}$ of \break
$E\oplus F$ transforms $S$ into $S'\oplus S''$, where $S'$ is the 
projection of $S$ in $E$ and $S''=S\cap F$.  Using this and descending 
induction on $i$ we can find an element $h\in H_L\subset G_L$, see 
\ref{subgrp}, such that the truncations $\theta_{L,i}(h)\in \GL(P_i)$, 
see \ref{bigG}, map $P'_i$ onto a direct sum $M'_i\oplus X_i$ for any 
$i$.  Since $\xi_i(P'_i\otimes A_{i,i-1})\subset P'_{i-1}$ we easily 
derive that 
\[ 
\underset{i<j}\oplus\ M'_j\otimes A_{ji}\ \subset\   
X_i\ \subset\ \underset{i<j}\oplus\ M_j\otimes A_{ji} 
\] 
for all possible $i$.  We put \[ \rho_i=\codim(\underset{i<j}\oplus 
M'_j\otimes A_{ji}, X_i)=\codim(P_i(M'), P'_i).  \] Note that $\rho_r=0$.

{\em Step 2}: \hskip4mm Let $M''_1, \ldots M''_r$ be subspaces of $M_1, \ldots,
M_r$ respectively such that
\[
P_i(M'')\supset P'_i
\]
is minimal over $P'_i$ for any $i$.
Then $M'_i\subset M''_i$ since these spaces are the first 
components of $P'_i\subset P_i(M'')$ respectively and we have $M'_1=M''_1$.
We let 
\[ 
\sigma_i=\sum\limits_{i\leq j}(m''_j-m'_j)a_{ji}=\codim(P_i(M'), P_i 
(M'')).  \] 

{\em Step 3}:\hskip4mm We are going to define the subspaces $N'_l\subset 
N''_l\subset N_l$ as images.

Let $P_1\otimes H^\ast_{l1}\xrightarrow{\gamma_l(w)} N_l$ be the map which is
  the sum of the composed maps
\[
M_i\otimes A_{i1}\otimes H^\ast_{l1}\to M_i\otimes 
H^\ast_{li}\xrightarrow{\phi_{li}} N_l.  
\] 
Then we define 
\[ 
N'_l=\gamma_l(w)(P'_1\otimes H^\ast_{l1})=\phi_{l1}(M'_1\otimes 
H^\ast_{l1})+\gamma_l(w)(X_1\otimes H^\ast_{l1}) 
\] 
and 
\[ 
N''_l=\gamma_l(w)(P_1(M'')\otimes H^\ast_{l1})=\phi_{l1}(M''_1\otimes 
H^\ast_{l1})+\sum\limits_{2\leq j} \phi_{lj}(M''_j\otimes 
H^\ast_{lj}).  \] It follows $N'_l\subset N''_l$ for any $l$.

{\em Step 4}:\hskip4mm If the weights $\beta_l$ are supposed to be 
positive, we may assume that 
\[ \gamma(w)(P'_1\otimes 
H^\ast_{s1})=Q'_s\quad \text{ and }\quad \eta_l(Q'_{l+1}\otimes 
B_{l+1,l})=Q'_l 
\] 
for $l<s$.  Otherwise we could choose subspaces 
$\bar{Q}'_{l}\subset Q'_l$ by descending induction as images.  Then 
$-\sum\limits_l\beta_l \bar{q}'_l\geq -\sum \beta_lq'_l$  would improve 
the assumption on the choice of the spaces $P'_i$ and $Q'_l$.  Now it 
follows that for any $l$ 
\[ 
Q'_l\subset Q_l(N'')
\] 
because $P'_1\otimes H^\ast_{s1}$ is mapped 
to $\underset{l\leq s}\oplus N''_l\otimes B^\ast_{sl}$ and the maps 
$\eta_l$ are the identity on the spaces $N''_m$. Note that we even have
$Q'_l\subset Q_l(N')$ since $\gamma_l\ |\ P'_1\otimes H^\ast_{s1}$ factorises
through $\underset{l\leq s}\oplus N'_L\otimes B^\ast_{sl}$ as follows from the
  definition of $N'_l$.
\end{sub}
\sepprop

\begin{subsub}\label{estimate} {\bf Lemma}: Suppose that all $\beta_1,\ldots,
  \beta_s>0$, and let $\Delta=\sum\limits_i\lambda_i m''_i-\sum\limits_l\mu_l
  n''_l.$ Then
\[
\Delta>\sum\limits_l\beta_l q'_l-\sum\limits_l\mu_ln'_l+\sum\limits_i\alpha_i(\sigma_i-\rho_i)-\sum\limits_l\mu_l
c_l(m_2, \ldots, m_r)(\sigma_1-\rho_1).
\]
\end{subsub}

\begin{proof} Let 
\[
Y_l=\delta_l(X_1\otimes H^\ast_{l1})\ \subset\ Z_l=\underset{2\leq i}\oplus\ 
M''_i\otimes H^\ast_{li}.
\]
Since $X_1$ is not contained in a direct sum with spaces smaller than $M''_i$
we get 
\[
\codim(Y_l, Z_l)\leq c_l(m''_2, \ldots, m''_r) \codim(X_1, 
\underset{2\leq j}\oplus M''_j\otimes A_{j1}).
\] 
By Lemma \ref{monoton} and above definitions we get 
\[ \codim(Y_l, Z_l)\leq c_l(m''_2, \ldots, m''_r)
(\sum\limits_{1\leq i} m''_i a_{i1}-p'_1)=c_l(m_2, \ldots, m_r) 
(\sum\limits_{2\leq i}(m''_i-m'_i)a_{i1}-\rho_1).  
\] 
The map $\sum\limits_j \phi_{lj}$ sends $(M''_1\otimes H^\ast_{l1})\oplus Z_l$ 
onto $N''_l$ by definition of $N''_l$ and also maps $(M'_1\otimes 
H^\ast_{l1})\oplus \delta_l(X_1\otimes H^\ast_{l1})$ onto $N'_l$.  
Therefore, since $M'_1=M''_1$, we have a surjection 
\[ Z_l/Y_l\to N''_l/N'_l 
\] 
and the dimension estimate 
\[ n''_l-n'_l\leq c_l(m_2,\ldots, m_r)(\sum\limits_{2\leq i} 
(m''_i-m'_i)a_{i1}-\rho_1).  
\] 
Now we can derive the estimate of the Lemma. If there is no summation 
condition it is understood that the sum has to be taken over all 
indices of the given interval.  We have
 \[
\begin{array}{lcl}
\Delta & = & \sum\limits_i \lambda_i m''_i-\sum\limits_l \mu_ln''_l\\
       & = & \sum\limits_i \lambda_i m'_i-\sum\limits_l 
\mu_l n'_l+\sum\limits_j\lambda_j(m''_j-m'_j)-\sum\limits_l \mu_l(n''_l-n'_l).
\end{array}
\]
Substituting for $\lambda_j$ in the third sum and replacing the first by
\[
\sum\limits_i \lambda_i m'_i=\sum\limits_i\alpha_i \dim(\underset{i\leq j}
\oplus M'_j\otimes A_{ji})=\sum\limits_i\alpha_i(p'_i-\rho_i) 
\] 
and using the definition of $\sigma_i$ we get 
\[ 
  \Delta=\sum\limits_i\alpha_i p'_i-\sum\limits_l\mu_ln'_l 
  +\sum\limits_i\alpha_i(\sigma_i-\rho_i)-\sum\limits_l\mu_l(n''_l- 
  n'_l).  \] Now using the assumed estimate for the first sum and the 
  derived estimate for $n''_l-n'_l$ we get \[ 
  \Delta>\sum\limits_l\beta_l 
  q'_l-\sum\limits_l\mu_ln'_l+\sum\limits_i\alpha_i(\sigma_i-\rho_i) 
  -\sum\limits_l\mu_l c_l(m_2, \ldots, m_r)(\sigma_1-\rho_1).  \]
\end{proof}
\sepprop

\begin{subsub} \label{sone}{\bf Corollary}: Suppose that $s=1$, let 
$\Lambda=(\lambda_1, \ldots,\lambda_r, -\frac{1}{n_1})$ and let 
$\tilde{\Lambda}$ be the associated polarization $(\alpha_1,\ldots, 
\alpha_r, -\frac{1}{n_1}).$ If all $\alpha_i >0$ and if 
\[ 
\lambda_2\geq \dfrac{a_{21}}{n_1} c_1(m_2, \ldots, m_r) 
\] 
then 
\[
\zeta^{-1}\bwss = \wss\quad \text{and}\quad\quad 
\zeta^{-1}\bws=\ws.
\]
\end{subsub}
\sepprop

{\em Remarks}: (1) Note that by the normalization of the polarizations we must
have $\mu_1 n_1=1$ such that $1/n_1$ is the only possible value for
$\mu_1=\beta_1$.

(2) If all $\alpha_i>0$, then the necessary conditions for $\ws\neq \emptyset$
and $\bws\neq\emptyset$ are both satisfied, see \ref{AssPolar}. The condition
of the corollary is an extra condition.
\vskip5mm

\begin{proof}
Let us first assume that $\zeta(w)$ is not semi--stable and let the
spaces $P'_i$ and $Q'_1$ be as at the beginning of \ref{conv1}.  The only 
$\beta_1=1/n_1$ is positive.  Let the other spaces be chosen as in 
\ref{conv1}.  The difference $\sum\beta_lq'_l-\sum \mu_l n'_l$ reduces 
to $q'_1/n_1-n'_1/n_1$, and since $N'_1=\gamma(w)(P'_1\otimes 
H^\ast_{11})=Q'_1$, this difference is zero.  Therefore 
\[ 
\Delta>\sum\limits_i \alpha_i(\sigma_i-\rho_i)-\dfrac{1}{n_1} c_1(m_2, 
\ldots, m_r)(\sigma_1-\rho_1).
 \] 
Since all the $\alpha_i$ are positive we have 
\[
\sum\limits_i \alpha_i(\sigma_i-\rho_i)\geq 
\alpha_1 (\sigma_1-\rho_1)+\alpha_2(\sigma_2-\rho_2).  
\] 
Moreover, $\xi_2$ induces a surjection 
\[
P_2(M'')\otimes A_{21}/P'_2\otimes A_{21}\to 
P_1(M'')/P'_1 
\] 
because $M'_1=M''_1$.  Therefore we obtain the dimensions estimate 
$(\sigma_2-\rho_2)a_{21} \geq \sigma_1-\rho_1$.  It follows that 
\[ 
\Delta>(-\dfrac{1}{n_1} c_1(m_2, \ldots, 
m_r)+\alpha_1+\dfrac{\alpha_2}{a_{21}})(\sigma_1-\rho_1).  \] Since 
$\lambda_2=a_{21}\alpha_1+\alpha_2\geq \frac{a_{21}}{n_1}c_1(m_2, 
\ldots, m_r)$ the last expression is non--negative.  This proves the 
case of semi--stability.  For the case of stability we assume that $w$ is
stable and that $\zeta(w)$ is already semi--stable. If $\zeta(w)$ were not 
stable, we would find subspaces $P'_i$ and $N'_1$ as in \ref{conv1} such that 
$\sum
\alpha_ip'_i-\mu_1n'_1=0$ and such that at least one $P'_i$ is different from
$P_i$. Now let the spaces $M''_i$ and $N''_I$ be constructed as above. Then we
have
\[
\Delta\geq \sum\limits_i \alpha_is_i-\frac{c_1}{n_1} s_1\geq \sum\limits_{2<i}
\alpha_i s_i+(\lambda_2-\frac{c_1}{n_1} a_{21})\frac{s_1}{a_{21}}\geq 0,
\]
where $s_i=\sigma_i-\rho_i=\dim(P_i(M'')/P'_i)$, and where we use that $s_2
a_{21}\geq s_1$. If the family $M''$ is different from $M$, then $0>\Delta$,
and if it is equal, then $\Delta=0$. In order to obtain a contradiction we
have to show that $M''$ is different from $M$. Assume that it is not. Then
$s_i=\dim(P_i/P'_i)$ and we must have $s_i=0$ for $i\geq 3$ and
$s_1(\lambda_2-\frac{c_1}{n_1} a_{21})=0$. If also $s_1=0$, then by the above
  estimate also $s_2=0$, contradicting the choice of the $P'_i$. Therefore
  $s_1\neq 0$ and $\lambda_2=\frac{c_1}{n_1} a_{21}.$ But then
  $\Delta=\alpha_2(s_2-\frac{s_1}{a_{21}})$ and we have $s_2a_{21}=s_1$. From
  this it is easy to see that $P'_i=P_i(\tilde{M})$ where
  $\tilde{M}_i=M_i$ for $i\neq 2$ and $\tilde{M}_2=M'_2\neq M_2$. Then
  we have
\[
\sum\limits_i\alpha_i\widetilde{m}_i-\mu_1 n'_1=
\sum\limits_i\lambda_i p'_i-\mu_1n'_1=0
\]
which contradicts the stability of $w$.
\end{proof} 

\sepsub

\begin{sub}\label{conv2} {\sc\small Study of the converse II}\rm

We keep the notation of \ref{conv1} and compare the (semi--)stability of
points in $W$ and ${\mathbf W}$ in two steps, each reducing to the case
$s=1$.  We consider the intermediate space
\[{\mathbf V}={\mathbf W}_{L}\oplus\ \underset{1\leq l\leq s}\oplus\ 
\Hom(P_1\otimes H_{l1}^\ast, N_l) \] and the maps 
\[W\xrightarrow{\zeta_1}{\mathbf V}\xrightarrow{\zeta_2}{\mathbf W}.\] 
Here $\zeta_1$ is defined by 
\[ w\mapsto (\xi_2, \ldots, \xi_r, 
\gamma_1(w), \ldots, \gamma_s(w)), \] 
where $\gamma_l(w)$ is the map 
defined by $w=(\phi_{li})$ as in \ref{conv1}.  The map $\zeta_2$ is 
defined by \[ (x_2, \ldots, x_r, \gamma_1, \ldots, \gamma_s)\mapsto 
(x_2, \ldots, x_r, \gamma, \eta_1, \ldots, \eta_{s-1}), \] 
where now 
$\gamma: P_1\otimes H_{s1}^\ast\to Q_s$ is induced by the tuple 
$(\gamma_1, \ldots, \gamma_s)$ as the sum of the compositions \[ 
P_1\otimes H_{s1}^\ast\to N_l\otimes H_{l1}\otimes H_{s1}^\ast\to 
N_l\otimes B_{sl}^\ast \] which are induced by the $\gamma_l$ and the 
pairings $B_{sl}\otimes H_{l1}\to H_{s1}$.  It is obvious that \[ 
\zeta=\zeta_2\circ \zeta_1.  \] Note that both $\zeta_1$ and $\zeta_2$ 
are injective by the same reason as for $\zeta$.

On ${\mathbf V}$ the group ${\mathbf G}_L \times G_R$ acts naturally and we
have the embedding 
\[
G=G_L\times G_R\overset{\theta_L\times id}\hookrightarrow {\mathbf 
G}_L\times G_R, 
\] see \ref{bigG}.  It follows as in section 5 that 
$\zeta_1$ is compatible with the group actions and that $w, w'\in W$ 
are on the same $G$--orbit if and only if $\zeta_1(w), \zeta_1(w')$ 
are on the same ${\mathbf G}_L\times G_R$ orbit.  Similarly we have 
the group embedding ${\mathbf G}_L\times G_R\hookrightarrow {\mathbf 
G}_L\times {\mathbf G}_R={\mathbf G}$ and $\zeta_2$ is equivariant and 
satisfies the analogous statements for the orbits.  Given the 
polarization $\Lambda=(\lambda_1, \ldots, \lambda_r, -\mu_1, \ldots, 
-\mu_s)$ for $(W, G)$ we consider the polarization 
$\overline{\Lambda}=(\alpha_1, \ldots, \alpha_r, -\mu_1, \ldots, 
-\mu_s)$ for $({\mathbf V}, {\mathbf G}_L \times G_R)$ where the 
$\alpha_i$ are defined as in \ref{AssPolar}.  As in \ref{ssrelat}, 
\ref{srelat} it is easy to show that 
\[ \zeta_1^{-1}{\mathbf 
V}^{ss}({\mathbf G}_L\times G_R, \overline{\Lambda})\subset \wss \quad 
\text{ and }\quad \zeta_1^{-1}{\mathbf V}^s({\mathbf G}_L\times G_R, 
\overline{\Lambda})\subset \ws \] 
and similarly that 
\[ 
\zeta_2^{-1}\bwss \subset {\mathbf V}^{ss}({\mathbf G}_L\times G_R, 
\overline{\Lambda})\quad \text{ and }\quad \zeta_2^{-1} \bws 
\subset{\mathbf V}^{s}({\mathbf G}_L\times G_R, \overline{\Lambda}) .\] 
Note that as for $W^{ss}, W^s$, we have unipotent sub-orbits in 
${\mathbf V}^{ss}$ and ${\mathbf V}^s$, see \ref{defstab}.  We are 
going to show that in all 4 cases equality holds under suitable 
conditions on the weights of the polarizations.  Then the same is true 
for $\zeta$.
\end{sub}
\sepsub

\begin{sub}\label{estim1}{\sc\small  Estimate for \rm $\boldsymbol \zeta_1$}\rm

Let $w=(\phi_{li})$ in $W$ be given and assume that $\zeta_1(w)$ is not
semi--stable. Then there are linear subspaces $P'_i\subset P_i$ and
$N'_l\subset N_l$ and a unipotent element $h\in H_R$ such that for
$(\gamma'_1,\ldots,\gamma'_s)=h_{.}(\gamma_1, \ldots, \gamma_s)$ we have 
\[
\xi_i(P'_i\otimes A_{i,i-1})\subset P'_{i-1}\quad \text{ and }\quad 
\gamma'_l(P'_l\otimes H_{l1}^\ast)\subset N'_l \] for all $i\geq 2$ 
and all $l$, and such that \[ \sum\limits_i\alpha_ip'_i-\sum\limits_l 
\mu_l n'_l>0.  \] 
We may assume that $h=id$ because $H_R$ acts on $W$ 
in the same way and we can replace $w$ by $h.w$.  Moreover, we may assume 
that all $N'_l$ are equal to $\gamma_l(P'_1\otimes H_{l1}^\ast)$ since 
all $\mu_l>0$.  Now we proceed as in \ref{conv1} replacing the spaces 
$Q_l$ by $N_l$.  Therefore we find subspaces $M'_i \subset 
M''_i\subset M_i$ such that $M'_1=M''_1$ and such that 
\[ P'_i=M'_i \oplus X_i\ ,\quad P_i(M')\subset P'_i\subset P_i(M'') \] 
and the family $M''$ is 
minimal with this property.  We denote 
\[ \rho_i= 
\codim(P_i(M'),P'_i)\quad,\quad \sigma_i=\codim(P_i(M'), P_i(M'')) \] 
and let 
\[ N''_l=\gamma_l(P_1(M'')\otimes H_{l1}^\ast)\supset N'_l.  \] 
As 
in \ref{estimate} we consider the surjection \[ Z_l/Y_l\to N''_l/N'_l, 
\] where $Y_l\subset Z_l$ are the same, and we get the estimate 
\[ 
n''_l-n'_l\leq c_l(m_2,\ldots,m_r)(\sigma_1-\rho_1) \] for any $l$.  
The estimation of the discriminant $\Delta$ is now simpler than in 
\ref{conv1}.
\end{sub}
\sepprop 

\begin{subsub}\label{discrim1} {\bf Lemma}: With the above 
notation 
\[ \Delta:=\sum\limits_i \lambda_im''_i-\sum\limits_l\mu_l 
n''_l\ >\ \sum\limits_i\alpha_i(\sigma_i-\rho_i)-\sum\limits_l \mu_l 
c_l(m)(\sigma_1-\rho_1) \] 
where $c_l(m)=c_l(m_2, \ldots, m_r)$.

\begin{proof}By replacing dimensions and inserting the estimate for
  $n''_l-n'_l$ as in \ref{conv1} we get
\[
\begin{array}{lcl}
\Delta & = & \sum\limits_i\alpha_ip'_i-\sum\limits_l\mu_l
n'_l+\sum\limits_i\alpha_i(\sigma_i-\rho_i)-\sum\limits_l\mu_l (n''_l-n'_l)\\
       & > & \sum\limits_i\alpha_i(\sigma_i-\rho_i)-\sum\limits_l\mu_l
       c_l(m)(\sigma_1-\rho_1).
\end{array}
\]
\end{proof}
\end{subsub}
\sepprop

 \begin{subsub}\label{equal1} {\bf Corollary}: Let 
$\Lambda=(\lambda_1, \ldots,\lambda_r, -\mu_1,\ldots, -\mu_s)$ be a 
polarization for $W$ and let $\bar{\Lambda} =(\alpha_1, \ldots, 
\alpha_r, -\mu_1,\ldots, -\mu_s)$ be the associated polarization for 
${\mathbf V}$ as in \ref{conv2}.  If all $\alpha_i > 0$ and 
\[ 
\lambda_2\geq a_{21}\sum\limits_l\mu_lc_l(m) \] then 
\[ 
\zeta_1^{-1}\bvss = \wss\quad \text{ and }\quad \zeta_1^{-1}\bvs=\ws.  
\]

\begin{proof} The proof is the same as for \ref{sone}, because the spaces $P'_i$
  and $P_i(M'')$ are defined in the same way and we thus get the estimate
  $(\sigma_2-\rho_2) a_{21}\geq \sigma_1-\rho_1$.
\end{proof}
\end{subsub}
\sepsub

\begin{sub} \label{estim2} 
{\sc\small Estimate for ${\boldsymbol \zeta_2}$}\rm

The analogous estimate for $\zeta_2$ follows by duality while we can assume
that $s=1$ or $r=1$. The proof could be done by formally transform it into a
dual situation which is similar to that of \ref{estim1}, but it is better to
keep direct track of the weights. Let $(x_2, \ldots, x_r, \gamma_1, \ldots
\gamma_s)$ be given in ${\mathbf W}_L\oplus V$ and assume that its image under
$\zeta_2$ is not semi--stable. Then there are subspaces $P'_i\subset P_i$ and
$Q'_l\subset Q_l$ such that
\[
x_i(P'_i\otimes A_{i,i-1})\subset P'_{i-1}\ ,\quad \gamma(P'_i\otimes 
H_{s1}^\ast)\subset Q'_s\ ,\quad \eta_l(Q'_{l+1}\otimes 
B_{l+1,l})\subset Q'_l, \] where $\gamma$ is defined as in 
\ref{conv2}, and such that 
\[ 
\sum\limits_i\alpha_ip'_i-\sum\limits_l\beta_lq'_l>0.  \] We assume 
that all $\alpha_i\geq 0$, and then we may assume that $P'_i$ is 
maximal, i.e.  the inverse image of $P'_{i-1}\otimes A_{i,i-1}^\ast$ 
under $P_i\to P_{i-1}\otimes A_{i,i-1}^\ast$ for $i\geq 2$, and 
similarly $P'_1$ in $P_1$ under $P_1\to Q_s\otimes H_{s1}.$ As in 
\ref{estim1} we can find subspaces $N'_l\subset N_l$ such that 
\[ 
Q'_l=N'_l\oplus X'_l \quad \text{ and hence } 
(Q_l/Q'_l)^\ast=(N_l/N'_l)^\ast\oplus X_l.  \] We choose subspaces 
$N''_l\subset N'_l$ which are maximal such that \[ Q_l(N'')\subset 
Q'_l\subset Q_l(N').  \] We have $N''_s = N'_s$.  We let $P''_1$ be 
the inverse image of $Q_s(N'')$ under $P_1\to Q_s\otimes H_{s1}$.  
Then $P''_1\subset P'_1$.  Furthermore we let inductively 
$P''_i\subset P'_i$ be the inverse images for $i\geq 2$.  Then we have 
injections 
\[ (P'_i/P''_i)\otimes A_{i,i-1}\to P'_{i-1}/P''_{i-1}
\] 
and induced by factorization the images 
\[
P'_i/P''_i\otimes A_{i, i-1}\otimes\ldots\otimes
A_{21}\twoheadrightarrow(P'_i/P''_i)\otimes A_{i1} \rightarrowtail P'_1/P''_1.
\]
The induced injections 
\[
P'_i/P''_i\rightarrowtail(P'_1/P''_1)\otimes A^\ast_{i1}
\]
imply the dimension estimates
\[ p'_i-p''_i\leq a_{i1} (p'_1-p''_1) 
\] 
for $i\geq 2$.  Next we consider the 
homomorphism \[ Z_1=\underset{l<s}\oplus\ (N_l/N''_l)^\ast\otimes 
H_{l1}^\ast\xleftarrow{\delta_1^\vee}\underset{l<s}\oplus\ 
(N_l/N''_l)^\ast\otimes B_{sl}\otimes H_{s1}^\ast.  \] We have 
$X_s\subset \underset{l<s}\bigoplus(N_l/N''_l)^\ast \otimes B_{sl}$ 
and consider the subspace \[ Y_1=\delta_1^\vee(X_s\otimes 
H_{s1}^\ast)\subset Z_1.  \] By the definition of the constant 
$d_1(n)=d_1(n_1, \ldots, n_{s-1})$ we get \[ \dim(Z_1/Y_1)\leq d_1(n) 
\codim(X_s)=d_1(n)(\sigma_s-\rho_s) \] where \[ 
\sigma_l=\codim(Q_l^\ast(N/N'), Q_l^\ast(N/N''))\quad \text{ and } 
\quad \rho_l=\codim(Q_l^\ast(N/N'), (Q_l/Q'_l)^\ast).  \] Further we 
have a surjective map \[ Z_1/Y_1\to (P_1/P''_1)^\ast/(P_1/P'_1)^\ast 
\] which is induced by the map $Q_s^\ast\otimes H_{s1}^\ast\to 
P_1^\ast$ and the induced surjection $Q_s^\ast(N/N'')\otimes 
H_{s1}^\ast\to (P/P''_1)^\ast$, since $N''_s=N'_s$.  So we get \[ p'_1 
-p''_1\leq d_1(n)(\sigma_s -\rho_s).  \] 
Now we can estimate the discriminant in
\sepprop

\begin{subsub}\label{discrim2} {\bf Lemma}: Let all the 
$\alpha_{i}$ be non-negative and let 
$\Delta:=\sum\limits_i\alpha_ip''_i-\sum\limits_l \mu_l n''_l.$
Then 
\[
\Delta >\ \sum\limits_l\beta_{l}(\sigma_{l} - \rho_{l}) - 
\sum\limits_i\alpha_i a_{i1} d_1(n)(\sigma_{s} - \rho_{s}).
\]

\begin{proof} Since $\sum\limits_i\alpha_ip_i=\sum\limits_l \mu_l n_l$ 
we also have 
$$\Delta =\sum\limits_l\mu_l(n_l-n''_l)-\sum\limits_i\alpha_i(p_i-p''_i).$$ 
with the same steps as in the previous proofs we get $$\Delta = 
\sum\limits_i\alpha_ip'_i-\sum\limits_l \beta_l q'_l + \sum\limits_l 
\beta_l(\sigma_{l} - \rho_{l})-\sum\limits_i\alpha_i(p'_i-p''_i)$$ 
Inserting the assumption on the first difference and the estimate for 
$p'_i-p''_i$ we get the result.  
\end{proof}\vskip3mm\rm As in the previous cases we obtain the
\end{subsub} 
\sepprop

\begin{subsub}\label{equal2} {\bf Corollary}: In the above notation let all
$\alpha_{i} > 0$,\ and all $\beta_{l} > 0$, 
and let 
\[
\mu_{s-1}\geq b_{s,s-1}d_1(n)\sum\limits_i \alpha_i a_{i1}.
\]
Then 
\[
\zeta_{2}^{-1}\bwss = \bvss \quad \text{and}\quad 
\zeta_{2}^{-1}\bws = \bvs.
\]

\begin{proof}: In the notation of \ref{estim2} there is a surjection
$(Q_{s-1}(N')/Q'_{s-1})^\ast\otimes B_{s, s-1}\to (Q_s(N')/Q'_s)^\ast$ because
$N''_s=N'_s$. Therefore $(\sigma_{s-1}-\rho_{s-1}) b_{s,s-1}\geq
\sigma_s-\rho_s$. If the condition of the Corollary is satisfied, then
$\Delta>0$ follows, where we use $\mu_{s-1}=\beta_sb_{s,s-1}+\beta_{s-1}.$
\end{proof}\end{subsub}
\rm Combining the results of \ref{equal1} and \ref{equal2} we get the
\end{sub}
\sepprop

\begin{subsub}\label{equality} {\bf Proposition}: Let
$\Lambda=(\lambda_{1},\ldots,\lambda_{r},-\mu_{1},\ldots,-\mu_{s})$ be a
polarization for $(W,G)$ and let 
$\widetilde{\Lambda}=(\alpha_{1},\ldots,\alpha_{r},-\beta_{1}, 
\ldots,-\beta_{s})$ be the associated polarization for $(\bw,\bg)$.  
Suppose that all $\alpha_{i} > 0$, and all $\beta_{l} > 0$ and that
\[ 
\lambda_2\geq a_{21}\sum\limits_l\mu_lc_l(m)\quad \text{ and }\quad
\mu_{s-1}\geq b_{s,s-1}d_1(n)\sum\limits_i \alpha_i a_{i1}.
\]
Then 
$$\zeta^{-1}\bwss = \wss\quad \text{and}\quad \zeta^{-1}\bws = \ws.$$
\end{subsub}

\sepsec
%\newpage
\section {Projectivity conditions}

The projectivity of the quotient in \ref{goodqu} depends on the second 
condition in (2),
i.e. whether the boundary $\bar{Z}\smallsetminus Z$ of the saturated set
contains no semi--stable points of ${\bf W}.$ Again this condition depends on
the chosen polarization and conditions for the weights. In order to derive
these conditions in some cases we describe the boundary in terms independent
of the group action.
\sepsub

\begin{sub}\label{satbd}{\sc\small   Saturated boundary.}

\rm The elements of ${\bf W}$ are tuples ${\bf w}=(x_2, \ldots,
x_r, \gamma, y_1, \ldots, y_{s-1})$ of linear maps 
\[
P_i\otimes A_{i,i-1}\xrightarrow{x_i}P_{i-1},\quad P_1\otimes
H^\ast_{s1}\xrightarrow{\gamma} Q_s,\quad Q_{l+1}\xrightarrow{y_l}
B^\ast_{l+1,l}\otimes Q_l
\]
If ${\bf w}\in Z$, there are an element $w\in W$ and automorphisms $\rho_i\in
\Aut(P_i)\ ,\ \sigma_l\in \Aut(Q_l)$ such that
\[
x_i=\rho_{i-1}\circ \xi_i\circ(\rho_i^{-1}\otimes id),\quad
\gamma=\sigma_1\circ \gamma(w)\circ (\rho_1^{-1}\otimes id),\quad 
y_l=(id\otimes\sigma_l)\circ \eta_l\circ \sigma^{-1}_{l+1}.
\]
Here $id$ stands for the different identities of the spaces $A, B\ \text{and}\ H.$ 
We let $\widetilde{x_i}$ respectively $\widetilde{\xi_i}$ be the
mapping
\[
P_i\otimes A_{i, i-1}\otimes \ldots\otimes A_{21}\to P_{i-1}\otimes
A_{i-1,i-2}\otimes\ldots\otimes A_{21}
\]
induced by $x_i$ respectively $\xi_i$ for $i\geq 3$. From the relations
between the $x_i$ and $\xi_i$ it follows easily that for each $i\geq 3$ the
composition $x_2\circ\widetilde{x_3}\circ\ldots\circ \widetilde{x_i}$ has a
factorization
\[
\xymatrix{
P_i\otimes A_{i,i-1}\otimes\cdots\otimes A_{21} \ar@{>>}[d] \ar[r] & P_1\\
P_i\otimes A_{i1}\ar[ur]_{x_{1i}} }
\]
where the vertical map is the surjection induced by the pairings. This follows
from the commutative diagrams induced by the automorphism $\rho_i$ and  because
$\xi_2\circ \widetilde{\xi_3}\circ\cdots\circ \widetilde{\xi_i}$ admits such a
factorization for each $i\geq 3$. We put $x_{21}=x_2$. By the dual description
for the maps $y_l$ we are given factorizations 
\[
\xymatrix{
 & B^\ast_{sl}\otimes Q_l \ar@{>->}[d]\\
Q_s \ar[ur]^{y_{ls}} \ar[r] & B^\ast_{s,s-1}\otimes \ldots\otimes 
B^\ast_{sl}\otimes Q_l
}
\]
 of the maps $\widetilde{y}_l\circ\ldots\circ\widetilde{y}_{s-2}\circ y_{s-1}$
 for $l\leq s-2$. By similar arguments there are also factorizations
\[\hspace{25mm}
\xymatrix
{P_i\otimes A_{i1}\otimes H^\ast_{s1} \ar@{>>}[d] \ar[r]^-{\scriptscriptstyle
    {x_{1i}\otimes id}} & P_1\otimes H^\ast_{s1}\ar[r]^{\gamma} & Q_s\\
P_i\otimes H^\ast_{si}\ar[urr]_{\gamma_{si}}
}\hspace{4cm}\raisebox{-7mm}{$(L_i)$}
\] 
for all $i\geq 2$ and dually factorizations
\[\hspace{25mm}
\xymatrix{ & & Q_l\otimes H_{l1} \ar@{>->}[d]\\ 
  P_1\ar[urr]^{\gamma_{l1}}\ar[r]_-{\gamma} & Q_s\otimes
  H_{s1}\ar[r]_-{\scriptscriptstyle{y_{ls}\otimes id}} & B^\ast_{sl}\otimes
  Q_l\otimes H_{s1} } \hspace{4cm}\raisebox{-7mm}{$(R_l)$}
\]
for all $l$. Moreover, there are further factorizations of the induced
composed maps
\[\hspace{25mm}
\xymatrix{ 
P_i\otimes H^\ast_{si}\otimes B_{sl} \ar@{->>}[d]
\ar[r]^-{\scriptstyle{\gamma_{si}\otimes id}} & Q_s\otimes B_{sl}
  \ar[r]^-{\widetilde{y}_{ls}} & Q_l\\
P_i\otimes H^\ast_{li}\ar[urr]_{\Phi_{li}}
}\hspace{4cm}\raisebox{-7mm}{$(L_{li})$}
\]
and dually
\[\hspace{25mm}
\xymatrix{
& & Q_l\otimes H_{li}\ar@{>->}[d]\\ 
P_i \ar[urr]^{\Psi_{li}} \ar[r]_-{\widetilde{x}_{1i}} & P_1\otimes A^\ast_{i1}
\ar[r]_-{\scriptstyle{\gamma_{l1}\otimes id}} & Q_l\otimes H_{l1}\otimes A^\ast_{i1}
}\hspace{4cm}\raisebox{-7mm}{$(R_{li})$} 
\]
All these factorizations are based on mappings induced by the pairings. All
factorization conditions are independent of the chosen automorphisms. One can
rediscover the original components $\phi_{li}$ of $w$ from $\Phi_{li}$ or
$\Psi_{li}$ if $x_j=\xi_j$ and $y_l=\eta_l$ for all $j$ and all $l$. In fact
we have 
\end{sub}
\sepprop

\begin{subsub}\label{descrZ} {\bf Lemma}: Let ${\bf w}=(x_2, \ldots, x_r,
  \gamma, y_1, \ldots, y_{s-1}) \in {\bf W}$. Then ${\bf w} \in Z$ if and only if 
  \begin{enumerate}\setlength{\itemsep}{+0.5ex}
  \item [(1)] rank $x_i = \sum_{i\leq j} m_j a_{j, i-1}\quad$ for $i\geq 2$
  \item [(1$^\ast$)] rank $y_l=\sum_{k\leq  l} b_{l+1,k}n_k$\quad for $l\leq s-1$
  \item [(2)] $x_2\circ\widetilde{x_3}\circ\ldots\circ \widetilde{x_i}$ has a
    factorization $P_i\otimes A_{i1}\xrightarrow{x_{i1}} P_1$ for $i\geq 3$
  \item [(2$^\ast$)] $\widetilde{y}_l\circ\ldots \circ\widetilde{y}_{s-2}\circ
  y_{s-1}$ has a factorization
$Q_s\xrightarrow{y_{ls}} B^\ast_{sl}\otimes Q_l\text{ for } l\leq s-2$
\item [(3)] $\gamma\circ (x_{1i}\otimes id)$ has factorizations $(L_i)$ and
    $(L_{li})$ 
  \item [(3$^\ast$)] $(y_{ls}\otimes id)\circ\gamma$ has factorizations $(R_l)$
    and $(R_{li})$.
  \end{enumerate}

\begin{proof} If ${\bf w}\in Z$, the three conditions are satisfied by the
  above, where $rank\ x_i$ is the dimension of the image of $\xi_i$ and $rank\
  y_l$ is the $rank$ of $\eta_l$ as the map $Q_{l+1}\to B^\ast_{l+1,l}\otimes
  Q_l$. Let   conversely ${\bf w}$ satisfy these conditions. We proceed by
  descending   induction to find automorphisms $\rho_i$ by which the $x_i$ can
  be identified with the $\xi_i$. Note that the factorization conditions are
  maintained under automorphisms. Since $x_r$ has maximal rank it is an
  injection $M_r\otimes A_{r, r-1}\to M_{r-1}\oplus M_r\otimes
  A_{r,r-1}=P_{r-1}$. Hence we can find an automorphism $\rho_{r-1}$ of
  $P_{r-1}$ such that $\rho_{r-1}\circ x_r$ becomes $\xi_r$. Let us assume now
  that modulo some automorphisms $\rho_{r-1},\ldots, \rho_i$ we have
  $x_j=\xi_j$ for $j>i$. We are going to find an automorphism $\rho_{i-1}$
  such that $\rho_{i-1}\circ x_i=\xi_i$. Because of the rank condition we can
  assume that $\oplus_{i\leq j} M_j\otimes A_{j,i-1}$ is the image of $x_i$ in
  $P_{i-1}$. Now using all the $x_i\circ\widetilde{\xi}_{i+1}\circ \ldots\circ
  \widetilde{\xi_k}$ we find that $x_i$ has a factorization through the
  standard map 
\[
P_i\otimes A_{i,i-1}\to \underset{i\leq j}\oplus M_j\otimes
A_{j,i-1}\xrightarrow{\bar{x}_i} M_{i-1}\oplus \underset{i\leq j}\oplus
  M_j\otimes A_{j,i-1}.
\]
induced by the pairings.
Now the rank condition implies that $\bar{x_i}$ induces an automorphism on
$\oplus_{i\leq j} M_j\otimes A_{j,i-1}$. This can be used to make $\bar{x_i}$
the identity via an automorphism $\rho'_{i-1}$. Now $x_i=\xi_i$. By the
analogous dual procedure we can also find automorphism $\sigma_l\in \Aut(Q_l)$
such that we can assume that $y_l=\eta_l$. Finally the factorizations
$(L_{li})$ or $(R_{li})$ resulting from (3) and (3$^\ast$) yield mappings
$\Phi_{li}$ or $\Psi_{li}$ from which we get $\phi_{li}$ as composition
\[
M_i\otimes H^\ast_{li}\rightarrowtail P_i\otimes H_{li}\xrightarrow{\Phi_{li}}
Q_l\twoheadrightarrow N_l.
\]
It follows from the special type of the $\xi_i$ and $\eta_l$ that these are
original components of an element $w=(\phi_{li})$ inducing $\gamma(w)=\gamma$.
\end{proof}
\end{subsub}
\sepprop 

\begin{subsub} \label{descrbdZ}{\bf Corollary}: With the same notation as in
\ref{descrZ}, if ${\bf w}\in\bar{Z}\smallsetminus Z$, then
  
(1) $rank\ x_i\leq rank\ \xi_i$ and $rank\ y_l\leq rank\ \eta_l$ with strict
inequality for at least one $i$ or $l$, and 

(2), (2$^\ast$), (3), (3$^\ast$) of \ref{descrZ} are satisfied.
 
\begin{proof} All conditions are closed and thus hold for points in
  $\bar{Z}$. If ${\bf w}\in\bar{Z}\smallsetminus Z$ then by \ref{descrZ}
  equality in (1) cannot hold for all $i$ and $l$.
\end{proof}\end{subsub}

\sepprop

We are going to derive effective sufficient conditions for the projectivity of
the quotient in the cases $(2,1), (2,2), (3,1)$.
\sepprop

\begin{subsub}\label{proj1} {\bf Proposition}: Let the polarizations $\Lambda$
  and $\widetilde{\Lambda}$ be as in proposition \ref{equality} and let
  $Z={\bf G}\zeta(W)$. Then $\bar{Z}\smallsetminus Z$ contains no semi--stable
point in the following cases
\begin{enumerate}
\item [(i)] $(r,s)=(2,1)\and \lambda_2\geq c_1(m_2)
  a_{21}\mu_1$

\item [(ii)] $(r,s)=(2,2)\quad\text{and}$
$$\lambda_2\geq  (\mu_1c_1(m_2)+\mu_2(c_2(m_2)-b_{21}c_1(m_2))a_{21},
\qquad
\mu_1\geq (\lambda_1(d_1(n_1))-d_2(n_1)a_{21})+\lambda_2 d_2(n_1))b_{21}.$$
\end{enumerate}

\begin{proof} We present only the case (ii), case(i) is an easier version of
  (ii). Let $(x_2, \gamma, y_1)\in\bar{Z}\smallsetminus Z$ and let us assume
  that $rank\ x_2$ is not maximal. Let $K$ be the kernel of $M_2\otimes
  A_{21}\xrightarrow{x_2} P_1$ and let $M'_2\subset M_2$ be the smallest
  subspace such that $K$ is contained in $M'_2\otimes A_{21}$. We  put
  $P'_2=M'_2, $
\[
P'_1=x_2(M'_2\otimes A_{21})\ ,\ Q'_2=\gamma(P'_1\otimes
  H^\ast_{21})\ \text{ and }\ Q'_1=y_1(Q'_2\otimes B_{21})
\]
and consider
\[
\Delta=\alpha_1 p'_1+\alpha_2p'_2-\beta_1 q'_1-\beta_2q'_2.
\]
By definition $p'_1 = \dim(M'_2\otimes A_{21}/K)$. Diagram $(L_2)$ reduces in 
our case, with $M_2$ replaced by $M'_2$, to
\[
\xymatrix{M'_2\otimes A_{21}\otimes
  H^\ast_{21}\ar@{->>}[d]_{\delta_2}\ar@{->>}[r]^-{\scriptstyle{x_2\otimes id}} & P'_1\otimes
  H^\ast_{21}\ar@{->>}[r]^{\gamma} & Q'_2\\
M'_2\otimes H^\ast_{22}\ar@{->>}[urr]_{\gamma_{22}}
}
\]

and $\gamma_{22}$ vanishes on $\delta_2(K\otimes H^\ast_{21})$ because $K$ is
the kernel of $x_2$. Therefore
\[
q'_2\leq \dim(M'_2\otimes H^\ast_{22}/\delta_2(K\otimes H^\ast_{21}))
\leq c_2(m'_2)p'_1.
\]
In order to estimate $q'_1$ we consider diagram $(L_{21})$ enlarged by the
commutative square of induced pairings
\[
\xymatrix{
M'_2\otimes A_{21}\otimes H^\ast_{21}\otimes B_{21}\ar@{->>}[d]\ar@{->>}[r] & M'_2\otimes
H^\ast_{22}\otimes B_{21}
\ar@{->>}[d]\ar@{->>}[r]^-{\scriptstyle{\gamma_{22}\otimes id}} & Q'_2\otimes
B_{21}\ar@{->>}[r]^-{y_1} & Q'_1\\
M'_2\otimes A_{21}\otimes H^\ast_{11} \ar@{->>}[r]^-{\delta_1} & M'_2\otimes
H^\ast_{12} \ar@{->>}[urr]_{\Phi_{12}}
}.
\]
Again the map $\Phi_{12}$ vanishes on $\delta_1(K\otimes H^\ast_{11})$ and we
  get
\[
q'_1\leq \dim(M'_2\otimes H^\ast_{12}/\delta_1(K\otimes H^\ast_{11}))
\leq c_1(m'_2)p'_1.
\]
Now we have the estimate
\[
\Delta\geq \alpha_2p'_2+(\alpha_1-\beta_1 c_1(m_2)-\beta_2c_2(m_2))p'_1.
\]

Therefore the condition $\alpha_1\geq \beta_1c_1(m_2)+\beta_2c_2(m_2)$ would
be sufficient, because $\alpha_2p'_2>0$. We modify the last estimate as
follows. Since the weights in case $(2,2)$ are related by 
\[
\begin{array}{lcl}
\lambda_1 & = & \alpha_1\\
\lambda_2 & = & a_{21}\alpha_1+\alpha_2
\end{array}\quad \text{ and }\quad
\begin{array}{lcl}
\mu_2 & = & \beta_2\\
\mu_1 & = & \beta_1+\beta_2b_{21}
\end{array}
\]
and since we have
\[
\lambda_2-a_{21}\lambda_1>0\quad\text{ and }\quad p'_2 a_{21}-p'_1>0,
\]
we get the estimate
\[
\Delta>(\frac{\lambda_2}{a_{21}}-\mu_1c_1(m_2)-\mu_2 c_2(m_2)+\mu_2c_1(m_2)
b_{21})p'_1.
\]
This shows that $\Delta>0$ if $x_2$ is degenerate and the first condition of 
(ii) is satisfied. In case $rank\ y_1$ is not maximal the second condition 
follows by the dual procedure.
\end{proof}
\end{subsub}
\sepsub

\begin{sub}\label{proj2} {\sc\small The case (3,1)}\rm 

In order to derive a similar result in case $(3,1)$ we introduce the
additional constant $c'_3(m_3)$ analogous to $c_3(m_3):=c_1(0, m_3)$ in
\ref{constants}. Let
\[
M_3\otimes A_{32}\otimes H^\ast_{12}\xrightarrow{\tau} M_3\otimes H^\ast_{13}
\]
be the linear map induced by the pairing and let $\kk$ be the family of all
proper
subspaces $K\subset M_3\otimes A_{32}$ which are not contained in $M'_3\otimes
A_{32}$ for any subspace $M'_3\subset M_3$ different from $M_3$. We put
\[
c'_3(m_3)=\underset{K\in\kk}{sup}\ \dfrac{\codim(\tau(K\otimes
  H^\ast_{12}))}{\codim(K)}
\]
For brevity we write $c'_3=c'_3(m_3),\ c_3=c_3(m_3)=c_1(0,m_3)$ and
$c_1=c_1(m_2, m_3).$
\sepprop

\begin{subsub}\label{projestim}{\bf Proposition}:  Let $(r,s)=(3,1)$, let
  $\Lambda=(\lambda_1, \lambda_2, \lambda_3, -\mu_1)$ be a polarization for
$(W,G)$ and $\widetilde\Lambda=(\alpha_1, \alpha_2, \alpha_3, -\mu_1)$ be
the associated polarization for $({\bf W},{\bf G})$, and assume that all 
$\alpha_i>0.$ (In this case $\mu_1=\frac{1}{n_1}$.) If 
\begin{enumerate}
    \item [(1)] $\alpha_2 c_3 + \lambda_1 c'_3\geq \mu_1 c_3 c'_3$
    \item [(2)] $\lambda_2\geq a_{21}\mu_1c_1$
    \item [(3)] $\lambda_3\geq a_{31}\mu_1 c_1$
\end{enumerate}
then $\bar{Z}\smallsetminus Z$ contains no semi--stable point.

Moreover, condition (1) may be replaced by any of the conditions
\begin{enumerate}
\item [(i)] $\lambda_3\geq \mu_1 c'_3 a_{32}+a_{31}\lambda_1$
\item [(ii)] $\lambda_3\geq \mu_1 c_3a_{31}+a_{32}\alpha_2$ 
\item [(iii)] $\lambda_3\geq \mu_1 c_3 a_{32}a_{21}$
\end{enumerate}
\end{subsub}

\sepprop

\rm{\bf Remark}: $\bar{Z}\smallsetminus Z$ contains no semi--stable point also 
in each of the following cases
\begin{enumerate}
\item [(a)] $\lambda_1\geq \mu_1 c_3$
\item [(b)] $\alpha_2\geq \mu_1 c'_3$
\item [(c)] $\alpha_3\geq \mu_1c_3a_{31}\quad \text{or}
\quad \alpha_3\geq \mu_1c'_3 a_{32}$.
\end{enumerate}
This can be seen by a direct estimate of the discriminant $\Delta$ after
substituting for $q'_1$ in the following proof.

\begin{proof} Let $(x_2, x_3, \gamma)\in \bar{Z}\smallsetminus Z$. We
distinguish the following cases of degeneracy of $x_2$ and $x_3$.

{\em case 1: $x_3$ is injective}: Then by the proof of \ref{descrZ} we can
assume that $x_3=\xi_3$ is the canonical embedding and that $x_{13}$ and $x_2$
have a factorization $\bar{x}_2$ in the following diagram
\[
\xymatrix{
M_3\otimes A_{32}\otimes A_{21}\ar@{->>}[d]\ar@{>->}[r]^-{\widetilde{\xi}_3} & (M_2\oplus
M_3\otimes A_{32})\otimes A_{21} \ar[d]\ar[r]^-{x_2} & P_1\\
M_3\otimes A_{31}\ar[urr]^{x_{13}}\ar@{>->}[r]_-{\xi'_3} & M_2\otimes
  A_{21}\oplus M_3\otimes A_{31}\ar[ur]_{\bar{x}_2}
}.
\]
Here also $\xi'_3$ is the canonical embedding. Moreover it is easy to
verify that in this case also the composed map $\gamma\circ(\bar{x}_2\otimes
id)$ admits a decomposition
\[
(M_2\otimes A_{21}\otimes H^\ast_{11})\oplus (M_3\otimes A_{31}\otimes
H^\ast_{11})\xrightarrow{\delta_1} (M_2\otimes H^\ast_{12})\oplus (M_3\otimes
H^\ast_{13})\xrightarrow{\bar{\gamma}} Q_1.
\]

Here $K=Ker(\bar{x}_2)\neq 0$ since $\bar{x}_2$ cannot be injective by the
assumption on its rank. 
We choose subspaces $M'_2, M'_3$ such that
\[
K\subset M'_2\otimes A_{21}\oplus M'_3\otimes A_{31}
\]
and such that these subspaces are minimal with this property. Now we consider
the spaces 
\[
P'_3=M'_3, \quad P'_2=M'_2\oplus (M'_3\otimes A_{32}), \quad
P'_1=x_2(P'_2\otimes A_{21}),\quad Q'_1=\gamma (P'_1\otimes H^\ast_{11})
\]
and their discriminant
\[
\Delta=\alpha_1p'_1+\alpha_2p'_2+\alpha_3p'_3-\beta_1q'_1.
\]

By the definition of the constant $c_1(m'_2, m'_3)$ and the diagram
\[
\xymatrix{
(M'_2\otimes A_{21} \oplus M'_3\otimes A_{31})\otimes
H^\ast_{11}\ar@{>>}[d]\ar@{>>}[r] & P'_1\otimes H^\ast_{11}\ar@{>>}[r] &
Q'_1\\
M'_2\otimes H^\ast_{12} \oplus M'_3\otimes H^\ast_{13}\ar@{>>}[urr]
}
\]
we obtain the estimate
\[
q'_1\leq c_1(m'_2, m'_3)p'_1\leq c_1(m_2, m_3) p'_1,
\]
where by the definition of $P'_1$ we have $p'_1=m'_2 a_{2}+m'_3
a_{31}-k.$ Inserting this we obtain
\[
\Delta\geq (\mu_1 c_1-\lambda_1)k +(\lambda_2-\mu_1 c_1 a_{21})
m'_2+(\lambda_3-\mu_1 c_1 a_{31})m'_3
\]
If $\mu_1 c_1-\lambda_1>0$, conditions (2) and (3) imply that $\Delta>0$. If,
however, $\lambda_1\geq \mu_1 c_1$ we have the direct estimate
\[
\Delta\geq (\lambda_1-\mu_1 c_1) p'_1+\alpha_2 p'_2+\alpha_3 p'_3 >0.
\]
This proves the proposition in the first case.

{\em case 2: $x_3$ is not injective}

Here we let $K$ denote the kernel of $x_3$ and we choose a subspace $M'_3
\subset M_3$ such that $K\subset M'_3\otimes A_{32}$ and $M'_3$ is minimal 
with this property. Then we consider the subspaces
\[
P'_3=M'_3,\quad P'_2=x_3(M'_3\otimes A_{32}),\quad P'_1=x_2(P'_2\otimes A_{21}),
\quad Q'_1=\gamma(P'_1\otimes H^\ast_{11}).
\]
We have the exact sequences
\[
\begin{array}{lclclclcl}
0 & \to & K & \to & M'_3\otimes A_{32} & \xrightarrow{x_2} & P'_2 & \to 0\\
0 & \to & L &\to & M'_3\otimes A_{31} & \xrightarrow{x_{13}} & P'_1 & \to & 0
\end{array}
\]
where $L$ denotes the kernel of $x_{13}$. From the factorization properties
restricted to the spaces $P'_i$ and $Q'_1$ we extract the following
commutative diagram of surjections
\[
\xymatrix{
& M'_3\otimes A_{31}\otimes H^\ast_{11}\ar[d]^{\delta_1}\ar[r] & P'_1\otimes
H^\ast_{11}\ar[d]^{\gamma}\\
M'_3\otimes A_{32}\otimes A_{21}\otimes H^\ast_{11}\ar[ur] \ar[dr] & M'_3\otimes
H^\ast_{13}\ar[r]^{\gamma_{13}} & Q'_1\\
& M'_3 \otimes A_{32}\otimes H^\ast_{12}\ar[r] \ar[u]_{\tau} & P'_2\otimes
H^\ast_{12}\ar[u]_{\gamma_{12}}
}\raisebox{-3cm}{.}
\]
From this we get again the estimates
\[
q'_1\leq c_3(m'_3)p'_1\leq c_3(m_3)p'_1\quad \text{and}\quad 
 q'_1\leq c'_3(m'_3)p'_2\leq c'_3(m_3)p'_2,
\]
where $p'_1=m'_3a_{31}-l$ and $p'_2=m'_3 a_{32}-k$.
Let $0<t<1$ be a real number. Then we have
\[
q'_1\leq t c'_3p'_2+(1-t) c_3p'_1.
\]
Substituting this into the discriminant we get
\[
\Delta\geq(\lambda_1-(1-t)\mu_1 c_3)p'_1+(\alpha_2-t\mu_1
c'_3)p'_2+\alpha_3m'_3.
\]
Now condition (1) enables us to find $t$ with
\[
1-\frac{\lambda_1}{\mu_1c_3}\leq t\leq \frac{\alpha_2}{\mu_1c'_3}\ ,
\]
such that the first two terms of the estimate are non--negative. Therefore
$\Delta>0$, and again $(x_2, x_3, \gamma)$ is not semi--stable.

In order to show that (1) can be replaced by one of (i), (ii) or (iii) we
substitute $\alpha_i$ and $p'_i$ and get after cancelation
\[
\begin{array}{lcl}
\Delta & = & -\lambda_1 l -\alpha_2k+\lambda_3 m'_3-\mu_1q'_1\\
       & \geq & -\lambda_1 l-\alpha_2k+\lambda_3 m'_3-\mu_1c'_3(m'_3a_{32}-k)\\
       & = &-\lambda_1 l +(\mu_1 c'_3-\alpha_2)k+( \lambda_3-\mu_1c'_3 a_{32})m'_3
\end{array}\hskip3mm\raisebox{-6mm}{.}
\]
If $\alpha_2\geq \mu_1c'_3$, then by a direct estimate we get
$\Delta>0$. Therefore we may assume that $\mu_1 c'_3-\alpha_2>0$. Since in
addition $l\leq m'_3a_{31}$, we get
\[
\Delta>(\lambda_3-\mu_1 c'_3 a_{32}- a_{31}\lambda_1)m'_3.
\]
This shows that (1) can be replaced by (i). In the same way one shows that (1)
can be replaced by (ii), using the other estimate of $q'_1$. That finally (1)
can be replaced by (iii) can be shown by substituting first $m'_3\geq
\frac{p'_2}{a_{32}}$ and canceling $\alpha_2p'_2$ and then substituting
$p'_2\geq \frac{p'_1}{a_{21}}$ to get
\[
a_{32} a_{21}\Delta \geq \lambda_1p'_1(a_{32}
a_{21}-a_{31})+(\lambda_3-\mu_1 c_3a_{32} a_{21}) p'_1.
\]
\end{proof}
\end{sub}

\sepsub

\begin{sub}{\sc\small Proof of theorems \ref{theo_main} and 
\ref{theo_main2}}\rm
\end{sub}

Theorem \ref{theo_main} is an immediate consequence of proposition 
\ref{goodqu}, corollary \ref{sone} and proposition \ref{proj1}.
Theorem \ref{theo_main2} follows immediately from theorem \ref{theo_main}
and \ref{constant_calc}. 

\sepsec
%\newpage
\section{Examples}

\begin{sub}\label{constant_calc} {\sc\small Constants}\rm

We give here some constants (cf. \ref{constants}) used in the examples.
The following result is proved in \cite{other}, prop. 6.1.

\begin{subsub}\label{formula1}{\bf Proposition} : For homomorphisms of type
\[(M_1\otimes\ko(-2))\oplus(M_2\otimes\ko(-1))\to N_1\otimes\ko\]
on a projective space of dimension $n$ we have
\[
c_1(m)=\frac{m(m-1)}{2(m(n+1)-1)}\txtem{if}m\leq n+1,\quad\quad
c_1(m)=\frac{n+1}{2(n+2)}\txtem{if}m\geq n+1,
\]
\end{subsub}

\sepprop

\begin{subsub}\label{symbd2} {\bf Lemma}: For homomorphisms of type 
\[
M_1\otimes \ko(-d)\oplus\ko(-2)\oplus\ko(-1)\to N_1\otimes \ko
\]
on the projective space $\P V$ the constant $c_1(1,1)$ is 
$\dim(V)/\dim(S^{d-1}V)$.

\begin{proof} We put $s(p)=\dim(S^pV)$. The homomorphisms $\delta_1$ of
 \ref{constants} reduces here to the canonical map
\[
(S^{d-2} V^\ast\oplus S^{d-1}V^\ast)\otimes S^dV\to S^2V\oplus V.
\]
If $K$ is a proper subspace of $S^{d-2} V^\ast\oplus S^{d-1} V^\ast$ which is
not contained in one of the summands, it contains elements $(f,g)$ with $f\neq
0$ or elements $(f,g)$ with  $g\neq 0$. But since $f\otimes S^dV\to S^2V$ is
surjective, the map $\delta(K)\to S^2V$ is surjective. Hence 
$\codim(\delta(K))\leq s(1)$. 
If $K$ contains an element $(0,g)$ with $g\neq 0$, then
$\delta(K)=S^2V\oplus V$. For then $\delta(K)$ contains $V$, and since
$\delta(K)\to S^2V$ is surjective, if follows that $\delta (K)=S^2 V\oplus
V$. Therefore, if $\codim(\delta(K))>0$, there is a basis $(f_1, g_1), \ldots,
(f_k, g_k)$ of $K$ with $f_1, \ldots f_k$ linearly independent, i.e. 
$\dim(K)\leq s(d-2)$ or $\codim(K)\geq s(d-1)$. Therefore $c_1(1,1)\leq
s(1)/s(d-1)$. But now we can find subspaces which realize this bound. For any
$z\in V^\ast$ we let $K$ be the space of all $(f, fz), f\in
S^{d-2}V^\ast$. Then $K\cong S^{d-2}V^\ast$ and it follows also that in this
case $\delta(K) \cong S^2V$. Then $\codim(\delta(K))/\codim(K)=s(1)/s(d-1)$. 
\end{proof}
\end{subsub}
\end{sub}

\sepsub

\begin{sub}\label{ex1} {\sc First example of type $(2,1)$}

\rm We use the abbreviation $m\kf$ for $\C^m\otimes \kf$ for a sheaf
and a positive integer and consider here homomorphisms
\[
2\ko(-2)\oplus \ko(-1)\xrightarrow{(\phi_1,\phi_2)} 3\ko
\]
over $\P_2$ of type $(2,1)$. The polarization $\Lambda=(\lambda_1,
\lambda_2, -\mu_1)$ is supposed to be proper for $W$ and ${\mathbf W}$,
i.e. $\lambda_i>0$ and $\alpha_i>0$ for all $i$. The only constant involved
here is $c_1(m_2)=c(1)=0$. Therefore the conditions of \ref{sone} and
\ref{proj1} are automatically satisfied by
$\alpha_2=\lambda_2-3\lambda_1>0$. Hence all the quotients of $\wss$ will be
good and projective under this condition. Since $2\lambda_1+\lambda_2=1$ and
$3\mu_1=1$, we can replace the polarization by the rational number
$t=\lambda_2 >\frac{3}{5}$ (cf. \ref{ex1gen}). The numerical condition for
(semi--)stability then becomes
\[
\Delta=\frac{1-t}{2}m_1+tm_2-\frac{1}{3}n<0\ (\leq 0),
\]
where $(m_1, m_2, n)$ is the dimension vector of a $(\phi_1,
\phi_2)$--invariant sub-family of vector spaces, such that $m_1\leq 2,\ 
m_2\leq 1,\ n\leq 3$. One can easily check that $t=\frac{2}{3}$ is the only value for
which $\Delta$ might be zero, and this is the case for the values $(0,1,2)$
and $(2,0,1)$. And indeed, the homomorphisms $\phi$ given by matrices
\[
\left(
  \begin{array}{ccc}
\ast & \ast & 0\\
\ast & \ast & z_2\\
\ast & \ast & z_3
  \end{array}
\right)\qquad \text{ and }\qquad \left(
  \begin{array}{ccc}
0 & 0 & z_1\\
0 & 0 & z_2\\
\ast & \ast & z_3
  \end{array}
\right)
\]
with generically chosen entries and linear forms $z_i$ are semi--stable and
not stable for $t=\frac{2}{3}$.
\end{sub}
\sepsubsub

\begin{subsub}\label{ex1plus} The case $t>\frac{2}{3}$

\rm It is easy to show that in this case $(\phi_1, \phi_2)$ is $t$--stable
if and only if
\begin{itemize}
\item $\phi_2$ is nowhere zero
\item for any $(\psi_1, \psi_2)=h.(\phi_1, \phi_2)$ with $h\in H$ and any
  $1$--dimensional subspace $M_1\subset {\mathbf C}^2 \text{ we have }\psi_1(M_1(\otimes
  \ko(-2))\neq 0$.
\end{itemize}

The first condition says that $\coker(\phi_2)$ is isomorphic to the universal
quotient bundle $Q$ on $\P_2$, and the second that $\phi_1$ induces a
$2$--dimensional subspace of $H^0 Q(2)$. It follows that the sets $W^s(t)$
of stable points are the same for $t>\frac{2}{3}$, which we denote by
$W^s_+$. Moreover, from the above characterization of stable homomorphism we
deduce that the geometric quotient $M_+=W^s_+/G$ is isomorphic to the
Grassmannian
\[
M_+\cong Gr(2, H^0 Q(2))
\]
which is smooth of dimension $26$. There is an interesting subvariety
$Z\subset M_+$ which consists of the images of the homomorphisms 
$$
\left(
  \begin{array}{ccc}
0 & 0 & z_1\\
0 & 0 & z_2\\
\ast & \ast & z_3
  \end{array}
\right)\eqno(1)
$$
which belong to $W^s_+$. These are those $(\phi_1, \phi_2)$ for which the
induced homomorphism $2\ko(-2)\to Q$ is not injective. We will see next that
$Z$ is isomorphic to the non--stable locus of $M_0$ below and is smooth of
dimension $10$.
\end{subsub}
\sepsubsub

\begin{subsub}\label{ex1zero} The case $t=\frac{2}{3}$

\rm We write $W^{ss}_0$ for $W^{ss}(\frac{2}{3})$. When considering the matrix
representations we find that $W^s_+\subset W^{ss}_0$ and that the remaining
part $W^{ss}_0\smallsetminus W^s_+$ consists of those homomorphisms for which
$\phi_2$ is zero in exactly one point. Such homomorphisms are equivalent to
matrices
$$\left(
  \begin{array}{ccc}
\ast & \ast & z\\
\ast & \ast & w\\
f & g & 0
  \end{array}
\right)\eqno(2)
$$
where $z,w $ are independent linear and $f,g$ are independent quadratic
forms. Note, however, that $W^s_+$ intersects the non--stable locus of
$W^{ss}_0$ in matrices equivalent to those of type (1). But the orbit closures
in $W^{ss}_0$ of both types (1) and (2) of matrices contain the direct sums
$$
\left(
  \begin{array}{ccc}
0 & 0 & z\\
0 & 0 & w\\
f & g & 0
  \end{array}
\right)\eqno(3)
$$
of independent linear and quadratic forms. From that it follows that the
induced morphism
\[
M_+\to M_0
\]
of the quotients is bijective and moreover an isomorphism by Zariski's main
theorem, because both spaces are normal. The points of the non--stable locus
$M_0\smallsetminus M_0^s$ are represented by matrices of type (3). It is again
routine to deduce from this observation that
\[
M_0\smallsetminus M_0^s\cong\P_2\times Gr(2, H^0\ko(2)).
\]
The subvariety $Z\subset M_+$ corresponds to this set under the
isomorphism. We can also identify the set $M^s_0$ of stable points with $Gr(2,
H^0Q(2))\smallsetminus Z$.
\end{subsub}

\sepsubsub

\begin{subsub}\label{ex1minus} The case $\frac{3}{5}<t<\frac{2}{3}$

\rm Similarly to the case $W^s_+$ we find that here $W^s_{-}=W^s(t)$ is
independent of $t$ and that $W^s_{-}\subset W^{ss}_0$. The remaining part
consists now of all homomorphisms which are equivalent to a matrix of type
(1). Note that now homomorphisms of type (2) are contained in $W^s_{-}$. The
induced morphism
\[
M_{-}\to M_0
\]
is again surjective but not injective over $M_0\smallsetminus M_0^s$. Let $Y$
be the inverse image of $M_0\smallsetminus M_0^s$. Then $Y$ consists of the
points which are represented by matrices of type (2) which are not equivalent
to matrices of type (3). It is easy to check that the restricted morphism
\[
M_{-}\smallsetminus Y\underset{\approx}{\longrightarrow} M_0^s
\]
is bijective and therefore also an isomorphism by Zariski's main theorem. We
are going to verify that $Y$ is a divisor in $M_{-}$. There is a morphism
\[
Y\xrightarrow{p}{\mathbf P}_2
\]
which assigns to the class of $(\phi_1, \phi_2)$ the point $x$ at which
$\phi_2$ is degenerate. In this case
\[
\coker(\phi_2)\cong \ko\oplus \ki_x (1)
\]
where $\ki_x$ is the ideal sheaf of $x$. For such $(\phi_1, \phi_2)$ we are
given an exact diagram
%\[
%\begin{array}{ccccccccc}
%& & 0 & & & & & &\\
%& &\uparrow & & & & & & \\
%& & 2\ko(-2)& & & & \ko & & \\
%& &\uparrow & & & & \uparrow & & \\
%0 & \to & 2\ko(-2)\oplus \ko(-1) & \xrightarrow{(\phi_1, \phi_2)} & 3\ko & \to &
%\kf & \to & 0\\
%&  & \uparrow & & \| & &\uparrow & &\\
%0 & \to & \ko(-1) &\xrightarrow{\phi_2} & 3\ko & \to & \ko\oplus\ki_x(1) & \to
%& 0 \\
%& & \uparrow & & & & \uparrow & & \\
%& & 0 & & & & 2\ko(-2)& &\\
%& &   & & & & \uparrow & &\\
%& &   & & & & 0& &\\
%\end{array}
%\]
\[\xymatrix{
& 0 \\
& 2\ko(-2)\ar[u] & & & \ko\\
0\ar[r] & 2\ko(-2)\oplus\ko(-1)\ar[u]\ar[rr]^-{(\Phi_1,\Phi_2)} & & 
3\ko\ar[r] & \kf\ar[r]\ar[u] & 0 \\
0\ar[r] & \ko(-1)\ar[u]\ar[rr]^-{\Phi_2} & & 3\ko\fleq[u]\ar[r] & 
\ko\oplus\ki_x(1)\ar[u]\ar[r] & 0\\
& 0\ar[u] & & & 2\ko(-2)\ar[u]\\
& & & & 0\ar[u]
}\]
such that $(\phi_1, \phi_2)$ corresponds to a $2$--dimensional subspace
$\Gamma\subset H^0(\ko(2)\oplus\ki_x(3))$. The condition of defining a element
of $Y$ is that $\Gamma$ is neither contained in $H^0\ki_x(3)$ nor in
$H^0(\ko(2)) s$ for any section $s$ of $\ko\oplus\ki_x(1)$. We let $U_x\subset
Gr(2, H^0(\ko(2)\oplus\ki_x(3))$ denote the open subvariety of such
$\Gamma$. By assigning to $\Gamma$ the class of $(\phi_2, \phi_2)$ where
$\phi_1$ is defined by a lifting in the above diagram, we get a morphism
$U_x\to M_{-}$ whose image is the fibre $Y_x=p^{-1}(x)$. The morphism
\[
U_x\twoheadrightarrow Y_x
\]
is nothing but the quotient of $U_x$ by the algebraic group
$\Aut(\ko\oplus\ki_x(1))$. It follows that $Y_x$ is a variety of dimension
23. Using the techniques of this paper for this quotient, we can even prove
that $Y$ is smooth. Finally $Y$ has dimension 25 and thus is a divisor in
the irreducible and normal variety $M_{-}$.
\end{subsub}
\vskip5mm

{\em Remarks}: (1) One would like to interpret the matrices of type (2) as
representing extensions of the sheaves $\coker(f,g)$ and
$\ki_x(1)=\coker\binom{z}{w}$. Indeed a matrix of type (2) defines such an
extension, but this extension is isomorphic to the direct sum.

(2) The above correspondence between $(\phi_1, \phi_2)$ and $\Gamma$ indicates
that the quotient spaces considered here are spaces of coherent systems as 
in \cite{potier1}. 

\sepsubsub
%\newpage
\begin{subsub}\label{flip1} The flip

\rm The diagram $M_{-}\to M_0\overset{\approx}{\leftarrow} M_+$ can be
interpreted as a flip. It is induced by the inclusions $W^s_{-}\subset
W^{ss}_0\supset W^s_+$. The orbits of stable points of type (2) in $W^s_{-}$ and
of type (1) in $W^s_+$ don't intersect in $W^{ss}_0$ but so do their closures
in $W^{ss}_0$. Thus the fibres of $M_{-}\to M_0$ and $M_0\leftarrow{M_+}$
correspond to the two different types of semi--stable orbits in $W^{ss}_0$
defining the same points in $M_0\smallsetminus M_0^s$.
\end{subsub}

\sepsub

\begin{sub}\label {ex1gen} {\sc\small General homomorphisms of type $(2,1)$}

\rm In a more general situation of type $(2,1)$ we consider homomorphisms
\[
m_1\ko(-2)\oplus m_2\ko(-1)\to n_1\ko
\]
over $\P_n$.  A polarization in this
case is determined by the rational number $t=m_2\lambda_2$ with $0<t<1$ and
$1-t=m_1\lambda_1,\ \mu_1=1/n_1$. A $\Lambda$--(semi--)stable homomorphism is
then called $t$--(semi--)stable. We write $W^{ss}(t)$ and $W^s(t)$ for $\wss$
and $\ws$. In terms of $t$ the conditions are
\[
1>t>\frac{(n+1)m_2}{(n+1)m_2+m_1}\quad \text{ and }\quad
t\geq\frac{(n+1)m_2}{n_1} c_1(m_2).
\]
The constant $c_1(m_2)$ is given in proposition \ref{formula1}.
Such polarizations exist if and only if 
\[
n_1>(n+1) m_2 c_1(m_2).
\]

In order to measure $t$--stability we introduce the numbers
\[
r_1=\frac{m'_1}{m_1},\quad r_2=\frac{m'_1}{m_2},\quad s_2=\frac{n'_1}{n_1}
\]
and call $(r_1, r_2, s_1)$ $\phi$--admissible if there are subspaces
$M'_1\subset M_1, M'_2\subset M_2, N'_1\subset N_1$ of dimensions $m'_1, m'_2,
n'_1$ such that $\phi$ maps $M'_1\otimes \ko(-2)\oplus M'_2\otimes \ko(-1)$
into $N'_1\otimes \ko$. Then $\phi$ is $t$--(semi--)stable if and only if for
any $\phi$--admissible proper triple $(r_1, r_2, s_1)$, i.e. a triple which is
neither $(0,0,0)$ or $(1,1,1)$, we have
\[
\Delta_t=(1-t)r_1+tr_2-s_1<0\ (\leq 0).
\]
A polarization $t$ is called critical if there are proper triples with
$\Delta_t=0$. Thus the critical values of $t$ are the rational numbers
\[
\frac{s_1-r_1}{r_2-r_1}\ ,
\]
where we may assume $s_1\neq 0, 1$ and thus $r_1\neq r_2$. We let $t_{max}$ be
the maximal critical value if there are such with $0<t<1$ and put $t_{max}=0$
otherwise. If $t$ is not critical we have $W^s(t)=W^{ss}(t)$. 
\end{sub}
\sepprop

\begin{subsub}\label{ex1stab} {\bf Lemma}: Suppose that $m_2$ and $n_1$ are
  relatively prime and that $t_{max}<t<1$. Then $\phi=(\phi_1, \phi_2)$ is
  $t$--stable if and only if

\begin{enumerate}
\item [(1)] $\phi_2$ is stable with respect to the group $\GL(M_2)\times
  \GL(N_1)$.
\item [(2)] For any $1$--dimensional subspace $\C\overset{j}{\hookrightarrow} 
M_1$, and any $h\in \Hom(M_1\otimes \ko(-2),
  M_2\otimes \ko(-1))$ the map $(\varphi_1+h\circ\varphi_2)\circ j\ :\ 
\ko(-2)\to N_1\otimes \ko$ is not zero.
\end{enumerate}

\begin{proof}By the characterization of stability in section 3 the
  homomorphism $\phi_2$ is stable if and only if for any proper pair
  $M'_2\subset M_2, N'_1\subset N_1$ of $\phi_2$--admissible subspaces
  $r_2<s_1$. Now let $(\phi_1, \phi_2)$ be stable. If $\phi_2$ were not stable
  there would be a  proper\ $\phi_2$--admissible pair $(r_2, s_1)$ with
  $s_1\leq   r_2$. But then $s_1 <r_2$ because $m_2, n_1$ are supposed to be
  relatively  prime. Then $s_1/r_2<t$ because $s_1/r_2$ is a critical value
  and thus $\Delta_t=r_2t-s_1>0$, contradicting the   stability of $(\phi_1,
  \phi_2)$. The condition (2) is trivially satisfied if   $(\phi_1, \phi_2)$
  is $t$--stable, because otherwise $(1,0,0)$ would be admissible with
  $\Delta_t=1-t>0$. We have to show now that conversely (1), (2) imply that
  $(\phi_1, \phi_2)$ is $t$--stable. For this let $(r_1, r_2, s_1)$ be a
  proper $(\phi_1, \phi_2)$--admissible triple. If $r_1\leq r_2$ and   $r_2=0$,
  there is nothing to prove. If $r_2>0$ then $r_2<s_1$ by (1) and we   have
  $t(r_2-r_1)<s_1-r_1$ and hence $\Delta_t<0$. If however $r_2<r_1$ we   have
  $\Delta_t<0$ in case $r_1\leq s_1$. Since the case $s_1\leq r_2$ is only
  possible if $s_1=r_2=0$ and then $r_1=0$ by (2), we can assume that 
  $r_2<s_1<r_1$. But then 
\[
\frac{r_1-s_1}{r_1-r_2}<t
\]
because the fraction is a critical value, and last inequality is the
inequality $\Delta_t<0$.
\end{proof}
\end{subsub}
\sepprop

Now we are able to describe the space $M_+=W^s(t)/G$ for $t_{max}<t$ which is
independent of $t$. According to the lemma $W^s(t)$ can only be non--empty if
there are stable morphisms $\phi_2$. This is the case if and only if
\[
\frac{1}{\sigma(n)} <\frac{n_1}{m_2}<\sigma (n)
\]
where $\sigma(n)=\frac{1}{2}(n+1+\sqrt{(n+1)^2-4})$, see \cite{dr1}. We
restrict ourselves now to the case where in addition to the previous
conditions on $n_1, m_2$ we have $n_1\geq n m_2$ and $(n_1, m_2)=1$. Then a
stable $\phi_2$ is injective and a subbundle (except at finite number of
points in case $n_1=nm_2$, see \cite{dr1}, \cite{dr3}). The quotient space of
this space of stable homomorphisms by $\GL(M_2)\times \GL(N_1)$ is denoted by
$N=N(n+1, m_2, n_1)$. It is a smooth projective variety and there is a
universal sheaf $\ke$ on $N\times \P_n$. For $x\in N$ let $\ke_x$ denote the 
fibre
sheaf representing $x$. Since it is the cokernel of the
representing homomorphism $\phi_2$, we get
\[
h^0\ke_x(2)=(n+1)\left(\frac{n_1(n+2)}{2}-m_2\right).
\]
Therefore $p_\ast\ke(2)$ is locally free on $N$ where $p$ denotes the first
projection of $N\times \P_n$. Now $M_+$ can be non--empty only if 
\[
m_1\leq (n+1)\left(\frac{n_1(n+2)}{2}-m_2\right).
\]
If conversely this is the case for any stable $\phi_2$ and any subspace
$M_1\subset H^0\ke_x(2)$ where $x=[\phi_2]$, there is a lifting $\phi_1:
M_1\otimes \ko(-2)\to N_1\otimes \ko$ of $M_1\otimes \ko(-2)\to\ke_x$, and
$(\phi_1, \phi_2)$ satisfies (1), (2) of the lemma. It follows now easily by
considering corresponding families that
\[
M_+\cong Gr_N(m_1, p_\ast\ke(2))
\]
where $Gr_N$ denotes the relative Grassmannian. It is more difficult to
characterize the other moduli spaces $M(t)=W^{ss}(t)/G$ for the intervals
between the critical values or for the critical values and to interpret the
flips between them.
\sepsub

\begin{sub}\label{ex2} {\sc Example of type $(2,2)$}

\rm We consider now a simple example of type $(2,2)$ on $\P_3$ of
homomorphisms
\[
\ko(-2)\oplus\ko(-1)\xrightarrow{\phi} \ko\oplus 3\ko(1).
\]
Again the polarizations $\Lambda=(\lambda_1, \lambda_2, -\mu_1, -\mu_2)$ are
supposed to be proper for $W$ and ${\mathbf W}$ such that we have
$\lambda_i>0,\ \mu_l>0$ and
\[
\lambda_2>4 \lambda_1\quad \text{ and }\quad \mu_1>4\mu_2.
\]
All constants $c_l(m_2)$ and $d_i(n_1)$ are again zero, because
$m_2=n_1=1$. Then by the above conditions also the conditions for proposition
\ref{equality} and proposition \ref{proj1} are satisfied, such that there
exists a good and projective quotient $\wss//G$ for any polarization
satisfying the conditions. Since we have $\lambda_1+\lambda_2=1$ and
$\mu_1+3\mu_2=1$, the polarization $\Lambda$ is determined already by
$\lambda_2$ and $\mu_1$, for which the above conditions become
$$
1>\lambda_2>\frac{4}{5}\quad \text{ and }\quad \frac{3}{7}>1-\mu_1>0.\eqno(1)
$$
Next we derive the conditions for the occurrence of true semi--stable
points. If $(m_1, m_2, n_1, n_2)$ is the dimension vector of a
$\phi$--invariant sub-family we have to consider the equation
\[
\Delta=(1-\lambda_2)m_1+\lambda_2m_2-\mu_1n_1-\frac{1}{3}(1-\mu_1)n_2=0.
\]
By inserting all possible dimension vectors we get the 6 conditions
$$
1-\mu_1=\frac{3}{k}\lambda_2,\quad
1-\mu_1=-\frac{3}{k}\lambda_2+\frac{3}{k}\eqno(2)
$$
for $k=1,2,3$. If one of these is satisfied, there might be non--stable points
in $\wss$. In the following figure 1 the lines with the equations (2) are shown
together with the rectangle (1) (lower right), for the points of which we get 
good and projective quotients.

\vskip5mm
%\newpage
\setlength{\unitlength}{0.006in}%
\begin{picture}(960,500)(80,220)
%\thicklines
\put( 310,560){\line( 1, 0){490}}
\put(800,560){\line( 0,-1){490}}
\put(320,80){\vector(1, 0){600}}
\put (300,60){$0$}
\put(860, 100){$\lambda_2$}
\put( 330,80){\line(-1, 0){ 10}}
\put( 320,80){\vector( 0, 1){600}}
\put(290,560){$1$}
\put(800,40){$1$}
\put(340,650){$1-\mu_1$}
\put( 320,560){\line( 1,-1){480}}
\put(800,80){\line( 0, 1){  5}}
\put(800,80){\line(-2, 3){320}}
%\put(800,85){\makebox(0.4444,0.6667){\SetFigFont{10}{12}{rm}.}}
\put(800,80){\line(-1, 3){160}}
\put(640,560){\line(-2,-3){320}}
\put(800,560){\line(-1,-1){480}}
\put( 320,85){\line( 0,-1){  5}}
%\put( 320,85){\makebox(0.4444,0.6667){\SetFigFont{10}{12}{rm}.}}
%\put( 320,80){\makebox(0.4444,0.6667){\SetFigFont{10}{12}{rm}.}}
\put( 320,80){\line( 1, 3){160}}
\multiput(800,285)(-10,0){49}{\line(-1,0){5}}
%\put(800,285){\line(-1, 0){490}}
\multiput(705,560)(0,-10){49}{\line(0,-1){5}}
%\put(705,560){\line( 0,-1){490}}
\put(270,285){$3/7$}
\put(705,40){$4/5$}
\end{picture}

\vskip3cm
\begin{center}{\small Figure 1}\end{center}
%\newpage
\sepsub

The homomorphism $\phi$ defined by the matrix
\[
\left(
\begin{array}{c|c}
z_2^2-z_1z_3 & z_0\\\hline
z_0^3 & z_1^2\\
z_1^3 & z_2^2\\
z_2^3 & z_3^2
\end{array}
\right)\raisebox{-4ex}{,}
\]
where the $z_i$ are homogeneous coordinates of $\P_3$, is easily verified to be
$G$--stable for each polarization $\Lambda$ in the rectangle (1). Therefore
the moduli spaces are not empty. On each of the 3 lines in the rectangle (1)
each point defines one and the same open set $\wss$ and hence one and the same
moduli space with semi--stable and non--stable points. Similarly, on each of
the 4 open triangles we have one and the same moduli space, which is a smooth
projective geometric quotient. Each of the 7 spaces has dimension $77$. The
reader may also verify that the moduli space for an open triangle admits a
morphism to the moduli space of each of its edges, thereby defining a chain of 
flips.
\end{sub}

%\sepsub
\newpage
\begin{sub}\label{ex2gen} {\sc\small More general homomorphisms of type 
$(2,2)$}

\rm More general homomorphisms for which we know the constants explicitly are
homomorphisms of type
\[
m_1\ko(-2)\oplus2\ko(-1)\to 2\ko\oplus n_2\ko(1)
\]
over $\P_3$, say. By proposition \ref{constant_calc} the constants are  here
\[
c_1(2)=d_2(2)=\frac{1}{7}\quad \text{ and }\quad c_2(2)=d_1(2)=\frac{4}{7}.
\]
Let $W$ be the space of those homomorphisms. A proper polarization\\
$\Lambda=(\lambda_1,\lambda_2,-\mu_1,-\mu_2)$ for $W$ satisfies
\[
m_1\lambda_1+2\lambda_2=1\ ,\quad 2\mu_1+n_2\mu_2=1
\]
with $\lambda_1, \lambda_2,\mu_1,\mu_2$ positive. We will also assume that
$\alpha_2>0,\ \beta_1>0$, i.e. $\lambda_2>4\lambda_1$ and
$\mu_1>4\mu_2$. These four conditions can be replaced by
$$
\frac{4}{8+m_1}<\lambda_2<\frac{1}{2}\quad \text{ and }  \quad
\frac{4}{8+n_2}<\mu_1<\frac{1}{2}\eqno(1)
$$
\sepprop

\begin{subsub}\label{exclaim} {\bf Claim}: There are polarizations $\Lambda$
  such that $\wss$ admits a good and projective quotient in the following
  cases
  \begin{enumerate}
  \item [(i)] $m_1<6$ and $n_2<8$
\item [(i')] $m_1\leq 6$ and $n_2=8$
\item [(ii)] $8\leq m_1+3\leq n_2$ and $8m_1+8<7n_2$
  \end{enumerate}

\begin{proof} The conditions of \ref{sone} for the equivalence of
  (semi--)stability become 
$$
\lambda_2\geq \frac{4}{7}(\mu_1+4\mu_2)\quad \text{ and }\quad \mu_1\geq
\frac{16}{7}(4\lambda_2-15\lambda_1)\eqno(2)
$$
and the conditions of \ref{proj1} for the projectivity of the quotient become
$$
\lambda_2\geq \frac{4}{7}\quad \text{ and }\quad \mu_1\geq
\frac{4}{7}\lambda_2\ . \eqno(3)
$$
The first condition of (3) follows already from the first of (2). After
replacing $\lambda_1$ and $\mu_2$ conditions (2) and (3) are equivalent to
$$
\begin{array}{rcl}
\frac{7}{4}n_2\lambda_2 & \geq & (n_2-8)\mu_1+4\\[1.0ex]
\frac{7}{16}m_1\mu_1 & \geq & (4m_1+30)\lambda_2-15\\[1.0ex]
\mu_1 & \geq & \frac{4}{7}\lambda_2
\end{array}
\eqno(4)
$$
Using (1) for $\lambda_2$, we find that (4) has a solution $(\lambda_2,
\mu_1)$ if the system
$$
\begin{array}{rcl}
\frac{7n_2}{8+m_1} & \geq & (n_2-8)\mu_1+4\\[1.5ex]
\mu_1 & > & \frac{16}{7(8+m_1)}
\end{array}
\eqno(5)
$$
has a solution $\mu_1$. For this we distinguish the cases $n_2<8,\ n_2=8,\ 
8<n_2$. If $n_2<8$ the first inequality of (5) has a solution
$\mu_1<\frac{1}{2}$ if $m_1<6$. If $n_2=8$, then $m_1\leq 6$, which is case
(i'). If $n_2>8$, the first inequality of (5) reduces to 
$$
\frac{7n_2-4m_1-32}{(n_2-8)(m_1+8)}\geq \mu_1>\frac{4}{n_2+8}.\eqno(6)
$$
Then (5) has a solution $\mu_1$ if and only if 
\[
\begin{array}{rcl}
7n_2-4m_1-32 & > & 0\\
(7n_2-4m_1-32)(n_2+8) & > & 4(n_2-8)(m_1+8)\\
7(7n_2-4m_1-32) & > & 16(n_2-8)
\end{array}
\]
These inequalities reduce to
\[
\begin{array}{rcl}
7n_2 & > & 4m_1+32\\
7n_2 & > & 8m_1+8\\
33n_2 & > & 28m_1+96
\end{array}
\]
They are all satisfied if we suppose (ii) of the claim.
\end{proof}\end{subsub}

\sepprop
In {\em figure 2} the lines of the critical values of the
polarizations i.e. of the pairs $(\lambda_2, \mu_1)$ are shown together with
the small region of those pairs which satisfy the sufficient conditions (4)
for the existence of a good and projective quotient, based on the values
$m_1=3$ and $n_2=5$.
\end{sub}
\sepsub

\begin{sub}\label{ex3} {\sc\small Example of type $(3,1)$}
\rm 

As an example of type $(3,1)$ we consider only the space of homomorphisms
\[
\ko(-4)\oplus\ko(-2)\oplus\ko(-1)\to 5\ko
\]
over $\P_3$. We assume again that all $\lambda_i$ and
all $\alpha_i$ are positive. Then the conditions of \ref{sone} together with
the normalization of the polarization are
\[
\begin{array}{l}
\lambda_1+\lambda_2+\lambda_3=1\\
\lambda_2>10\lambda_1\\
\lambda_3-4\lambda_2+20\lambda_1>0.
\end{array}
\quad
\begin{array}{l}
\mu_1=\frac{1}{k}\\[1.0ex]
\lambda_2\geq \frac{4}{5} c_1(1,1)\\
\phantom{xx}
\end{array}
\]
As additional condition for the projectivity of the quotient we use condition
(a) of the remark following proposition \ref{projestim}. Since in this case 
both the constants $c_3(1)$ and $c'_3(1)$ are zero, this condition is just 
$\lambda_1\geq 0$ and is already satisfied by our assumption.

For homomorphisms of the above type the condition $\lambda_3 < \frac{4}{5}$
is necessary if $\ws\neq \emptyset.$
For if $\phi=(\phi_1, \phi_2, \phi_3)$ is an element of $W$ then $\phi_3$ has 
degree $1$ and thus contains at most $4$ independent components. Then 
$m_1=m_2=0$ and $m_3=1,\ n_1= 4$ is a choice of dimensions of $\phi$--invariant 
subspaces and the discriminant becomes $\Delta=\lambda_3-\frac{1}{5}$.

By \ref{symbd2} the value of $c_1(1,1)$ is $\frac{1}{5}$.
Now it is easy to see that there exist polarizations 
$\Lambda$ which satisfy the above inequalities. That $\ws$ is then indeed 
non--empty follows from the existence of generic matrices as in \ref{ex2}.
Moreover there are again regions of polarizations for which the sets $\wss$
are the same and which are responsible for flips.
\end{sub}
\vspace{1cm}
%\newpage
\psfig{figure=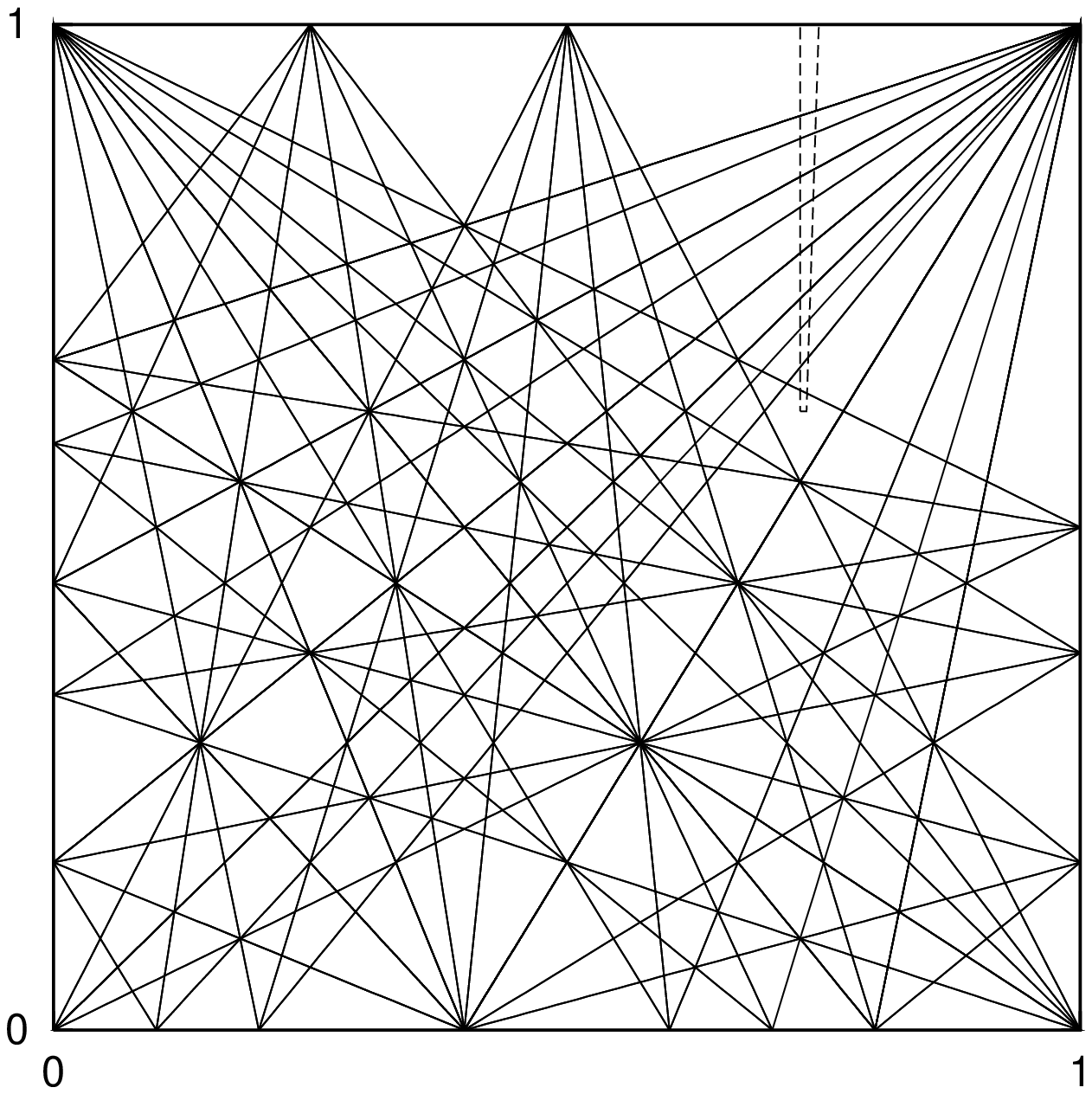}

\begin{center}{\small Figure 2\\
Here the horizontal axis represents $m_2\lambda_2$
and the vertical axis represents  $n_1\mu_1$\\ for $m_1=3$ and $n_2=5$.}
\end{center}

\sepsec
%\newpage
\section{Construction of fine moduli spaces of torsion free sheaves}
\label{Fine_Mod}

Let $n$,$k$ be integers such that $n\geq 2$ and 
\[\frac{(n+1)(n+2)}{2} \ < \ k \ \leq \ (n+1)^2 . \]
Let $V$ be a vector space of dimension $n+1$, $\P_n=\P(V)$.
We will study in this chapter morphisms of sheaves on $\P_n$ of type
\[\Phi=(\Phi_1,\Phi_2) : \ko(-2)\ot\C^2\lra
\ko(-1)\oplus(\ko\ot\C^k) .\]
Let
\[f_1 : \C^2\lra V^*\]
the linear map induced by $\Phi_1$. For semistable morphisms (with respect
to a given polarization) $f_1$ is non zero. So it is of rank 1 or 2.
Morphisms $\Phi$ such that $f_1$ is of rank 2 are called {\em generic},
and  those such that $f_1$ is of rank 1 are called {\em special}.

\sepsub

\begin{sub}{\sc\small Generic morphisms}\label{fine_gen}\rm

Suppose that $\Phi=(\Phi_1,\Phi_2)$ is a generic morphism.
Let \ $P=\imm(f_1)$ \ and \ $\P_{n-2}\subset\P_n$ be the linear 
subspace of zeroes of linear forms in $P$. Then $\Phi_1$ is isomorphic to
the canonical morphism
\[\ko(-2)\ot P\lra\ko(-1)\]
hence we have \ $\ker(\Phi_1)\simeq\ko(-3)$, and \
 $\imm(\Phi_1)\simeq\ki_{\P_{n-2}}(-1)$ \ (the ideal sheaf of $\P_{n-2}$
 twisted by $\ko(-1)$). Let
\[\Phi' : \ko(-3)\lra\ko\ot\C^k\]
be the restriction of $\Phi_2$ to $\ker(\Phi_1)$. It vanishes on $\P_{n-2}$
and induces a linear map
\[f' : {\C^k}^*\lra H^0(\ki_{\P_{n-2}}(3)) .\]

\sepprop

\begin{subsub}\label{fine_lem1}
{\bf Lemma} : If $\Phi$ is semi-stable (for some polarization)
then $f'$ is injective.
\rm

\begin{proof}
Let \ $K_0=\ker(f')^\bot\subset\C^k$. Then \ 
\m{\imm(\Phi')\subset\ko\ot K_0}. The morphism
\[\ko(-2)\ot\C^2\lra\ko(-1)\oplus(\ko\ot\C^k/K_0)\]
induced by $\Phi$ vanishes on $\ko(-3)=\ker(\Phi_1)$. Hence it induces a
morphism
\[(\psi_1,\psi_2) : \ki_{\P_{n-2}}(-1)
\lra\ko(-1)\oplus(\ko\ot\C^k/K_0)\]
where $\psi_1$ is the inclusion. Since \ \m{\Hom(\ki_{\P_{n-2}}(-1),\ko)=
\Hom(\ko(-1),\ko)}, we can (by replacing $\Phi$ by an element of its
$\Hom(\ko(-1),\ko\ot\C^k)$-orbit) suppose that $\psi_2=0$. It follows that
\m{\imm(\Phi)\subset\ko(-1)\oplus(\ko\ot K_0)}, and since $\Phi$ is
semi-stable, we have \ \m{K_0=\C^k}, i.e. $f'$ is injective. 
\end{proof}
\end{subsub}

\sepprop

Note that we have taken \ $k\leq (n+1)^2=h^0(\ki_{\P_{n-2}}(3))$, to allow the
injectivity of $f'$.

Suppose that $f'$ is injective. Let \ \m{K=\imm(f')}. then
$\Phi'$ is isomorphic to the canonical morphism
\[\phi_K : \ko(-3)\lra\ko\ot K^* .\]
It is easy to see that $P$ and $K$ depend only on the $G$-orbit of $\Phi$. 
Conversely, suppose $P$ and $K$ are given. We can define an element
$(\Phi_1,\Phi_2)$ of $W$ associated to $P$ and $K$ as follows :
let $(z_1,z_2)$ be a basis of
$P$. Let \m{(z_1q_{1i}+z_2q_{2i})_{1\leq i\leq k}} be a basis of $K$, with
\m{q_{1i}, q_{2i}\in S^2V^*}. Using this basis we can identify $K$ and $K^*$
with $\C^k$. We define
\[\Phi_1 : \ko(-2)\ot\C^2\lra\ko(-1)\]
by
\[\xymatrix@R=6pt{
\C^2\ar[r] & V^*\\ (\lambda,\mu)\ar@{|-_{>}}[r] & 
\lambda z_1-\mu z_2 }\]
and
\[\Phi_2 : \ko(-2)\ot\C^2\lra\ko\ot K^*\simeq\ko\ot\C^k\]
over $x\in\P_n$ by
\[\Phi_{2x}(x^2\ot(\lambda,\mu)) \ = \ 
(\lambda q_{2i}(x)+\mu q_{1i}(x))_{1\leq i\leq k} .\]

\sepprop

\begin{subsub}\label{fine_lem2}
{\bf Lemma} : Let $K\subset H^0(\ki_{\P_{n-2}}(3))$ a linear 
subspace of dimension $k$. Then $\Phi_K$ is injective outside of a
closed subvariety of codimension 2. \rm

\begin{proof}
Let $x\in\P_n$. Then $\Phi_K$ is non injective at $x$ if and only if all
the elements of $K$ (which are homogeneous polynomials of degree 3) vanish
at $x$. Suppose that $\Phi_K$ is non injective on an irreducible hypersurface
$S$. Then all the polynomials in $K$ vanish on $S$. Let $f$ be an irreducible
equation of $S$. Then all the elements of $K$ are multiple of $f$. It follows
that $f$ is of positive degree $d\leq 3$, and $K\subset S^{3-d}V^*.f$. But
this in impossible since
\[\dim(K)>\frac{(n+1)(n+2)}{2}=\dim(S^2V^*) .\]
\end{proof}
\end{subsub}

\sepprop

\begin{subsub}\label{fine_lem3}
{\bf Lemma} : Let $\Phi=(\Phi_1,\Phi_2)\in W$ be defined by
$P\subset V^*$ and \m{K\subset H^0(\ki_{\P_{n-2}}(3))}. Suppose that there
exists a polarization such that $\Phi$ is semi-stable. Then $\Phi$ is
generically injective and $\coker(\Phi)$ has no torsion.
Moreover, if $K$ is generic, $\Phi$ is injective.
\rm

\begin{proof}
Lemma \ref{fine_lem2} implies that $\Phi$ is injective outside a closed
subvariety of codimension $\geq 2$. It follows that $\Phi$ is generically
injective and that $\coker(\Phi)$ has no torsion. To prove that $\Phi$ is
injective for a generic $K$, it suffices to find a $K$ such that $\Phi$ is
injective. Let \m{(z_1,z_2)} be a basis of $P$. Let \m{q_1,\cdots,q_r},
(resp. \m{q'_1,\cdots,q'_s}) be linearly independant
elements of $S^2V^*$ that have no
common zeroes in $\P_n$, with \m{r+s=k} (this is possible since \
\m{2n+2\leq k\leq (n+1)^2}). Let
\[K \ = \ \Big(\som_{1\leq i\leq r}z_1q_i\Big)\oplus\Big(
\som_{1\leq j\leq s}z_2q'_j\Big) .\]
It is easy to see that for such a $K$, $\Phi$ is injective.
\end{proof}
\end{subsub}
\end{sub}

\sepsub

\begin{sub}{\sc\small The obvious moduli space of morphisms
and its universal sheaf}\label{fine_obv}\rm

Let $P\subset V^*$ a plane, $\P_{n-2}\subset\P_n$ the subspace defined
by $P$ and \m{K\subset H^0(\ki_{\P_{n-2}}(3))} a linear subspace of
dimension $k$. Let \ \m{\ke(P,K)=\coker(\Phi)}, where $\Phi$ is a
morphism associated to $P$ and $K$. Since the $G$-orbit of $\Phi$ is
determined by $P$ and $K$, $\ke(P,K)$ is well defined. We will give another
construction of $\ke(P,K)$.

Let \ \m{\kf_K=\coker(\Phi_K)}. It is a torsion free sheaf according to
lemma \ref{fine_lem2}. 

\sepprop

\begin{subsub}\label{fine_lem4}
{\bf Lemma} : We have \ $\Ext^1(\ko_{\P_{n-2}}(-1),\kf_K) \simeq \C$, and
the non-trivial extension of $\kf_K$ by $\ko_{\P_{n-2}}(-1)$ is
isomorphic to $\ke(P,K)$.
\rm

\begin{proof}
The exact sequence
\[(*) \ \ \ \ \
0\lra\ko(-3)\lra\ko\ot K^*\lra\kf_K\lra 0\]
implies  \ $H^0(\kf_K(1))\simeq V^*\ot K^*$, $H^1(\kf_K(1))=\nsp$.
Using the exact sequence
\[0\lra\ki_{\P_{n-2}}\lra\ko\lra\ko_{\P_{n-2}}\lra 0\]
we obtain the exact sequence
\[0\to\Hom(\ko(-1),\kf_K)=V^*\ot K^*\lra\Hom(\ki_{\P_{n-2}}(-1),\kf_K)
\lra\Ext^1(\ko_{\P_{n-2}}(-1),\kf_K)\to 0 .\]
From $(*)$ we get \ \m{H^0(\kf_K(2))\simeq S^2V^*\ot K^*},
\m{H^0(\kf_K(3))\simeq S^3V^*\ot K^*/\C i_K}, where $i_K$ is the
inclusion \m{K\in S^3V^*}. From the exact sequence
\[0\lra\ko(-2)\lra\ko(-1)\ot P^{\bot*}\lra \ki_{\P_{n-2}}\lra 0\]
we deduce the exect sequence
\[\xymatrix{
0\ar[r] & \Hom(\ki_{\P_{n-2}}(-1),\kf_K)\ar[r] & S^2V^*\ot P^{\bot}\ot K^*
\ar[rr]^-\theta & & S^3V^*\ot K^*/\C i_K
}\]
where $\theta$ comes from the multiplication
\[\mu : S^2V^*\ot P^\bot\subset S^2V^*\ot V^*\lra S^3V^* .\]
The kernel of $\mu$ is canonically isomorphic to $\wedge^2P^\bot\ot V^*$
and it is easy to see that $i_K$ is contained in the image of 
$\mu\ot I_{K^*}$. It follows that we have an exact sequence 
\[0\lra V^*\ot K^*\lra\ker(\theta)\lra\C i_K\lra 0\]
and that \ $\Ext^1(\ko_{\P_{n-2}}(-1),\kf_K) \simeq \C$.

The last assertion follows from the commutative diagram with exact rows and
columns :
\[\xymatrix{
 & 0\ar[d] & 0\ar[d] & 0\ar[d] \\
0\ar[r] & \ko(-3)\ar[r]\ar[d] & \ko\otimes K^*\ar[r]\ar[d] & \kf_K\ar[r]\ar[d]
& 0 \\
0\ar[r] & \ko(-2)\otimes\C^2\ar[r]\ar[d] & \ko(-1)\oplus(\ko\ot K^*)\ar[r]\ar[d]
& \ke(P,K)\ar[r]\ar[d] & 0 \\
0\ar[r] & \ki_{\P_{n-2}}(-1)\ar[r]\ar[d] & \ko(-1)\ar[r]\ar[d] & 
\ko_{\P_{n-2}}(-1)\ar[r]\ar[d] & 0 \\
& 0 & 0 & 0
}\]
\end{proof}
\end{subsub}

\sepprop

Let ${\bf M}$ be the projective variety of pairs $(P,K)$, where $P$ is a 
plane of $V^*$ and $K\subset H^0(\ki_{\P_{n-2}}(3))$ is a vector subspace of
dimension $k$ ($\P_{n-2}$ beeing the codimension 2 linear subspace of $\P_n$
defined by $P$). We can view ${\bf M}$ as a moduli space for generic 
morphisms. We will give a construction of a {\em universal sheaf} $\E$ on
\ ${\bf M}\times\P_n$, i.e $\E$ is flat on ${\bf M}$ and for every 
$(P,K)\in {\bf M}$, $\E_{(P,K)}$ is
isomorphic to the cokernel of a generic morphism associated to $(P,K)$.
It is also possible to define a {\em universal morphism} whose cokernel is
isomorphic to $\E$, but we will see this more generally in \ref{f_mod}.

Let $Gr(2,V^*)$ be the grassmannian of planes in $V^*$ and
\ \m{q : {\bf M}\to Gr(2,V^*)} \ be the obvious projection. Let $U$ be
the universal subsheaf of $\ko\times V^*$ on $Gr(2,V^*)$. Let
\[p_{\bf M} : {\bf M}\times\P_n\to{\bf M}, \ \ \ \
p_2 : {\bf M}\times\P_n\to\P_n\]
be the projections. 
Then we have a canonical obvious morphism of vector bundles on 
${\bf M}\times\P_n$
\[p_2^*(\ko(-1))\ot U\lra\ko .\]
Let $\kp$ be its cokernel. It is a flat family of sheaves on $\P_n$. 
For every $(P,K)\in {\bf M}$ we have \ \m{\kp_{(P,K)}=\ko_{\P_{n-2}}}.
Let $\kk$ be the universal sheaf on ${\bf M}\times\P_n$, such that \
\m{\kk_{(P,K)} = K}. Then we have a canonical obvious morphism of vector
bundles on ${\bf M}\times\P_n$
\[p_2^*(\ko(-3))\lra \kk^* .\]
Let $\kf$ be its cokernel. Then for every $(P,K)\in{\bf M}$, $\kf_{(P,K)}$
is the sheaf that was noted $\kf_K$ before. By lemma \ref{fine_lem4}, the
sheaf $\EExt^1_{p_{\bf M}}(\kp\ot p_2^*(\ko(-1)),\kf)$ is a line bundle 
$L$ on ${\bf M}$. Then we have a universal extension
\[0\lra\kf\lra\E\lra
\kp\ot p_2^*(\ko(-1))\ot p_{\bf M}^*(L)\lra 0\]
on ${\bf M}\times\P_n$. Then using lemma \ref{fine_lem4} it is easy to see that 
for every $(P,K)\in{\bf M}$, $\E_{(P,K)}$ is
isomorphic to the cokernel of a generic morphism associated to $(P,K)$.

\end{sub}

\sepsub

\begin{sub}{\sc\small Special morphisms}\label{fine_spe}\rm

Let $\Phi=(\Phi_1,\Phi_2)$ be a special morphism. Let \ \m{f_1:\C^2\to V^*}
\ the associated application of rank 1. Let $H$ be the hyperplane of $\P_n$
defined by $\imm(f_1)$. We have an exact sequence
\[\xymatrix{
0\ar[r] & \ko(-2)\ar[r] & \ko(-2)\ot\C^2\ar[r]^-{\Phi_1} & \ko(-1)\ar[r]
& \ko_H(-1)\ar[r] & 0 .
}\]
Let
\[\ov{\Phi_2} : \C^2\lra H^0(\ko_H(2))\ot\C^k\]
be the linear map induced by $\Phi_2$.

\sepprop

\begin{subsub}\label{fine_lem5}
{\bf Lemma} : If $\Phi$ is semi-stable (for a given polarization) then
$\ov{\Phi_2}$ is injective. 
\rm

\begin{proof}
Let $\C_1$, $\C_2$ be the two factors $\C$ of \ \m{\C\oplus\C=\C^2}. We can
suppose thet \ \m{\ker(f_1)=\C_1}. Let
\[\Phi_{2i} : \ko(-2)\ot\C_i\lra\ko\ot\C^k\]
$i=1,2$, be the restrictions of $\Phi_2$, defined by \m{q_{1i},\cdots,
q_{ki}\in S^2V^*}. Let $(z_1,\cdots,z_{n+1})$ be a basis of $V^*$, such that
$z_1$ is an equation of $H$. By using the action of \m{\Hom(\ko(-1)\ot\C^2,
\ko\ot\C^k)} on $W$ we can assume that \m{q_{21},\cdots,q_{2k}\in
S^2<z_2,\cdots,z_{n+1}>}. 

Now $\ov{\Phi_2}$ is not zero on $\C_1$ : otherwise we would have
$q_{1i}\in z_1V^*$, and \ \m{\imm(\Phi_2)\subset\ko\ot\C^{k'}}, with
\[k' \ \leq \ n+1+\dim(S^2<z_2,\cdots,z_{n+1}>) \ \leq \dim(S^2V^*)
\ < \ k ,\]
and this would contradict the semi-stability of $\Phi$. Hence, by
considering the action of $\GL(2)$, it suffices to prove that $\ov{\Phi_2}$
does not vanish on $\C_2$. Suppose it does. Then $\Phi_2$ vanishes on
$\ko(-2)\ot\C_2$, because \m{q_{2i}\in S^2<z_2,\cdots,z_{n+1}>}, and again
\ \m{\imm(\Phi_2)\subset\ko\ot\C^{k'}}, with \
\m{k'\leq\dim(S^2V^*)<k}, which contradicts the semi-stability of $\Phi$.
\end{proof}
\end{subsub}

\sepprop

\begin{subsub}\label{fine_lem6}{\bf Lemma} : 
Suppose that $\Phi$ is semi-stable with respect to some polarization. Then it
is injective outside of a closed subvariety of codimension $\geq 2$, and
$\coker(\Phi)$ has no torsion.
\rm

\begin{proof}
It suffices to prove the first statement. Let $x\in\P_n$ and $u\in\C^2$ 
such that \ \m{\Phi_1(x^2\ot u)=0}. Then we have either $u\in\C_1$ or
$u\not\in\C_1$ and $x\in H$. 
Suppose that $\Phi$ is not injective at all points
of an irreducible hypersurface $D\not =H$. Then the same is true for
$\Phi_{\mid\ko(-2)\ot\C_1}$. Suppose that this morphism is defined by
quadratic forms $q_1,\cdots,q_k$. These forms vanish on $D$, hence they are 
all multiple of an equation of $D$. It follows as in the proof of
\ref{fine_lem5} that \m{\imm(\Phi_2)\subset\ko\ot\C^{k'}}, with $k'<k$,
which contradicts the semi-stability of $\Phi$.

Now it remains to prove that $\Phi_2$ is generically injective on $H$, but
this follows easily from the fact that $\Phi_2$ is defined by an injection
\ \m{\C^2\to H^0(\ko_H(2))\ot\C^k}.
\end{proof}
\end{subsub}
\end{sub}

\sepsub

\begin{sub}{\sc\small Fine moduli spaces of torsion-free sheaves}\label{f_mod}\rm

\begin{subsub} {\bf D\'efinition : }\label{fine_def}\rm
Let $S$ be a smooth variety, $\kf$ a coherent sheaf on $S\times\P_n$, flat
on $S$. We say that $S$ is a {\em fine moduli space of sheaves} with
{\em universal sheaf} $\kf$ if the following properties are verified :
\begin{enumerate}
\item[(i)] For every closed point $s\in S$ the Koda\"\i ra-Spencer map
$$\omega_s : T_sS\lra\Ext^1(\kf_s,\kf_s)$$
is bijective.
\item[(ii)] For every closed points $s_1,s_2\in S$ with $s_1\not= s_2$, 
$\kf_{s_1}$ and $\kf_{s_2}$ are not isomorphic.
\item[(iii)] For every flat family $\ke$ of coherent sheaves on $\P_n$
parametrized by an algebraic variety $T$, and for any closed points $s\in S$,
$t\in T$ such that $\kf_s\simeq\ke_t$, there exists an open neighbourhood
$U$ of $t$ in $T$, and a morphism \ $f:U\to S$ \ such that \ $f(t)=s$ \ 
and
$$(f\times I_{\P_n})^*(\kf) \ \simeq \ \ke_{\mid U} .$$
\end{enumerate}
(cf. \cite{dr5}).

For example moduli spaces of stables sheaves admitting a universal sheaf are
fine moduli spaces of sheaves.
\end{subsub}

\sepsubsub

\begin{subsub}\label{fine_app}Application of theorem \ref{theo_main2}\rm

Polarizations for morphisms
\[\ko(-2)\ot\C^2\lra\ko(-1)\oplus(\ko\ot\C^k)\]
are defined by pairs $(\lambda_1,\lambda_2)$ of positive rational numbers such
that \ \m{\lambda_2+\lambda_1k=1} (so here $\lambda_2$ is associated to 
$\ko(-1)$ and $\lambda_1$ to $\ko\ot\C^k$).

By theorem \ref{theo_main2}, there exists a projective good quotient of
the open subset $W^{ss}$ of semi-stable points as soon as
\[t \ = \ \lambda_2 \ > \ \frac{n+1}{n+1+k} .\]
The critical polarizations in our range are given by
\[\lambda_1=\frac{1}{2p}, \ \ \
t=\lambda_2=1-\frac{k}{2p}, \ \ \
\frac{n+1+k}{2}<p\leq \frac{(n+1)(n+2)}{2} .\]
Let
\[ q \ = \ \frac{(n+1)(n+2)}{2}-\left[\frac{n+1+k}{2}\right]+1\]
(where $[x]$ denotes the integer part of $x$). Then we obtain exactly $q$
moduli spaces of morphisms corresponding to non critical values :
$\M_1,\cdots,\M_q$, where for $1\leq i < q$
\[\M_i \ = \ M(t) \ \ \ \ {\rm for} \ \ \ t=1-\frac{1}{2p}-\epsilon\]
with \ $p=i+\left[\frac{n+1+k}{2}\right]$, $\epsilon$ beeing a sufficiently
small positive rational number. We have \ \m{\M_q={\bf M}} (cf. the end
of \ref{fine_obv}).
\end{subsub}

\sepsubsub

\begin{subsub}\label{fine_mod}Fine moduli spaces\rm

Suppose that we choose a polarization such that $t$ is not a critical value.
In this case we have $W^{ss}=W^s$, and the stabilizer in $G$ of the points of
$W^s$ is the canonical subgroup isomorphic to $\C$. Let \
\m{M(t)=W^s/G}, and \ \m{\pi : W^s\to M(t)} \ be the quotient map. 
On $W^s\times\P_n$ we have a universal morphism
\[\Psi : p_2^*(\ko(-2))\ot\C^2\lra p_2^*(\ko(-1))\oplus(\ko\ot\C^k)\]
(where $p_2$ is the projection \ $W^s\times\P_n\to\P_n$) such that
\ \m{\kf=\coker(\Psi)} \ is a flat family of torsion free sheaves on $\P_n$
parametrized by $W^s$ (this is a consequence of lemmas \ref{fine_lem3} and
\ref{fine_lem6}). There is a canonical action of $G$ on $\kf$ such that
$\C$ acts by multiplication. 

Recall that a $G$-sheaf $\ke$ on $W^s\times\P_n$ {\em descends} to 
$M(t)\times\P_n$
if there exists a coherent sheaf $\ke'$ on $M(t)\times\P_n$ and a
$G$-isomorphism \ \m{(\pi\times I_{\P_n})^*(\ke')\simeq\ke}.
\end{subsub}

\sepprop

\begin{subsub}\label{fine_theo}{\bf Theorem : }
There exists a $G$-line bundle $\kl$ on $M(t)\times\P_n$ such that
$\kf\ot\kl$ descends to $M(t)$. Let $\ke$ be the corresponding sheaf on
$M(t)\times\P_n$. Then $M(t)$ is a fine moduli space of sheaves on $\P_n$
with universal sheaf $\ke$.
\end{subsub}

\begin{proof}
On $W^s$ we have a canonical action of $G$ on the bundles $\ko_{W^s}\ot\C^2$,
$L=\ko_{W^s}$ and $\ko_{W^s}\ot\C^k$. On these bundles $\C$ acts as
ordinary multiplication by scalars. Let $\ka_0$, $\kb_0$ be the $G$-bundles
\[\ka_0 \ = \ (p_2^*(\ko(-2))\ot\C^2)\ot p_W^*(\kl^{-1}), \ \ \ \ 
\kb_0 \ = \ (p_2^*(\ko(-1))\oplus(\ko\ot\C^k))\ot p_W^*(\kl^{-1}) \]
(where $p_W$ is the projection \ $W^s\times\P_n\to W^s$).
On these bundles $\C$ acts trivially. We can multiply the universal
morphism with \m{p_W^*(L^{-1})} and we obtain a new universal morphism
\[\Psi_0 : \ka_0\lra\kb_0 .\]
Now it is easy to see that the bundles $\ka_0$, $\kb_0$ descend to
$M(t)\times\P_n$ either directly from our construction of the quotient,
or by using the more general results of \cite{dr4}, 2.3. Let \
\m{\ka=\ka_0/G}, \m{\kb=\kb_0/G}. The $G$-morphism $\Psi_0$ also descends
and we get a universal morphism of vector bundles on $M(t)\times\P_n$
\[\ov{\Psi} : \ka\lra\kb .\]
We define now \ \m{\ke=\coker(\ov{\Psi})}, and it is clear that
\ \m{\pi^*(\ke)\simeq\kf\ot L^{-1}}. 

Now we prove that the Koda\"\i ra-Spencer map of $\ke$ at $z\in M(t)$ is 
bijective. Let $w\in\pi^{-1}(z)$. Then we have a commutative diagram
\[\xymatrix{
T_wW\ar[r]^-{T\pi}\ar[d]^{\omega_w} & T_zM(t)\ar[d]^{\omega_z} \\
\Ext^1(\kf_w,\kf_w)\fleq[r] & \Ext^1(\ke_z,\ke_z)
}\]
The tangent map $T\pi$ is surjective because $M(t)$ is a geometric quotient.
So it suffices to prove that $\omega_z$ is surjective and that \
\m{\dim(\Ext^1(\ke_z,\ke_z))=\dim(M(t))}.
Consider the exact sequence
\[0\lra\ka_{0w}=\ko(-2)\ot\C^2\lra\kb_{0w}=
\ko(-1)\oplus(\ko\ot\C^k)\lra\kf_w\lra 0 .\]
It is well known that (up to a sign) $\omega_w$ is the composition
\[\Hom(\ka_{0w},\kb_{0w})\lra\Hom(\ka_{0w},\kf_w)\lra
\Ext^1(\kf_w,\kf_w)\]
of maps induced by the preceeding exact sequence. Now the result follows easily
from the exact sequence
\[0\lra\End(\kf_w)\lra\End(\kb_{0w})\lra
\Hom(\ka_{0w},\kb_{0w})/\End(\ka_{0w})\lra\Ext^1(\kf_w,\kf_w)\lra 0 
.\]

We must now verify that if $z_1,z_2\in M(t)$ are distinct closed points, then
$\ke_{z_1}$ and $\ke_{z_2}$ are not isomorphic. This follows from the more
general following result : if two injective morphisms of vector bundles on
$\P_n$
\[\ko(-2)\ot\C^{n_1}\lra(\ko(-1)\ot\C^{m_2})\oplus
(\ko\ot\C^{m_1})\]
have isomorphic cokernels, then they are in the same orbit. 

The property (iii) of the definition of a fine moduli space is easily
verified.
\end{proof}

\sepprop

It follows that the $q$ moduli spaces of morphisms $\M_1,\cdots,\M_q$, with 
their corresponding universal sheaves, are also fine moduli spaces of torsion
free sheaves on $\P_n$. The moduli space $\M_q$ is the same as the 
{\em obvious} one {\bf M} (cf. \ref{fine_obv}), and the corresponding
universal sheaf is the same (up to an element of $Pic({\bf M})$) as $\E$.

These examples are generalizations of the case of $\P_2$ (with $k=7$)
that was treated in \cite{dr5}. But in this case our results are not
needed, because we get only two fine moduli spaces : one is the obvious
moduli space and the other is the corresponding moduli space of stable sheaves
on $\P_2$.

On $\P_n$, $n\geq 3$, our moduli spaces are new. We don't know if the
corresponding moduli space of stable sheaves is among them.

\sepprop

\begin{subsub}{\bf Remark:} \rm it is not hard to prove that all the
moduli spaces $\M_1,\cdots,\M_q$ are distinct.
\end{subsub}
\end{sub}

\vskip 1.5cm

\end{document}